%% file: ModelFree_arXiv.tex
\documentclass[mnsc,nonblindrev]{informs3}

\pdfoutput=1

\OneAndAHalfSpacedXI 

\usepackage{epsfig}
\usepackage{centernot}
\usepackage{url}
\usepackage{multirow}
\usepackage{graphicx}
\usepackage{enumitem}
\usepackage{lscape}
\usepackage{dsfont}
\usepackage[ruled,vlined]{algorithm2e}
\usepackage{algorithmic}
\usepackage{placeins}
\usepackage{mathtools}

\newcommand{\CALL}{\Tc\Ta\Tl\Tl}
\newcommand{\PUT}{\Tp\Tu\Tt}

\input{mathsymbols}

\usepackage{natbib}
 \bibpunct[, ]{(}{)}{,}{a}{}{,}%
 \def\bibfont{\small}%
 \def\BIBand{and}%

\TheoremsNumberedBySection  
\ECRepeatTheorems

\EquationsNumberedBySection 

\MANUSCRIPTNO{MS-OPT-20-01985.R2}

\begin{document}


\RUNAUTHOR{Neufeld, Papapantoleon, Xiang}

\RUNTITLE{Model-free bounds for multi-asset options and their exact computation}

\TITLE{Model-free bounds for multi-asset options using option-implied information and their exact computation}

\ARTICLEAUTHORS{%
\AUTHOR{Ariel Neufeld}
\AFF{Division of Mathematical Sciences, Nanyang Technological University, 21 Nanyang Link, 637371 Singapore, \EMAIL{ariel.neufeld@ntu.edu.sg}} 
\AUTHOR{Antonis Papapantoleon}
\AFF{Delft Institute of Applied Mathematics, TU Delft, 2628 Delft, The Netherlands \& Institute of Applied and Computational Mathematics, FORTH, 70013 Heraklion, Greece, \EMAIL{a.papapantoleon@tudelft.nl}}
\AUTHOR{Qikun Xiang}
\AFF{Division of Mathematical Sciences, Nanyang Technological University, 21 Nanyang Link, 637371 Singapore, \EMAIL{qikun001@e.ntu.edu.sg}}
} 

\ABSTRACT{%
We consider derivatives written on multiple underlyings in a one-period financial market, and we are interested in the computation of model-free upper and lower bounds for their arbitrage-free prices. 
We work in a completely realistic setting, in that we only assume the knowledge of traded prices for other single- and multi-asset derivatives, and even allow for the presence of bid--ask spread in these prices. 
We provide a fundamental theorem of asset pricing for this market model, as well as a superhedging duality result, that allows to transform the abstract maximization problem over probability measures into a more tractable minimization problem over vectors, subject to certain constraints.
Then, we recast this problem into a linear semi-infinite optimization problem, and provide two algorithms for its solution. 
These algorithms provide upper and lower bounds for the prices that are $\varepsilon$-optimal, as well as a characterization of the optimal pricing measures.
These algorithms are efficient and allow the computation of bounds in high-dimensional scenarios (\textit{e.g.} when $d=60$). 
Moreover, these algorithms can be used to detect arbitrage opportunities and identify the corresponding arbitrage strategies.
Numerical experiments using both synthetic and real market data showcase the efficiency of these algorithms, while they also allow to understand the reduction of model risk by including additional information, in the form of known derivative prices.
}%



\maketitle

%


\section{Introduction}
\label{sec:intro}

The classical paradigm in finance and theoretical economics assumes the existence of a model that provides an accurate description of the evolution of asset prices, and all subsequent computations about hedging strategies, exotic derivatives, risk measures, and so forth, are based on this model. 
However, academics, practitioners, and regulators have realized that all models provide only a partially accurate description of this reality, thus, either methods need to be developed in order to aggregate the results of many models, or approaches have to be devised that allow for computations in the absence of a specific model. 
The first approach led to the introduction of robust methods in asset pricing and no-arbitrage theory, see \textit{e.g.}
 \citet{bayraktar2015hedging, beissner2017equilibrium, beissner2019equilibria, bouchard2015arbitrage, bouchard2019superreplication, dana2013intertemporal, epstein2013ambiguous, maruhn2009robust, neufeld2013superreplication, possamai2013robust, rigotti2005uncertainty}, and \citet{yan2018marginal}, while the second one led to model-free methods in asset pricing and no-arbitrage theory, see \textit{e.g.} \citet{acciaio2016model, bartl2019duality}, \citet{bartl2020pathwise, Beiglboeck_HenryLabordere_Penkner_2013, bertsimas2006option, burzoni2016universal}, \citet{burzoni2019pointwise}, \citet{burzoni2021viability, cheridito2017duality, davis2014arbitrage, dolinsky2014martingale, dolinsky2018super, galichon2014stochastic, henry2013automated, hobson1998robust, hu2019data, lutkebohmert2019tightening}, and \citet{riedel2015financial}.

In this work, we consider derivatives written on multiple underlyings in a one-period financial market, and we are interested in the computation of upper and lower bounds for their arbitrage-free prices. 
We work in a completely realistic setting, in that we only assume the knowledge of traded prices for other single- and multi-asset derivatives, and even allow for the presence of bid--ask spread in these prices. 
In other words, we work in a model-free setting in the presence of option-implied information, and make no assumption about the probabilistic evolution of asset prices (\textit{i.e.} their marginal distributions) or their dependence structure. 

The computation of bounds for the prices of multi-asset options, most often basket options, is a classical problem in the mathematical finance literature and has connections with several other branches of mathematics, such as probability theory, optimal transport, operations research, and optimization. 
In the most classical setting, one assumes that the marginal distributions are known and the joint law is unknown; this framework is known as dependence uncertainty. 
In this framework, several authors have derived bounds for multi-asset options using tools from probability theory, such as copulas and Fr\'echet--Hoeffding bounds, see \textit{e.g.} \citet{Chen_Deelstra_Dhaene_Vanmaele_2008, Dhaene_etal_2002_b, Dhaene_etal_2002_a}, and \citet{Hobson_Laurence_Wang_2005_2, Hobson_Laurence_Wang_2005_1}. 
These bounds turned out to be very wide for practical applications, hence recently there was an interest in methods that allow for the inclusion of additional information on the dependence structure, in order to reduce this gap. 
This led to the creation of improved Fr\'echet--Hoeffding bounds and the pricing of multi-asset options in the presence of additional information on the dependence structure, see \textit{e.g.}  \citet{tankov}, \citet{lux2016}, and \citet{puccetti2016}.

The setting of dependence uncertainty is intimately linked with optimal transport theory, and its tools have also been used in order to derive bounds for multi-asset option prices, see \textit{e.g.} \citet{bartl2017marginal} for a formulation in the presence of additional information on the joint distribution. 
More recently, \citet{aquino2019bounds,eckstein2019computation}, and \citet{eckstein2019robust} have translated the model-free superhedging problem into an optimization problem over classes of functions, by extending results in optimal transport, and used neural networks and the stochastic gradient descent algorithm for the computation of the bounds. 

Ideas from operations research and optimization have also been applied for the computation of model-free bounds in settings that are closer to ours, and do not necessarily assume knowledge of the marginal distributions (or, equivalently, knowledge of call option prices for a continuum of strikes). 
\citet{bertsimas2002on} consider the computation of the model-free bounds on a single-asset call option given the moments of the underlying asset price as well as the model-free bounds on a single-asset call option given other single-asset call and put option prices. In addition, they also consider specific conditions under which the model-free bounds on a multi-asset option can be theoretically computed in polynomial time.
\citet{dAspremont_ElGhaoui_2006} consider a framework where the prices of forwards and single-asset call options are known, and compute upper and lower bounds on basket options prices using linear programming. In the more general case where the prices of other basket options are also known, they derive a relaxation to the problem which can be solved using linear programming. 
This work was later extended by various authors.
\citet{Pena_Vera_Zuluaga_2010} improve the results of \cite{dAspremont_ElGhaoui_2006} when computing the lower bounds on basket options prices in two special cases: (i) when the number of assets is limited to two and prices of basket options are known; and (ii) when the prices of only a forward and a single-asset call option per asset are known. 
\citet{pena2012computing} develop a linear programming-based approach for the problem of computing the upper price bound of a basket option given bid and ask prices of vanilla call options.
\citet{pena2010computing} study the problem of computing the upper and lower bounds on basket and spread option prices when the prices of other basket and spread option prices are known. 
Their approach involves solving a large linear programming problem via the Dantzig--Wolfe decomposition in which the corresponding subproblem is solved using mixed-integer programming. 
Compared to \citet{dAspremont_ElGhaoui_2006, pena2010computing}, and \citet{Pena_Vera_Zuluaga_2010, pena2012computing}, the numerical methods we develop in Section~\ref{sec:numerical} apply to settings that are much more general, where the derivative being priced and the traded derivatives with known prices can be any continuous piece-wise affine function (including, but not limited to, vanilla, basket, spread, and rainbow options, as well as any linear combination of these options). 
Moreover, as we demonstrate in Section~\ref{sec:exp}, these methods are able to efficiently compute the price bounds in high-dimensional scenarios, \textit{e.g.} when 60 assets are considered. 
This is considerably higher compared to existing studies. 
\citet{daum2011novel} develop a discretization-based algorithm for solving linear semi-infinite programming problems that returns a feasible solution, and apply the algorithm to compute the upper bounds on basket or spread options prices when single-asset call, put, and exotic options prices are known. 
\citet{Cho_Kim_Lee_2016} develop methods similar to \cite{daum2011novel} but for lower bounds on basket or spread options prices. 
The algorithm we introduce in Section~\ref{ssec:ecp} takes a similar approach, but is able to solve the problem when the prices of multi-asset options with a more general class of payoff functions are known. 
\citet{kahale2017superreplication} uses a central cutting plane algorithm to compute the super- and sub-replicating prices of financial derivatives using hedging portfolios that consist of other financial derivatives in the multi-period discrete-time setting. 
The algorithm only works under the assumption that the underlying state space (\textit{i.e.} the space of asset prices) is finite. 
When the state space is infinite, it is discretized before applying the central cutting plane algorithm, and the discretization error is analyzed. 
However, the approach of discretizing the state space has limited applicability to the multi-dimensional settings (\textit{i.e.} with multiple underlying assets) due to the curse of dimensionality. 
The algorithm we develop in Section~\ref{ssec:accp} is also based on a central cutting plane algorithm, but it allows us to efficiently compute model-free price bounds in high-dimensional state spaces for financial derivatives that depend on multiple assets. 

Our contributions are three-fold:
Firstly, we provide a fundamental theorem of asset pricing for the market model described above, as well as a superhedging duality, that allows to transform the abstract maximization problem over probability measures into a more tractable problem over vectors, subject to certain constraints.
Secondly, we recast this problem into a linear semi-infinite optimization problem, and provide two algorithms for its solution. 
These algorithms provide upper and lower bounds for the prices of multi-asset derivatives that are $\varepsilon$-optimal, as well as a characterization of the optimal pricing measures.
These algorithms are efficient and allow the computation of bounds in high-dimensional scenarios (\textit{e.g.} when $d=60$), that were not possible by previous methods. 
Moreover, these algorithms can be used to detect arbitrage opportunities in multi-asset financial markets and to identify the corresponding arbitrage strategies.
Thirdly, we perform numerical experiments using synthetic data as well as real market data to showcase the efficiency of these algorithms.
These experiments allow us to understand the reduction of the no-arbitrage gap, \textit{i.e.} the difference between the upper and lower no-arbitrage bounds, by including additional information in the form of known derivative prices.
The no-arbitrage gap directly reflects the model-risk associated to a particular derivative and the information available in the market.
The numerical experiments show a decrease of the model-risk by the inclusion of additional information, although this decrease is not uniform and depends on the form of information and the specific structure of the payoff functions.

This paper is organized as follows:
In Section \ref{sec:mainresults}, we present the modeling framework, state the no-arbitrage theorem and the superhedging duality, and discuss a setting that is relevant for practical applications.
In Section \ref{sec:numerical}, we present the algorithms that have been developed for the computation of model-free bounds and state theorems that show the validity of these algorithms.
In Section~\ref{sec:exp},  we discuss various numerical experiments using both synthetic and real market data that show the efficiency of the algorithms and the reduction of model-risk by the inclusion of additional information in the form of known derivative prices.
We also perform a numerical experiment to show the ability of the algorithms to detect arbitrage opportunities.
Appendix~\ref{sec:proofsmain} contains the proofs of the main results of this paper. 
The online appendices contain additional remarks and discussions about the theoretical results,  the numerical methods, and the numerical experiments, as well as the proofs of the results in Section~\ref{sec:mainresults}.

\section{Duality in the presence of option-implied information}
\label{sec:mainresults}

In this section, we introduce a general framework for a model-free, one-period, financial market where multiple assets and several single- and multi-asset derivatives written on these assets are traded simultaneously.
Model-free means that we will not make any assumption about the probabilistic model that governs the evolution of asset prices.
Instead, we will utilize information available in the financial market and implied by the prices of single- and multi-asset derivatives.
We will provide both a fundamental theorem and a superhedging duality in this setting, where our results and proofs are inspired by \citet{bouchard2015arbitrage}.
Moreover, we will describe concrete examples of this framework that are of practical interest.

Throughout this work, all vectors are column vectors unless otherwise stated. 
We denote vectors and vector-valued functions by boldface symbols.
For a vector $\BIx$ in a Euclidean space, let $[\BIx]_j$ denote the $j$-th component of $\BIx$. 
For simplicity, we also use $x_j$ to denote $[\BIx]_j$ when there is no ambiguity.
Let $\|\BIx\|$ denote the Euclidean norm of $\BIx$. 
Let $\langle\BIx,\BIx'\rangle$ denote the Euclidean inner product of two vectors $\BIx$ and $\BIx'$. 
We denote by $\BIe_i$ the $i$-th standard basis vector of a Euclidean space, by $\veczero$ the vector with all entries equal to zero, \textit{i.e.} $\veczero=(0,\ldots,0)^{\mathsf{T}}$, and by $\mathbf 1$ the vector with all entries equal to one, \textit{i.e.} $\mathbf 1=(1,\dots,1)^{\mathsf{T}}$.
We call a subset of a Euclidean space a polyhedron or a polyhedral set if it is the intersection of finitely many closed half-spaces. 
We call a subset of a Euclidean space a polytope if it is a bounded polyhedron. 

Let $\Omega$ be a Polish space equipped with its Borel $\sigma$-algebra denoted by $\CB(\Omega)$. 
Let $\CP(\Omega)$ denote the set of Borel probability measures on $\Omega$. 
Let $g_j:\Omega\to\R$ be Borel measurable for $j=1,\ldots,m$, for some fixed $m\in\N$, and let $\BIg:\Omega\to\R^m$ denote the vector-valued Borel measurable function where the $j$-th component corresponds to $g_j$. 
Let $\underline{\pi}_j,\overline{\pi}_j\in\R$ be such that $\underline{\pi}_j\le\overline{\pi}_j$ for $j=1,\ldots,m$. 
Let $\BIy=(y_1,\ldots,y_m)^{\mathsf T}\in\R^m$ and define $\pi:\R^m\to\R$ as
\begin{align}
\pi(\BIy):=\sum_{j=1}^m y^+_j\overline{\pi}_j-y^-_j\underline{\pi}_j,
\label{eqn:pricepi}
\end{align} 
where $y^+_j:=\max\{y_j,0\}, y^-_j:=\max\{-y_j,0\}$. 
Let $\langle \BIy,\BIg\rangle$ denote the function $\sum_{j=1}^my_jg_j:\Omega\to\R$. 

We make the following \textit{no-arbitrage} assumption.
\begin{assumption}[No-Arbitrage]
The following implication holds for any $\BIy\in\R^m$:
\begin{align*}
\langle\BIy,\BIg\rangle-\pi(\BIy)\ge 0 \quad\Longrightarrow\quad \langle\BIy,\BIg\rangle-\pi(\BIy)=0,
\end{align*}
where the inequality and the equality are both understood as point-wise. 
\label{asp:na}
\end{assumption}

\begin{remark}
Assumption~\ref{asp:na} is inspired by the no-arbitrage assumption introduced in Definition~1.1 of \citet{bouchard2015arbitrage}, where the set of possible models for the market is $\CP(\Omega)$, \textit{i.e.} all Borel probability measures, and a single time step is considered. 
The difference between Assumption~\ref{asp:na} and the no-arbitrage assumption in \cite{bouchard2015arbitrage} is that the price of a financial derivative in the present work is not a singleton but can lie anywhere between the corresponding bid and ask prices. 
Note that there are other notions of no-arbitrage that are weaker than Assumption~\ref{asp:na}, for example the ``no uniform strong arbitrage'' assumption in Definition~2.1 of \citet{bartl2017marginal}. 
\end{remark}

Let $f:\Omega\to\R$ be a Borel measurable function, and define the functional $\phi(f)$ as follows:
\begin{align}
	\phi(f) := \inf\big\{c+\pi(\BIy):c\in\R,\BIy\in\R^m, c+\langle\BIy,\BIg\rangle\ge f\big\}.
\label{eqn:superreplicate}
\end{align}
Let $\CQ$ be defined as follows:
\begin{align*}
\CQ:=\left\{\mu\in\CP(\Omega):\underline{\pi}_j\le\int_{\Omega} g_j\Td\mu\le\overline{\pi}_j,\text{ for }j=1,\ldots,m\right\}.
\end{align*}

The main results of this section are the following fundamental theorem and superhedging duality, whose proofs are provided in Appendix~\ref{sec:proofduality}. 
\begin{theorem}[Fundamental Theorem]
\label{thm:equivalence}
The following are equivalent:
\begin{enumerate}
	\item[(i)] Assumption~\ref{asp:na} holds.
	\item[(ii)] For all $\nu\in\CP(\Omega)$, there exists $\mu\in\CQ$ such that $\nu\ll\mu$.  
\end{enumerate}
\end{theorem}



\begin{theorem}[Superhedging duality]
\label{thm:duality}
Under Assumption~\ref{asp:na}, the following statements hold.
\begin{enumerate}
\item[(i)] $\phi(f)>-\infty$.
\item[(ii)] There exists $\BIy\in\R^m$ such that $\phi(f)+\langle\BIy,\BIg\rangle-\pi(\BIy)\ge f$. Hence, the infimum in (\ref{eqn:superreplicate}) is attained when $\phi(f)<\infty$. 
\item[(iii)] We have the following superhedging duality result:
\begin{align}
\phi(f)=\sup_{\mu\in\CQ}\int_{\Omega} f\Td\mu.
\label{eqn:duality}
\end{align}
\end{enumerate}
\end{theorem}


\begin{remark}
No-arbitrage conditions and Fundamental Theorems of Asset Pricing are essential tools to understand and characterize the viability of a model in a financial market.
They have to be tailored to the modeling assumptions and the specific applications in mind, hence a multitude of comparable statements exist in the literature. 
Analogously to our no-arbitrage condition, the Fundamental Theorem presented in Theorem \ref{thm:equivalence} is closely related to the First Fundamental Theorem in \citet{bouchard2015arbitrage}.
The main difference is the presence of a bid--ask spread, which means that we cannot exactly reduce our results to their theorem and another proof is needed; this proof is motivated by the results in \cite{bouchard2015arbitrage}.
There are several other versions of a Fundamental Theorem in the presence of model uncertainty in the mathematical finance literature; see \textit{e.g.} \citet{acciaio2016model,Bayrakter_Zhang_2016}, \citet{Bayrakter_Zhang_Zhou_2014}, and \citet{burzoni2021viability} for discrete time models, and \citet{Biagini_Bouchard_Kardaras_Nutz_2017} and \citet{Dolinsky_Soner_2014_b} for continuous time ones. 

Let us point out that the Fundamental Theorem presented here plays a particular role in conjunction with the numerical methods developed in the next section.
More specifically, it provides a sufficient condition for the detection of arbitrage opportunities by numerically testing the violation of the no-arbitrage condition. 
Moreover, it provides a sufficient condition for repairing derivative prices by removing arbitrage opportunities from the market, in the same spirit as \citet{cohen2020detecting}. 
These results are novel in the related literature on multi-asset model-free price bounds, and are facilitated by the tailor-made Fundamental Theorem.
\end{remark}

\begin{remark}
Superhedging dualities are also classical and essential tools in mathematical finance, typically tailored to specific modeling assumptions and applications.
The superhedging duality presented in Theorem \ref{thm:duality} is motivated by the superhedging theorem in \citet{bouchard2015arbitrage}, with the main difference being once again that we are considering an interval of bid and ask prices instead of a single price.
There are multiple comparable duality results or superhedging theorems in various areas of mathematics.
In the mathematical finance literature, these results are known as superhedging theorems or (martingale) optimal transport dualities, see \textit{e.g.} \citet{Beiglboeck_HenryLabordere_Penkner_2013, acciaio2016model, Bayrakter_Zhang_Zhou_2014, cheridito2017duality}, and \citet{dolinsky2014martingale}. 
In the operations research literature, these results are known as perfect or strong dualities, see \textit{e.g.} \citet{bertsimas2002on}, \citet{dAspremont_ElGhaoui_2006}, \citet{nishihara2007duality}, and \citet{Pena_Vera_Zuluaga_2010}.
These latter dualities are typically based on classical results in mathematical programming, see \textit{e.g.} \citet{Karlin_Studden_1966}, and \citet{Hettich_Kortanek_1993}.

Let us point out that the superhedging duality (\ref{eqn:duality}) is crucial when verifying the $\varepsilon$-optimality of a measure in the numerical algorithms introduced in Section~\ref{sec:numerical}, see Theorem~\ref{thm:accpalg} and Corollary~\ref{coro:ecpprimal}. 
\end{remark}


The canonical way to interpret the framework developed above is as follows: when $\Omega=\R^d_+$, then there exist $d$ underlying risky assets that are traded in the financial market, and $\Omega$ represents the (non-negative) prices of the assets at a fixed future date. 
Investing into a unit of the asset $i$ then corresponds to the payoff function $g(\BIx)\equiv\mathrm{proj}_i(\BIx):=x_i$ for $\BIx\in\R^d_+$.
Moreover, there exist $m$ traded derivatives (typically $m\gg d$) with known bid and ask prices $(\underline\pi_j, \overline\pi_j)_{j=1:m}$, written either on single or on multiple assets. 
The payoff function $g_j$ of a single-asset derivative depends on the price of only a single asset, \textit{i.e.} $g_j=\widetilde{g}_j\circ\mathrm{proj}_i$ for some $i\in\{1,\ldots,m\}$ and $\widetilde{g}_j:\R_+\to\R$. For example, $g_j(\BIx)=(x_i-\kappa_j)^+$ corresponds to a call option with strike $\kappa_j$. The payoff function $g_j$ of a multi-asset derivative depends on the prices of multiple assets. For example, $g_j(\BIx) = (\langle\BIw,\BIx \rangle-\kappa_j)^+$ corresponds to a basket call option with weight $\BIw$ and strike $\kappa_j$. These derivatives encode all the information available in this market. Specifically, information about the marginals of the probability measures in $\CQ$ is \textit{implied} by the bid and ask prices of single-asset derivatives, while \textit{partial} information on the joint distribution is \textit{implied} by the bid and ask prices of multi-asset derivatives.

In this setting, the right-hand side of \eqref{eqn:duality} is the model-free upper bound for the price of a derivative with payoff function $f$ written on these $d$ assets.
The optimization takes place over all probability measures that are compatible with the option-implied information, \textit{i.e.} all probability measures that produce prices for a given option within its respective bid and ask prices.
The duality result in \eqref{eqn:duality} states that this model-free upper bound equals the least superhedging price achieved by trading in the $m$ derivatives according to the strategy $(c,\BIy)$, \textit{i.e.} holding $c$ units of cash and $y_j$ units of derivative $j$ for $j=1,\ldots,m$, where the minimization takes place over all $(c,\BIy)$ such that the payoff $f$ is dominated, \textit{i.e.} $c+\langle\BIy,\BIg\rangle \ge f$. 


Section~\ref{ecsec:duality} in the online appendices contains additional discussions about the duality result, including Example~\ref{ex:primalnotattained} which demonstrates that the supremum on the right-hand side of (\ref{eqn:duality}) is not necessarily attained, as well as Proposition~\ref{prop:examplena} which provides a specific setting in which Assumption~\ref{asp:na} holds.







\section{Numerical methods for the computation of bounds}
\label{sec:numerical}

The superhedging duality in Theorem~\ref{thm:duality} allows to transform the abstract maximization problem over probability measures into a more tractable minimization problem over vectors that satisfy certain constraints.
The aim of this section is to develop novel numerical methods for the exact and efficient computation of upper and lower bounds on $\phi(f)$.
More specifically, we will develop methods for the computation of upper and lower bounds $\phi(f)^{\TU\TB}$ and $\phi(f)^{\TL\TB}$ which are $\varepsilon$-optimal, \textit{i.e.}
\[
	\phi(f)^{\TL\TB} \le \phi(f) \le \phi(f)^{\TU\TB} \ \text{ and } \ \phi(f)^{\TU\TB} - \phi(f)^{\TL\TB} \le \varepsilon,
\]
for $\varepsilon>0$.
Our methods allow us to also characterize the optimal pricing measure associated with the primal maximization problem.
Therefore, we provide a complete solution to both optimization problems, and can characterize the solution both in terms of $\varepsilon$-optimal hedging strategies and in terms of the optimal pricing measure.

Let $\overline{\Bpi}:=(\overline{\pi}_1,\ldots,\overline{\pi}_m)^{\mathsf T}$ and $\underline{\Bpi}:=(\underline{\pi}_1,\ldots,\underline{\pi}_m)^{\mathsf T}$. 
The minimization problem $\phi(f)$ in \eqref{eqn:superreplicate} can be equivalently formulated as a linear semi-infinite programming (LSIP) problem, \textit{i.e.} as an optimization problem with a linear objective and an infinite number of linear constraints, one for each $\omega\in\Omega$,
\begin{align}
	\begin{split}
		\phi(f) = \quad \text{minimize }\quad& c+\langle\BIy^+,\overline{\Bpi}\rangle-\langle\BIy^-,\underline{\Bpi}\rangle\\
				\text{subject to }\quad&c+\langle\BIy^+-\BIy^-,\BIg(\omega)\rangle\ge f(\omega)\quad \forall \omega\in\Omega,\\
				&c\in\R,\;\BIy^+\ge\veczero,\;\BIy^-\ge\veczero.
	\end{split}
\label{eqn:lsipdef}
\end{align}
LSIP problems are classical optimization problems that have been thoroughly studied in the related literature; see \textit{e.g.} \citet{goberna1998linear, goberna2018recent}. 
More general semi-infinite programming problems, including non-linear semi-infinite programming problems and \textit{generalized} semi-infinite programming problems (where the index set can depend on the decision variable) have also been studied in the literature; see \textit{e.g.} \citet{Reemtsen_Rueckmann_1998, lopez2007semi}, and \citet{stein2012how}.
In this section, we develop novel algorithms tailored to solving~(\ref{eqn:lsipdef}) under different assumptions on the space $\Omega$, and the functions $\BIg$ and $f$. 

%

Let us first introduce the notion of CPWA functions and their radial functions. 
\begin{definition}[Continuous piece-wise affine function and its radial function]
We call a function $h:\R^d\to\R$ a CPWA function if it can be represented as
\begin{align}
	h(\BIx) &= \sum_{k=1}^{K}\xi_k\max\big\{\langle\BIa_{k,i},\BIx\rangle+b_{k,i}:1\le i\le I_k\big\},
\label{eqn:cpwadef}
\end{align}
where $K\in\N$, $I_k\in\N$ for $k=1,\ldots,K$, and $\BIa_{k,i}\in\R^d$, $b_{k,i}\in\R$, $\xi_k\in\{-1,1\}$ for $i=1,\ldots,I_k,k=1,\ldots,K$. The \emph{radial function} of $h$, denoted by $\widetilde{h}:\R^d_+\to\R$, is defined as
\begin{align*}
\widetilde{h}(\BIz):=\sum_{k=1}^{K}\xi_k\max\left\{\langle\BIa_{k,i},\BIz\rangle:1\le i\le I_k\right\}.
\end{align*}
\label{def:cpwa} 
\end{definition}
The class of CPWA functions contains many popular payoff functions in finance, including vanilla call and put options, basket options, spread options, call/put-on-max options, call/put-on-min options, best-of-call options, etc. 
We refer the reader to Section~\ref{ecssec:cpwa} in the online appendices for the CPWA representations of these payoff functions as well as some properties of CPWA functions.

\subsection{CPWA payoff functions on unbounded domains}
\label{ssec:ecp}

In the first setting, we work under the following assumptions. 
\begin{assumption}[Setting 1]
We assume the following:
\begin{enumerate}
\item[(i)] $\Omega=\R^d_+$;
\item[(ii)] $f$ and $(g_j)_{j=1:m}$ are CPWA functions on $\Omega$; 
\item[(iii)] $\phi(f)<\infty$ and $\phi(-f)<\infty$. 
\end{enumerate}
\label{asp:setting1}
\end{assumption}

In the sequel, for notational reasons, we use $\BIx$ in place of $\omega$ when $\Omega$ is a subset of the Euclidean space. 
Let us introduce the notion of the slack function for the LSIP problem in (\ref{eqn:lsipdef}). 
\begin{definition}[Slack function]\label{def:slackfunc}
Let $\BIy\in\R^m$ be fixed, and denote the slack function of the LSIP problem in (\ref{eqn:lsipdef}) by $s_{\BIy}:\R^d_+\to\R$, which is defined as: $s_{\BIy}(\BIx):=\langle\BIy,\BIg(\BIx)\rangle-f(\BIx)$. 
\end{definition}




\begin{algorithm}[t]
\KwIn{$\overline{\Bpi}$, $\underline{\Bpi}$, $(g_j)_{j=1:m}$, $f$, $X^{(0)}\subset\R^d_+$, $\underline{\phi}$, $\overline{\BIx}>\veczero$, $\varepsilon>0$, $\tau>0$, $0<\delta\le1$}
\KwOut{$\phi(f)^{\TU\TB}$, $\phi(f)^{\TL\TB}$, $c^\star$, $\BIy^\star$, $X$}
\nl Formulate the function $\widetilde{s}_{\BIy}(\cdot)$ and use Algorithm~\ref{alg:radcons} (see Section~\ref{ecssec:ecp}) to generate radial constraints, denoted by a system of linear inequalities with auxiliary variables $\widetilde{\sigma}$. \label{alglin:ecpinitaux}\\
\nl Construct a system of linear inequalities $\sigma$ which contains all variables and inequalities in $\widetilde{\sigma}$, additional variables $c\in\R$, $\BIy^+\ge\veczero$ and $\BIy^-\ge\veczero$, and an additional equality $\BIy=\BIy^+-\BIy^-$. \label{alglin:ecpinit2}\\
\nl Add the linear inequality $c+\langle\BIy^+,\overline{\Bpi}\rangle-\langle\BIy^-,\underline{\Bpi}\rangle\ge \underline{\phi}-\tau$ to $\sigma$. \label{alglin:ecpinitlb} \\
\nl $r\leftarrow 0$. \\
\nl \Repeat{${s}^{(r-1)}\ge-\varepsilon$}{
\nl \label{alglin:ecploopstart}\For{each $\BIx\in X^{(r)}$}{
\nl Add the linear inequality $c+\langle\BIy^+-\BIy^-,\BIg(\BIx)\rangle\ge f(\BIx)$ to $\sigma$. \label{alglin:ecpfeascut}\\ 
}
\nl Solve the LP problem: $\underline{\varphi}^{(r)}\leftarrow$ minimize $c+\langle\BIy^+,\overline{\Bpi}\rangle-\langle\BIy^-,\underline{\Bpi}\rangle$ subject to linear constraints $\sigma$, and denote the computed minimizer as $(c^{(r)},\BIy^{+(r)},\BIy^{-(r)})$. \label{alglin:ecplp}\\
\nl $\BIy^{(r)}\leftarrow\BIy^{+(r)}-\BIy^{-(r)}$. \\
\nl Formulate the global minimization problem: $\min_{\veczero\le\BIx\le\overline{\BIx}}s_{\BIy^{(r)}}(\BIx)$ into a mixed-integer linear programming (MILP) problem (see (\ref{eqn:milp})). \label{alglin:ecpformulatemilp}\\
\nl ${s}^{(r)}\leftarrow c^{(r)}+\min_{\veczero\le\BIx\le\overline{\BIx}}s_{\BIy^{(r)}}(\BIx)$ (solve the MILP problem via the BnB algorithm). \label{alglin:ecpmilp}\\
\nl $X^{(r+1)}\leftarrow\Big\{\BIx:(\BIx,(\lambda_k),(\zeta_k),(\delta_{k,i}),(\iota_{k,i}))$ is an integer feasible solution found by the BnB algorithm while solving the MILP problem such that $c^{(r)}+s_{\BIy^{(r)}}(\BIx)\le\delta s^{(r)}\Big\}$. \label{alglin:ecpintfeas}\\
\nl $r\leftarrow r+1$. \label{alglin:ecploopend}\\
}
\nl $\phi(f)^{\TL\TB}\leftarrow\underline{\varphi}^{(r-1)}$, $\phi(f)^{\TU\TB}\leftarrow\underline{\varphi}^{(r-1)}-{s}^{(r-1)}$, $c^\star\leftarrow c^{(r-1)}-{s}^{(r-1)}$, $\BIy^\star\leftarrow\BIy^{(r-1)}$, $X\leftarrow\bigcup_{l=0}^{r-1}X^{(l)}$. \label{alglin:ecpfeas}\\
\nl \If{$\phi(f)^{\TU\TB}<\underline{\phi}$}{
\nl \Return the problem (\ref{eqn:lsipdef}) is unbounded. \label{alglin:ecpunbounded}\\
}
\nl \Else{
\nl \Return $\phi(f)^{\TU\TB}$, $\phi(f)^{\TL\TB}$, $c^\star$, $\BIy^\star$, $X$.\\
}
    \caption{{\bf Exterior Cutting Plane (ECP) Algorithm} (under Assumption~\ref{asp:setting1})}
    \label{alg:ecpalgo}
\end{algorithm}

To numerically solve the LSIP problem (\ref{eqn:lsipdef}), let us now introduce the cutting plane discretization method, detailed in Algorithm~\ref{alg:ecpalgo}, which is inspired by the Conceptual Algorithm~11.4.1 in \citet{goberna1998linear}. 
In Line~\ref{alglin:ecpinitaux} of Algorithm~\ref{alg:ecpalgo}, the so-called ``radial constraints'' are generated using Algorithm~\ref{alg:radcons} in Section~\ref{ecssec:ecp} of the online appendices. The purpose of this step is to generate sufficient and necessary constraints on $\BIy$ to guarantee that the slack function $s_{\BIy}(\cdot)$ is bounded from below. 
We refer the reader to Section~\ref{ecssec:ecp} for a detailed explanation about this step. 
In Line~\ref{alglin:ecpformulatemilp} and \ref{alglin:ecpmilp}, the global minimization problem $\min_{\veczero\le\BIx\le\overline{\BIx}}s_{\BIy^{(r)}}(\BIx)$ is solved by formulating it into a mixed-integer linear programming (MILP) problem in (\ref{eqn:milp}), as discussed in Lemma~\ref{lem:milp}. 
The MILP problem can be solved efficiently by state-of-the-art solvers such as Gurobi (\cite{gurobi}), that uses the so-called Branch-and-Bound (BnB) algorithm. We refer the reader to Remark~\ref{rmk:bnb} for a brief description of the BnB algorithm. 


We name Algorithm~\ref{alg:ecpalgo} the exterior cutting plane (ECP) method, since every constraint (also known as cut) generated in Line~\ref{alglin:ecpfeascut} does not restrict the feasible set of (\ref{eqn:lsipdef}), hence is exterior to the feasible set. 
We refer the reader to Section~\ref{ecssec:ecp} for detailed discussions about various aspects of Algorithm~\ref{alg:ecpalgo}. Specifically, Remark~\ref{rmk:ecpalg} explains the inputs of Algorithm~\ref{alg:ecpalgo}, and Remark~\ref{rmk:ecpalgdiff} discusses the differences between Algorithm~\ref{alg:ecpalgo} and the Conceptual Algorithm~11.4.1 in \citet{goberna1998linear}. 
Under the assumption that the inputs of Algorithm~\ref{alg:ecpalgo} are specified according to Remark~\ref{rmk:ecpalg}, Theorem~\ref{thm:ecpalg} below shows the properties of Algorithm~\ref{alg:ecpalgo}, whose proof is provided in Appendix~\ref{ssec:proofsetting1main}.


%

\begin{theorem}[Properties of Algorithm~\ref{alg:ecpalgo}]
Let Assumption~\ref{asp:setting1} hold. 
Assume that $\overline{\BIx}$ and $\underline{\phi}$ are specified as stated in Remark~\ref{rmk:ecpalg}. 
Then, the following statements hold.
\begin{enumerate}
\item[(i)] If Assumption~\ref{asp:na} holds, then $\underline{\varphi}^{(r)}$ is non-decreasing in $r$. 
	At any stage of Algorithm~\ref{alg:ecpalgo}, $s^{(r)}\le0$ and $\underline{\varphi}^{(r)}\le\phi(f)\le\underline{\varphi}^{(r)}-{s}^{(r)}$.
\item[(ii)] If Assumption~\ref{asp:na} holds, then Algorithm~\ref{alg:ecpalgo} terminates after finitely many iterations with an $\varepsilon$-optimal solution $(c^\star,\BIy^\star)$ of (\ref{eqn:superreplicate}) and $\phi(f)^{\TL\TB}\le\phi(f)\le\phi(f)^{\TU\TB}$ with $\phi(f)^{\TU\TB}-\phi(f)^{\TL\TB}\le\varepsilon$. 
\item[(iii)] If Line~\ref{alglin:ecpunbounded} of Algorithm~\ref{alg:ecpalgo} is reached, then Assumption~\ref{asp:na} is violated and problem~(\ref{eqn:lsipdef}) is unbounded. 
\end{enumerate}
\label{thm:ecpalg}
\end{theorem}

\begin{remark}
In Section~\ref{ssec:accp}, under the more restrictive assumption that $\Omega=\{\BIx\in\R^d:\veczero\le\BIx\le\overline{\BIx}\}$ for some $\overline{\BIx}>\veczero$\footnote{$\BIx=(x_1,\ldots,x_d)^{\mathsf T}>\veczero  \quad\Leftrightarrow\quad  x_1>0,\ldots,x_d>0.$} (see Assumption~\ref{asp:setting2}), we show that Algorithm~\ref{alg:ecpalgo} also produces an $\varepsilon$-optimal solution to the right-hand side of (\ref{eqn:duality}), which corresponds to the most extreme pricing measure in the original model-free superhedging problem. 
This will be explained in detail in Corollary~\ref{coro:ecpprimal}. 
\end{remark}


\subsection{CPWA payoff functions on bounded domains}
\label{ssec:accp}

In the second setting, we adopt similar but more restrictive assumptions than in the first one. 

\begin{assumption}[Setting 2]
We assume the following:
\begin{enumerate}
\item[(i)] $\Omega=\{\BIx\in\R^d:\veczero\le\BIx\le\overline{\BIx}\}$ for $\overline{\BIx}=(\overline{x}_1,\ldots,\overline{x}_d)^{\mathsf T}>\veczero$;
\item[(ii)] $f$ and $(g_j)_{j=1:m}$ are CPWA functions on $\Omega$; 
\item[(iii)] $\phi(f)<\infty$ and $\phi(-f)<\infty$. 
\end{enumerate}
\label{asp:setting2}
\end{assumption}


\SetKwRepeat{RepeatNoUntil}{while $\overline{\varphi}^{(r)}-\underline{\varphi}^{(r)}>\varepsilon$}{\normalfont \textit{(continues on the next page)}}%
\SetKwFor{RepeatNoRepeat}{\normalfont \textit{(continued)}}{}%

\begin{algorithm} [p]
\KwIn{$\overline{\Bpi},\underline{\Bpi}$, $(g_j)_{j=1:m}$, $f$, $X^{(0)}\subset\Omega$, $\overline{\BIx}$, $\underline{\phi}$, $\overline{\phi}$, $c^{\star(0)}$, $\BIy^{\star(0)}$, $\overline{c}>0$, $\overline{\BIy}>\veczero$, $\varepsilon>0$, $\tau>\varepsilon$, $0\le\gamma<1$, $0<\zeta<1$, $0<\delta\le1$}
\KwOut{$\phi(f)^{\TU\TB}$, $\phi(f)^{\TL\TB}$, $c^\star$, $\BIy^\star$, $c^\dagger$, $\BIy^\dagger$, $X^\dagger$, $X$}
\nl $r\leftarrow 0$, $\underline{\varphi}^{(0)}\leftarrow\underline{\phi}-\tau$, $\overline{\varphi}^{(0)}\leftarrow\overline{\phi}$, flag $\leftarrow$ false. \\
\nl Mark all elements of $X^{(0)}$ as active and removable. \\
\nl \label{alglin:accploop}\RepeatNoUntil{}{
\nl $r\leftarrow r+1$, $\underline{\varphi}^{(r)}\leftarrow\underline{\varphi}^{(r-1)}$, $\overline{\varphi}^{(r)}\leftarrow\overline{\varphi}^{(r-1)}$, $\varphi^{(r)}\leftarrow\frac{\underline{\varphi}^{(r)}+\overline{\varphi}^{(r)}}{2}$, $c^{\star(r)}\leftarrow c^{\star(r-1)}$, $\BIy^{\star(r)}\leftarrow\BIy^{\star(r-1)}$. \label{alglin:accpbisect1}\\
\nl \If{flag is true}{
\nl $\varphi^{(r)}\leftarrow(\underline{\varphi}^{(r)}+{\varphi}^{(r)})/2$. \label{alglin:accpbisect2}\\
}
\nl $X\leftarrow\bigcup_{0\le l\le r-1}\{\BIx\in X^{(l)}:\BIx \text{ is marked active}\}$. \\
\nl Compute the Chebyshev center of $\sigma(\overline{c},\overline{\BIy},\underline{\varphi}^{(r)},{\varphi}^{(r)},X)$ by solving an LP problem. \label{alglin:accpcclp}\\
\nl \If{the LP problem in Line~\ref{alglin:accpcclp} is infeasible}{
\nl $\underline{\varphi}^{(r)}\leftarrow\min\{c+\langle\BIy^+,\overline{\Bpi}\rangle-\langle\BIy^-,\underline{\Bpi}\rangle:(c,\BIy^+,\BIy^-)\text{ satisfies }\sigma(\overline{c},\overline{\BIy},-\infty,\infty,X)\}$ (which is an LP problem). Let $(c^\dagger,\BIy^{+\dagger},\BIy^{-\dagger})$ be a minimizer of this LP problem. \label{alglin:accplblp}\\
\nl $\BIy^{\dagger}\leftarrow\BIy^{+\dagger}-\BIy^{-\dagger}$, $X^\dagger\leftarrow X$. \label{alglin:accplblp2}\\
\nl Mark all elements of $\bigcup_{1\le l\le r-1}X^{(r)}$ as removable.\\
\nl $\rho^{(r)}\leftarrow-1$, $X^{(r)}\leftarrow\emptyset$. \\
\nl Skip to the next iteration. \\
}
\nl Let $(c^{(r)},{\BIy}^{+(r)},{\BIy}^{-(r)})$ be the Chebyshev center and let $\rho^{(r)}$ be the radius of the largest inscribed ball of $\sigma(\overline{c},\overline{\BIy},\underline{\varphi}^{(r)},{\varphi}^{(r)},X)$. \label{alglin:accpccenter}\\
\nl $\BIy^{(r)}\leftarrow{\BIy}^{+(r)}-{\BIy}^{-(r)}$. \label{alglin:accpmodifyweights}\\
\nl Formulate the global minimization problem: minimize $c^{(r)}+s_{\BIy^{(r)}}(\BIx)$ subject to $\veczero\le\BIx\le\overline{\BIx}$ into an MILP problem (see (\ref{eqn:milp})). Solve it with relative gap tolerance $\zeta$. Let $\overline{s}^{(r)}$ be its approximate optimal value. Let $\underline{s}^{(r)}$ be its lower bound at termination. \label{alglin:accpmilp}\\
\nl $X^{(r)}\leftarrow\Big\{\BIx:(\BIx,(\lambda_k),(\zeta_k),(\delta_{k,i}),(\iota_{k,i}))$ is an integer feasible solution found by the BnB algorithm while solving (\ref{eqn:milp}) such that $c^{(r)}+s_{\BIy^{(r)}}(\BIx)\le\delta \overline{s}^{(r)}\Big\}$. Mark all elements of $X^{(r)}$ as active and removable.  \label{alglin:accpintfeas}\\
\nl \label{alglin:accpobjcutcheck}\If{$c^{(r)}+\pi(\BIy^{(r)})-\underline{s}^{(r)}<\overline{\varphi}^{(r)}-\varepsilon$}{
\nl $\overline{\varphi}^{(r)}\leftarrow c^{(r)}+\pi(\BIy^{(r)})-\underline{s}^{(r)}$, $c^{\star(r)}\leftarrow c^{(r)}-\underline{s}^{(r)}$, $\BIy^{\star(r)}\leftarrow\BIy^{(r)}$. \label{alglin:accpobjcut}\\
\nl \If{$\underline{s}^{(r)}\ge0$}{
\nl Mark all elements of $\bigcup_{1\le l\le r}X^{(l)}$ as removable. \label{alglin:accpfeasible}\\
\nl Skip to the next iteration (Line~\ref{alglin:accploop}). \\
}
}
}
    \caption{{\bf Accelerated Central Cutting Plane (ACCP) Algorithm} (under Assumption~\ref{asp:setting2})}
    \label{alg:accpalgo}
\end{algorithm}


\setcounter{algocf}{1}
\begin{algorithm}[t!]
\setcounter{AlgoLine}{23}
\RepeatNoRepeat{}{
\nl \If{flag is true}{
\nl flag $\leftarrow$ false. \\
\nl Mark all elements of $X^{(r)}$ as non-removable, and skip to the next iteration (Line~\ref{alglin:accploop}). \\
}
\nl flag $\leftarrow$ true. \\
\nl \For{each $0\le l\le r$ such that $\rho^{(r)}<\gamma\rho^{(l)}$}{
\nl \For{each $\BIx\in X^{(l)}$ marked as removable}{
\nl \If{$c^{(r)}+\langle\BIy^{(r)},\BIg(\BIx)\rangle-\big(1+\|\BIg(\BIx)\|_2^2\big)^{\frac{1}{2}}\rho^{(r)}>f(\BIx)$}{
\nl Set $\BIx$ as inactive. \label{alglin:accpremove}\\
}
}
}
}
\nl $\phi(f)^{\TU\TB}\leftarrow\overline{\varphi}^{(r)}$, $\phi(f)^{\TL\TB}\leftarrow\underline{\varphi}^{(r)}$, $c^\star\leftarrow c^{\star(r)}$, $\BIy^\star\leftarrow\BIy^{\star(r)}$, $X\leftarrow\bigcup_{0\le l\le r}\{\BIx\in X^{(l)}:\BIx \text{ is marked active}\}$. \\
\nl \If{$\phi(f)^{\TU\TB}<\underline{\phi}$}{
\nl \Return the problem (\ref{eqn:lsipdef}) is unbounded. \label{alglin:accpunbounded}\\
}
\nl \Else{
\nl \Return $\phi(f)^{\TU\TB}$, $\phi(f)^{\TL\TB}$, $c^\star$, $\BIy^\star$, $c^\dagger$, $\BIy^\dagger$, $X^\dagger$, $X$.\\
}
    \caption{{\bf Accelerated Central Cutting Plane (ACCP) Algorithm} (under Assumption~\ref{asp:setting2})}
    \label{alg:accpalgo2}
\end{algorithm}


Let us introduce a version of the accelerated central cutting plane (ACCP) method inspired by \citet{betro2004accelerated}, detailed in Algorithm~\ref{alg:accpalgo}. 
In Algorithm~\ref{alg:accpalgo}, we maintain and update a sequence of lower bounds $(\underline{\varphi}^{(r)})_{r\ge0}$ of $\phi(f)$, a sequence of upper bounds $(\overline{\varphi}^{(r)})_{r\ge0}$ of $\phi(f)$, and polytopes in $\R^{2m+1}$ that are denoted by $\sigma(\overline{c},\overline{\BIy},\underline{\varphi}^{(r)},{\varphi}^{(r)},X)$, which have the form
\begin{align}
\begin{split}
\sigma(\overline{c},\overline{\BIy},\underline{\varphi}^{(r)},{\varphi}^{(r)},X):=\Big\{(c,\BIy^+,\BIy^-):\;&|c|\le\overline{c},\;\veczero\le\BIy^+\le\overline{\BIy},\;\veczero\le\BIy^-\le\overline{\BIy},\\
&\underline{\varphi}^{(r)}\le c+\langle\BIy^+,\overline{\Bpi}\rangle-\langle\BIy^-,\underline{\Bpi}\rangle\le{\varphi}^{(r)},\\
&c+\langle\BIy^+-\BIy^-,\BIg(\BIx)\rangle\ge f(\BIx),\;\forall\BIx\in X\Big\},
\end{split}
\label{eqn:accppolytope}
\end{align}
where $\overline{c}$, $\overline{\BIy}$ specify a bounding box, $\underline{\varphi}^{(r)}$ and $\overline{\varphi}^{(r)}$ specify the lower and upper bounds on $c+\pi(\BIy^+-\BIy^-)$, and $X\subset\Omega$ specifies a collection of feasibility constraints. 
In Algorithm~\ref{alg:accpalgo}, ${\varphi}^{(r)}$ is between the lower bound $\underline{\varphi}^{(r)}$ and the upper bound $\overline{\varphi}^{(r)}$, and is used as a speculative upper objective cut. 
The idea of Algorithm~\ref{alg:accpalgo} is that $(\underline{\varphi}^{(r)})_{r\ge0}$ is a non-decreasing sequence of lower bounds that approach $\phi(f)$ from below while $(\overline{\varphi}^{(r)})_{r\ge0}$ is a non-increasing sequence of upper bounds that approach $\phi(f)$ from above. 
These facts will be made clear later in Theorem~\ref{thm:accpalg}. 
Algorithm~\ref{alg:accpalgo} has various advantages over Algorithm~\ref{alg:ecpalgo}. 
Most importantly, the MILP problem in Line~\ref{alglin:accpmilp} of Algorithm~\ref{alg:accpalgo} only needs to be solved approximately with a large error tolerance, and some linear constraints are removed in Line~\ref{alglin:accpremove} to make solving the LP problem in Line~\ref{alglin:accpcclp} faster. 
We refer the reader to \citet{betro2004accelerated} and Section~11.4 of \citet{goberna1998linear} for further discussions.

A crucial step of Algorithm~\ref{alg:accpalgo} is to compute the Chebyshev center, that is, the center of the largest inscribed ball, of the polytope $\sigma(\overline{c},\overline{\BIy},\underline{\varphi}^{(r)},{\varphi}^{(r)},X)$ in Line~\ref{alglin:accpcclp}. 
It is well-known that the Chebyshev center of a polytope can be computed by solving an LP problem (see, for example, the Conceptual Algorithm~11.4.2 of \citet{goberna1998linear}).

%

We refer the reader to Section~\ref{ecssec:accp} for detailed discussions about various aspects of Algorithm~\ref{alg:accpalgo}. Specifically, Remark~\ref{rmk:accpalg} explains its inputs, and Remark~\ref{rmk:accpalgdiff} discusses the differences between Algorithm~\ref{alg:accpalgo} and the ACCP algorithm by \citet{betro2004accelerated}. 
Assuming that the inputs of Algorithm~\ref{alg:accpalgo} are specified according to Remark~\ref{rmk:accpalg}, Theorem~\ref{thm:accpalg} below shows the properties of Algorithm~\ref{alg:accpalgo}, whose proof is provided in Appendix~\ref{ssec:proofsetting2}. 


%

\begin{theorem}[Properties of Algorithm~\ref{alg:accpalgo}]
Let Assumption~\ref{asp:setting2} hold. 
Assume that $\underline{\phi}$, $\overline{\phi}$, $c^{\star(0)}$, $\BIy^{\star(0)}$, $\overline{c}$, $\overline{\BIy}$ are specified as stated in Remark~\ref{rmk:accpalg}. 
Then, the following statements hold.
\begin{enumerate}
\item[(i)] If Assumption~\ref{asp:na} holds, then $\underline{\varphi}^{(r)}$ is non-decreasing in $r$, $\overline{\varphi}^{(r)}$ is non-increasing in $r$. 
	Moreover, at any stage of Algorithm~\ref{alg:accpalgo}, $\underline{\varphi}^{(r)}\le\phi(f)\le\overline{\varphi}^{(r)}$, and $c^{\star(r)}+\langle\BIy^{\star(r)},\BIg\rangle\ge f$ holds. 
\item[(ii)] If Assumption~\ref{asp:na} holds, then Algorithm~\ref{alg:accpalgo} terminates after finitely many iterations with an $\varepsilon$-optimal solution $(c^\star,\BIy^\star)$ of (\ref{eqn:superreplicate}) and $\phi(f)^{\TL\TB}\le\phi(f)\le\phi(f)^{\TU\TB}$ with $\phi(f)^{\TU\TB}-\phi(f)^{\TL\TB}{\le\varepsilon}$.
\item[(iii)] If Assumption~\ref{asp:na} holds, then $c^\dagger$, $\BIy^\dagger$ and $X^\dagger$ are defined when Algorithm~\ref{alg:accpalgo} terminates. 
	If $|c^\dagger|<\overline{c}$ and $-\overline{\BIy}<\BIy^\dagger<\overline{\BIy}$, then the following LP problem with decision variables $(\mu_\BIx)_{\BIx\in X^\dagger}$: 
\begin{align}
\begin{split}
\text{maximize }\quad& \sum_{\BIx\in X^\dagger}\mu_{\BIx}f(\BIx)\\
\text{subject to }\quad&\sum_{\BIx\in X^\dagger}\mu_{\BIx}=1,\quad\underline{\Bpi}\le\sum_{\BIx\in X^\dagger}\mu_{\BIx}\BIg(\BIx)\le\overline{\Bpi},\quad\mu_{\BIx}\ge0\;\;\forall\BIx\in X^\dagger
\end{split}
\label{eqn:accpdual}
\end{align}
has an optimal solution $(\mu^\star_{\BIx})_{\BIx\in X^\dagger}$. 
	Let $\mu^\star$ be a finitely supported measure defined by $\mu^\star:=\sum_{\BIx\in X^\dagger}\mu^\star_{\BIx}\delta_{\BIx}$. 
	Then, $\mu^\star$ is $\varepsilon$-optimal for the right-hand side of (\ref{eqn:duality}).
\item[(iv)] If Line~\ref{alglin:accpunbounded} of Algorithm~\ref{alg:accpalgo} is reached, then Assumption~\ref{asp:na} is violated and problem~(\ref{eqn:lsipdef}) is unbounded. 
\end{enumerate}
\label{thm:accpalg}
\end{theorem}

Theorem~\ref{thm:accpalg}(iii) explicitly provides a pricing measure which is an $\varepsilon$-optimal solution to the original model-free superhedging problem. 
The ECP method (Algorithm~\ref{alg:ecpalgo}) is also applicable in Setting~2 with Line~\ref{alglin:ecpinitaux} removed. It has the same property as Theorem~\ref{thm:accpalg}(iii), which is detailed in the next corollary and proved in Appendix~\ref{ssec:proofsetting2}. 

\begin{corollary}
Under Assumption~\ref{asp:na} and Assumption~\ref{asp:setting2}, Algorithm~\ref{alg:ecpalgo} (with Line~\ref{alglin:ecpinitaux} removed) terminates after finitely many iterations, and the following LP problem with decision variables $(\mu_\BIx)_{\BIx\in X}$: 
\begin{align}
\begin{split}
\text{maximize }\quad& \sum_{\BIx\in X}\mu_{\BIx}f(\BIx)\\
\text{subject to }\quad&\sum_{\BIx\in X}\mu_{\BIx}=1,\quad\underline{\Bpi}\le\sum_{\BIx\in X}\mu_{\BIx}\BIg(\BIx)\le\overline{\Bpi},\quad\mu_{\BIx}\ge0\;\;\forall\BIx\in X
\end{split}
\label{eqn:ecpprimallp}
\end{align}
has an optimal solution $(\mu^\star_{\BIx})_{\BIx\in X}$. 
Define the finitely supported measure $\mu^\star$ by $\mu^\star:=\sum_{\BIx\in X}\mu^\star_{\BIx}\delta_{\BIx}$. Then, $\mu^\star$ is $\varepsilon$-optimal for the right-hand side of (\ref{eqn:duality}). 
\label{coro:ecpprimal}
\end{corollary}

\begin{remark}
Under Assumptions~\ref{asp:setting1} and \ref{asp:setting2}, the payoff functions of the traded derivatives and the target derivative $f$ must be continuous piece-wise affine (CPWA) functions. 
Thus, the proposed Algorithms~\ref{alg:ecpalgo} and \ref{alg:accpalgo} are unable to directly treat derivatives with non-CPWA payoff functions, such as digital and power options. 
Hence, it would be necessary to first approximate these payoff functions by CPWA functions in order to treat such derivatives using the proposed algorithms. 
Existing studies have proposed methods to treat non-CPWA payoff functions under more restrictive assumptions. 
\citet{bertsimas2002on} develop a method to compute the price bounds on a single vanilla call option given moments of the underlying asset price. 
Moreover, in the multi-asset setting, they show that the theoretical time complexity to compute the price bounds is polynomial when all payoff functions are the sum of a CPWA function and a quadratic function, under the assumption that all of the CPWA functions share the same partition of $\R^d$ and that the number of polyhedra in the partition is polynomial in the number of assets and the number of traded derivatives. 
However, this assumption is rather restrictive, since the presence of a fixed number of traded vanilla options written on each asset incurs a partition of $\R^d$ in which the number of polyhedra is exponential in the number of assets. 
In \citet{daum2011novel}, a method is developed to compute the price bounds on a basket call option given the prices of vanilla call and put options and the prices of single-asset digital and power options. 
In this setting, due to the structure of the basket call option, the global optimization problem associated with the LSIP problem can be reduced to a sequence of one-dimensional global optimization problems which can then be efficiently solved. 

On the computational side, our assumption that all derivatives have CPWA payoff functions makes it possible to formulate the global minimization problem associated with the LSIP problem (\ref{eqn:lsipdef}) into an MILP problem, which then allows us to efficiently solve the LSIP problem in high-dimensional situations (\textit{e.g.} when $d=60$) using state-of-the-art solvers. 
On the practical side, the exclusion of derivatives with non-CPWA payoff functions does not harm the applicability of the methods we have developed, since most relevant financial derivatives have CPWA payoff functions; see Example~\ref{ex:CPWA-payoffs}. 
Let us point out, that many trading platforms of digital options are un-regulated, some adopt ethically questionable practices, while the involved risk is anyhow difficult to manage and hedge.
Because of these, the trading of digital options has been banned in many countries including Australia\footnote{\url{https://asic.gov.au/about-asic/news-centre/find-a-media-release/2021-releases/21-064mr-asic-bans-the-sale-of-binary-options-to-retail-clients/}, accessed: April 18, 2021.}, Canada\footnote{\url{https://www.investmentexecutive.com/news/from-the-regulators/binary-options-ban-takes-effect/}, accessed: April 18, 2021.}, Israel\footnote{\url{https://www.reuters.com/article/us-israel-binaryoptions-lawmaking/israel-ban-on-binary-options-gets-final-parliamentary-approval-idUSKBN1CS2L1}, accessed: April 18, 2021.}, and all EU countries\footnote{\url{https://www.esma.europa.eu/press-news/esma-news/esma-agrees-prohibit-binary-options-and-restrict-cfds-protect-retail-investors}, accessed: April 18, 2021.}.
Power options, on the other hand, are conceptual financial derivatives that are not traded in real markets. 
\end{remark}



\begin{remark}
Algorithms~\ref{alg:ecpalgo} and \ref{alg:accpalgo} both solve a sequence of constrained optimization problems which are relaxations of the LSIP problem (\ref{eqn:lsipdef}) where the semi-infinite constraint is reduced to finitely many constraints. This approach is known as the discretization method. 
Instead of enforcing these constraints strictly, there is an alternative discretization method in which the constraints are replaced by penalty functions, see \textit{e.g.} \citet{coope1998exact} and \citet{auslender2009penalty}. 
In this method, a sequence of unconstrained optimization problems are solved. 
We refer the reader to \citet{auslender2015comparative} for a comparison of different variants of this method. 
Another approach that is based on penalization is the integral-type penalization method, see \textit{e.g.} \citet{borwein1991convergence, borwein1991duality, lin1998unconstrained}, and \citet{auslender2009penalty}. 
This approach transforms the original LSIP problem into an unconstrained convex optimization problem with an integral-type penalty term, which can be subsequently solved by the stochastic (sub)gradient descent (SGD) algorithm. 
This is similar to a recent approach adopted for solving optimal transport, martingale optimal transport, and related problems, see \textit{e.g.} \citet{eckstein2018robust, eckstein2019computation, eckstein2019robust}, and \citet{aquino2019bounds}. 
We have empirically tested the \textit{integral-type penalization plus SGD} approach in our problem, and found that it is unstable and highly sensitive to the initialization and the hyper-parameter settings of the SGD algorithm. 

\label{rmk:penalty}
\end{remark}


\section{Numerical experiments and results}
\label{sec:exp}


In this section, we perform three experiments using synthetically generated derivative prices and one experiment using real market data, in order (i) to demonstrate the performance of the proposed approaches under these settings, (ii) to quantify the effect of the additional information from traded multi-asset options on the width of the no-arbitrage gap, \textit{i.e.} the difference between the upper and lower model-free bounds, and (iii) to show that the proposed algorithms are capable of detecting arbitrage opportunities in a financial market.



We refer the reader to Section~\ref{ecssec:expdetails} in the online appendices for details about the implementation of the proposed numerical algorithms. The code used in this work is available on GitHub\footnote{\url{https://github.com/qikunxiang/ModelFreePriceBounds}}. 
In the subsequent numerical experiments, we consider financial derivatives with the following CPWA payoff functions,
whose representations are discussed in Example \ref{ex:CPWA-payoffs}.

\begin{enumerate}[label=(\roman*)]
\item Trading in the $i$-th asset: $g(\BIx)=x_i$.
\item Vanilla call option on the $i$-th asset with strike $\kappa>0$: $g(\BIx)=\big(x_i-\kappa\big)^+$.
\item Basket call option with weights $\BIw\in\R^d_+$ and strike $\kappa>0$: $g(\BIx)=\big(\sum_i w_ix_i-\kappa\big)^+$.
\item Spread call option with weights $\BIw\in\R^d\setminus \R^d_+$ (\textit{e.g.} $\BIw=\BIe_i-\BIe_j$) and strike $\kappa\in\R$, \textit{e.g.}:~$g(\BIx)=\big(x_i - x_j - \kappa\big)^+$.
\item Call-on-max (rainbow) option on assets $i_1,\ldots,i_l$ with strike $\kappa\ge0$: $g(\BIx) = \big[(x_{i_1} \vee x_{i_2} \vee \dots \vee x_{i_l})-\kappa\big]^+$.
\item Call-on-min (rainbow) option of assets $i_1,\ldots,i_l$ with strike $\kappa\ge0$: $g(\BIx) = \big[(x_{i_1} \wedge x_{i_2} \wedge \dots \wedge x_{i_l})-\kappa\big]^+$.
\item Best-of-calls option of assets $i_1,\ldots,i_l$ with strikes $\kappa_1,\ldots,\kappa_l\ge0$: $g(\BIx) = (x_{i_1}-\kappa_1)^+ \vee (x_{i_2}-\kappa_2)^+ \vee \dots \vee (x_{i_l}-\kappa_{l})^+$.
\end{enumerate}

Moreover, in the numerical experiments with synthetically generated prices, we consider market models of the following type:
\begin{itemize}
\item The marginal distribution of the price of an asset at terminal time is a log-normal distribution.
	Under Setting~2, the $i$-th marginal distribution is truncated to $[0,\overline{x}_i]$ for $i=1,\ldots,d$.
\item The dependence structure among the marginals of the $d$ assets at terminal time is a $t$-copula with a positive definite correlation matrix $\BC$ and $\nu$ degrees of freedom. 
\end{itemize}

Given these market models, the prices of the single-asset derivatives listed above can be computed in closed-form by taking the discounted expectations of the corresponding payoff functions (w.r.t. a pricing measure). 
We have assumed that the interest rate is equal to zero, for the sake of simplicity.
For the multi-asset derivatives listed above, we approximate their prices via Monte Carlo integration by randomly generating one million independent samples from the copula model and subsequently using these samples to approximate the expectations of the payoff functions. 
The markets models are also used to compute reference prices for the target derivatives (with payoff $f$). However, they are not used in the computation of the model-free bounds. 
In the computation of the bounds, the only information used are the prices of single- and multi-asset derivatives that are synthetically generated from these market models.
In order to simulate an incomplete market with the presence of bid--ask spread, we specify multiple market models with different parameters and subsequently take the minimum (resp.\ maximum) price of a derivative among its prices under all models as the bid (resp.\ ask) price of the derivative. 

Notice that under the market models described above, the pricing measure $\widehat{\mu}$ has strictly positive density with respect to the Lebesgue measure on $\Omega$.
Moreover, the way the prices of derivatives are generated guarantees that $\widehat{\mu}\in\CQ$. Therefore, Assumption~\ref{asp:na} holds by Proposition~\ref{prop:examplena}.



\subsection{Experiment~1}
\label{ssec:exp1}

In this experiment, we consider a financial market with 5 assets ($d=5$). 
We consider Setting~2 (\textit{i.e.} Assumption~\ref{asp:setting2}) where $\Omega=[0,100]^5$.
Our goal is to compute the model-free lower and upper price bounds for a call-on-max option with payoff function $f(\BIx) = \big( x_{2} \vee x_{3} \vee x_{4} - \kappa \big)^+$, where the strike price $\kappa$ ranges from 0 to 10 with an increment of $0.2$. 
We assume that a total of 439 financial derivatives are traded in the market ($m=439$). 
These include:
\begin{itemize}
\item The 5 assets $x_1,x_2,x_3,x_4,x_5$. 
\item Vanilla call options on the 5 assets with strikes $1,2,\ldots,10$. 
\item Basket call options with the following weights and strikes $1,2,\ldots,10$:
\newline
	$(\frac{1}{5},\frac{1}{5},\frac{1}{5},\frac{1}{5},\frac{1}{5})^{\mathsf T}$, 
	$(\frac{1}{4},\frac{1}{4},\frac{1}{4},\frac{1}{4},0)^{\mathsf T}$, 
			$(\frac{1}{4},\frac{1}{4},\frac{1}{4},0,\frac{1}{4})^{\mathsf T}$, 
			$(\frac{1}{4},\frac{1}{4},0,\frac{1}{4},\frac{1}{4})^{\mathsf T}$, 
			$(\frac{1}{4},0,\frac{1}{4},\frac{1}{4},\frac{1}{4})^{\mathsf T}$, 
			$(0,\frac{1}{4},\frac{1}{4},\frac{1}{4},\frac{1}{4})^{\mathsf T}$,
	$(\frac{1}{3},\frac{1}{3},\frac{1}{3},0,0)^{\mathsf T}$, 
			$(0,\frac{1}{3},\frac{1}{3},\frac{1}{3},0)^{\mathsf T}$, 
			$(0,0,\frac{1}{3},\frac{1}{3},\frac{1}{3})^{\mathsf T}$, 
	$(0,\frac{1}{2},\frac{1}{2},0,0)^{\mathsf T}$, 
			$(0,\frac{1}{2},0,\frac{1}{2},0)^{\mathsf T}$, 
			$(0,0,\frac{1}{2},\frac{1}{2},0)^{\mathsf T}$.	
 \item Spread call options with the following weights and strikes $-5,-4,\ldots,0,1,\ldots,5$:\newline
	$\BIe_1-\BIe_2$,
	$\BIe_1-\BIe_3$, 
	$\BIe_1-\BIe_4$, 
	$\BIe_2-\BIe_3$, 
	$\BIe_2-\BIe_4$, 
	$\BIe_2-\BIe_5$, 
	$\BIe_3-\BIe_4$, 
	$\BIe_3-\BIe_5$, 
	$\BIe_4-\BIe_5$, 
	$\BIe_2-\BIe_1$, 
	$\BIe_3-\BIe_1$, 
	$\BIe_4-\BIe_1$, 
	$\BIe_3-\BIe_2$, 
	$\BIe_4-\BIe_2$, 
	$\BIe_5-\BIe_2$, 
	$\BIe_4-\BIe_3$, 
	$\BIe_5-\BIe_3$, 
	$\BIe_5-\BIe_4$. 
\item Call-on-max (rainbow) options on the following 6 groups of assets and strikes $0,1,\ldots,10$:
	$\{x_1,x_2,x_3,x_4,x_5\}$, 
	$\{x_1,x_2,x_3,x_4\}$, 
	$\{x_2,x_3,x_4,x_5\}$, 
	$\{x_2,x_3\}$, 
	$\{x_2,x_4\}$, 
	$\{x_3,x_4\}$. 
\end{itemize} 
The bid and ask prices of the assets and derivatives are synthetically generated from the market models specified in Section~\ref{ecssec:exp1}. 
We consider five cases, where we use certain subsets of the traded derivatives to compute the model-free lower and upper price bounds for the target derivative:
\begin{itemize}
\item Case 1 (denoted as \textit{V}): we use only vanilla options;
\item Case 2 (\textit{V+B}): we use vanilla and basket options;
\item Case 3 (\textit{V+B+S}): we use vanilla, basket, and spread options;
\item Case 4 (\textit{V+B+S+R}): we use vanilla, basket, spread and call-on-max (rainbow) options;
\item Case 5 (\textit{V+R}): we use vanilla and call-on-max options.
\end{itemize}

%

\begin{figure}[t]
\centering
\includegraphics[width=0.42\linewidth]{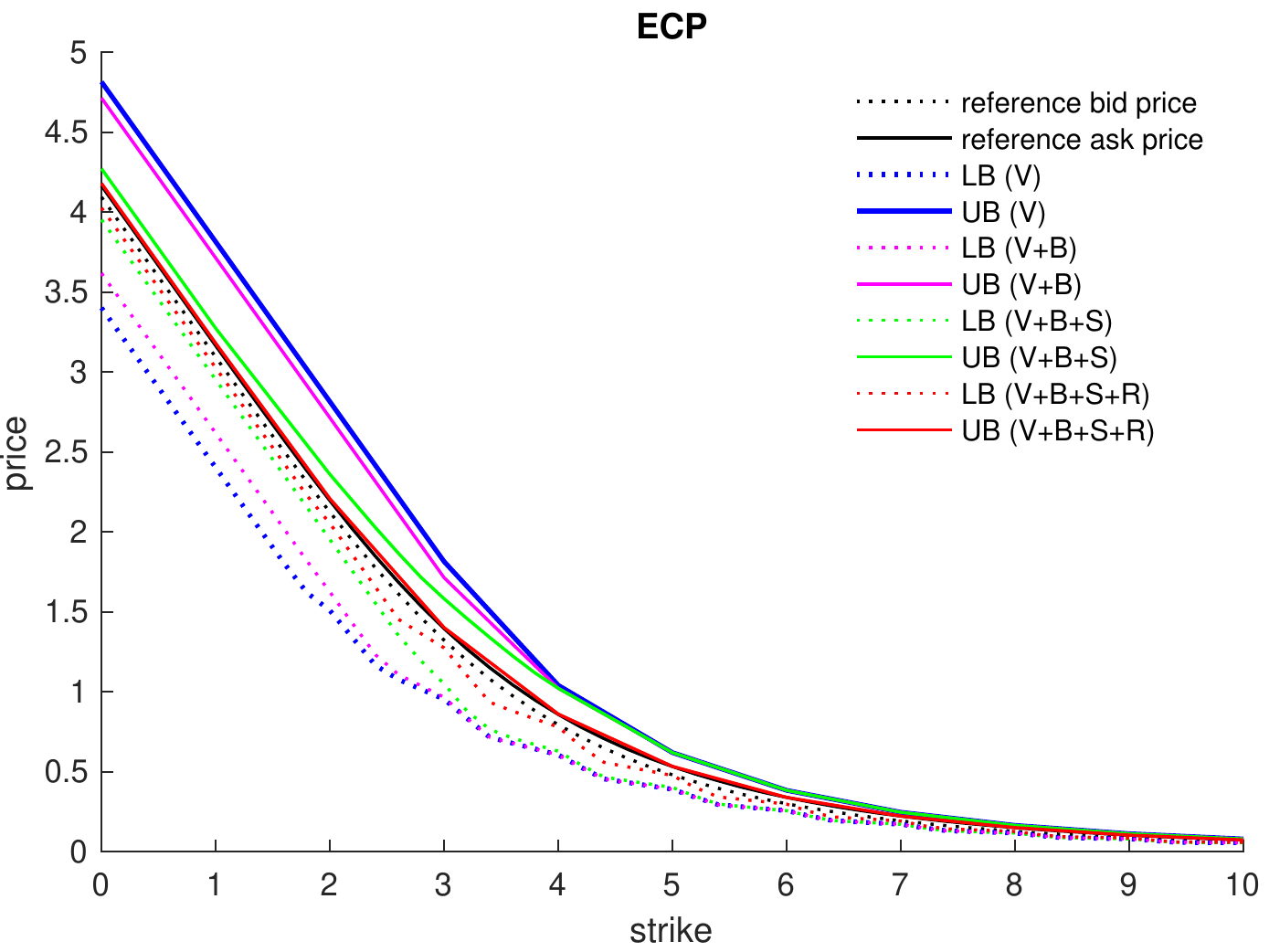}
~
\includegraphics[width=0.42\linewidth]{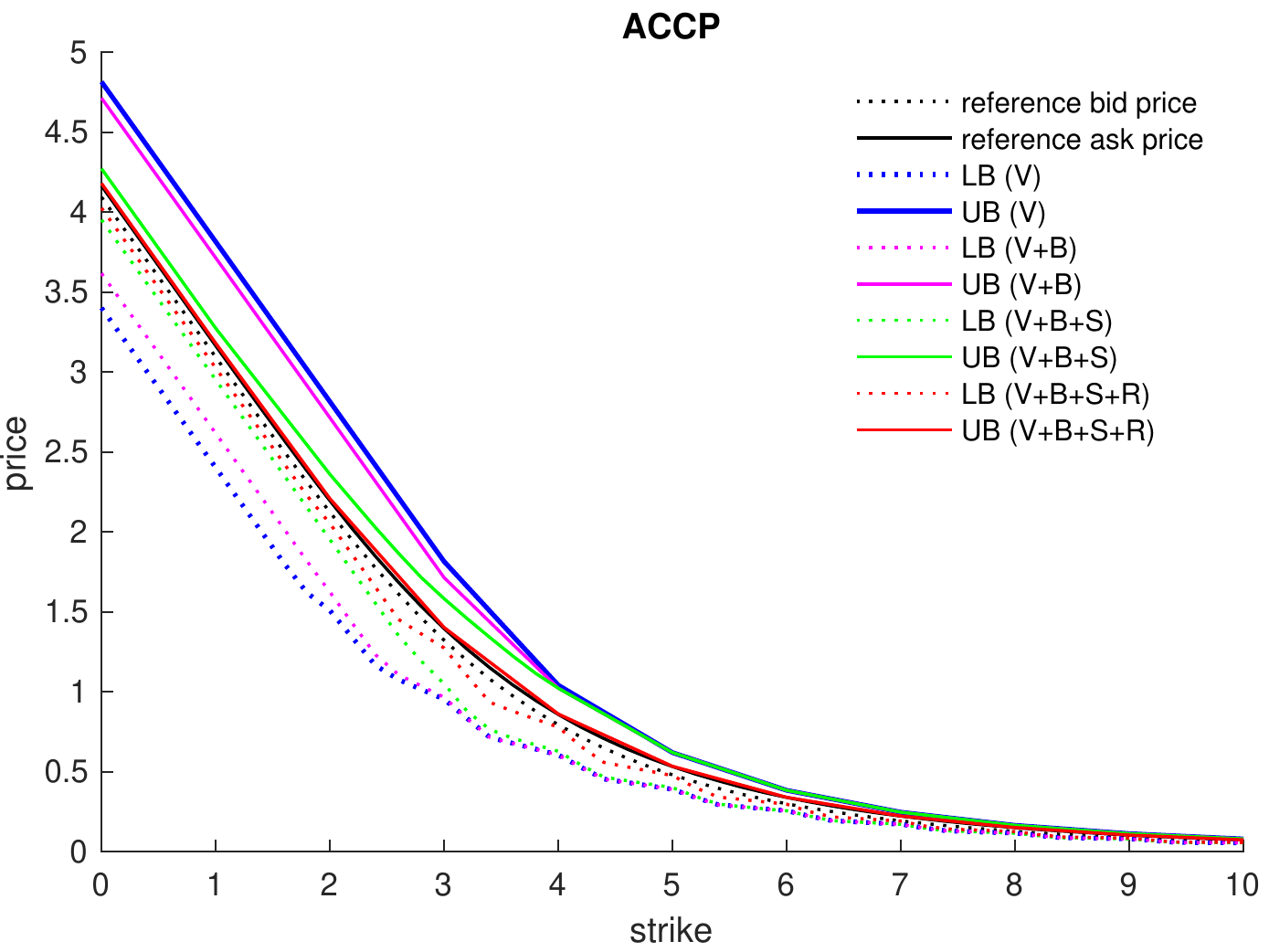}

\vspace{0.2cm}

\includegraphics[width=0.42\linewidth]{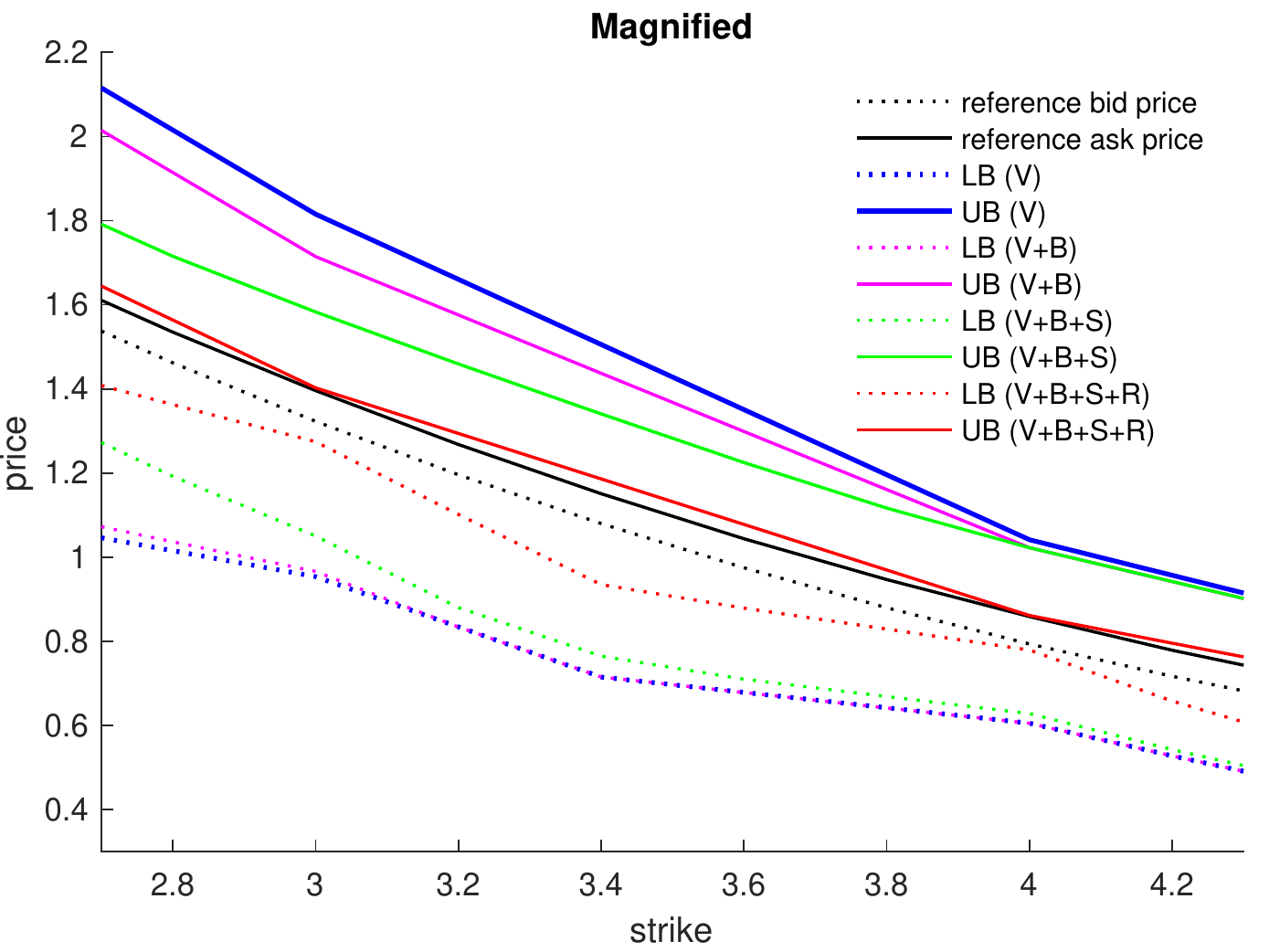}
~
\includegraphics[width=0.42\linewidth]{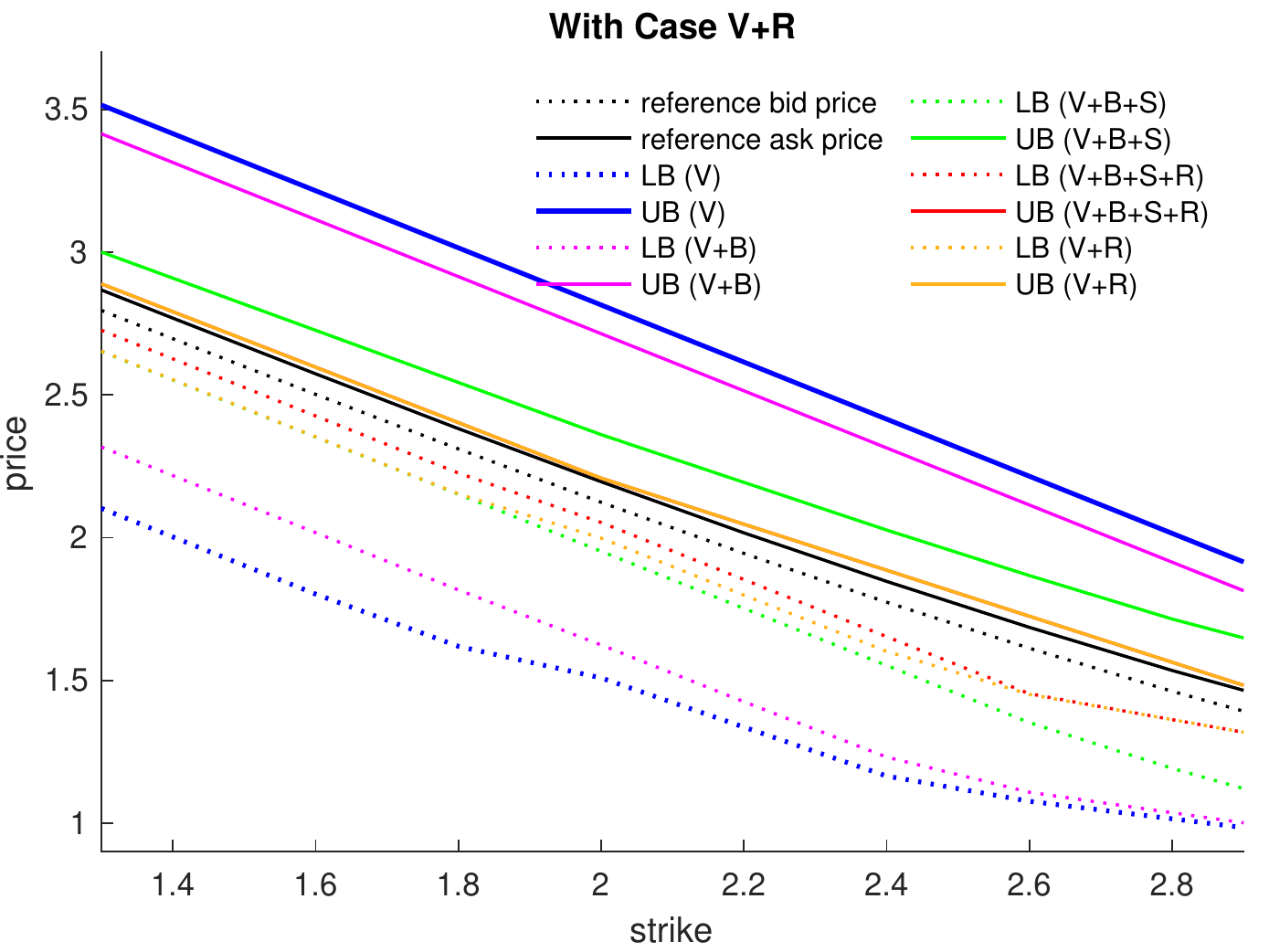}
\caption{Experiment 1 -- Model-free lower and upper price bounds of call-on-max options with strikes between 0 and 10. The bottom left panel shows a magnified version of a part of the top right panel. The bottom right panel shows a magnified version of a part of the top right panel with the Case~\textit{V+R} included.}
\label{fig:exp1}
\end{figure}

We compute the lower and upper bounds of the call-on-max option with payoff function $f$ using the ECP method (Algorithm~\ref{alg:ecpalgo} with Line~\ref{alglin:ecpinitaux} removed) and the ACCP method (Algorithm~\ref{alg:accpalgo}). The inputs of the two algorithms for this experiment are specified in Section~\ref{ecssec:exp1}. 

Figure~\ref{fig:exp1} shows the computed lower and upper price bounds of the call-on-max option with different strikes, along with their reference bid and ask prices.
Let us point out that the price bounds computed by the two algorithms are almost identical. 
Indeed, we have checked that all of the absolute differences between the bounds computed by the two algorithms are below $\varepsilon=0.001$. 
This is a consequence of Theorem~\ref{thm:ecpalg}(ii) and Theorem~\ref{thm:accpalg}(ii) and confirms the correctness of the computed price bounds.

The following observations ensue from the price bounds computed in this example; see Figure~\ref{fig:exp1} for illustrations.
(i) The price bounds in the Cases~1--4 are distinct, and the gap between the lower and upper bounds shrinks when the prices of more traded derivatives are added. 
This means that observing the market prices of more traded derivatives substantially restricts the class of possible pricing measures~$\CQ$ and reduces the no-arbitrage gap between the bounds.
On the dual side (\ref{eqn:superreplicate}), this can be equivalently interpreted as having the information about more traded derivatives provides more ways to sub-replicate and super-replicate the given payoff function and thus makes the gap between the sub-replication price and the super-replication price smaller.
(ii) The addition of rainbow options (\textit{R}) in Case~4 results in a significant reduction of the no-arbitrage gap.
For example, in Case~4, when the strike is 3.2, the upper bound is $2.07\%$ higher than the reference ask price and the lower bound is $7.88\%$ lower than the reference bid price. 
The respective percentages in Case~3 are $15.09\%$ and $26.39\%$ for the upper and lower bounds.
The reason is that the traded call-on-max options provide more information to determine the price of the target derivative, since they are similar in structure to the target derivative. 
This becomes concrete when one considers the dual optimization problem (\ref{eqn:superreplicate}), where these call-on-max options offer direct ways to sub-replicate and super-replicate the target payoff, \textit{e.g.} $\big(x_2 \vee x_3 - \kappa\big)^+ \le \big(x_2 \vee x_3 \vee x_4 - \kappa\big)^+ \le \big(x_1 \vee x_2 \vee x_3 \vee x_4 - \kappa\big)^+$.
(iii) On the other hand, the addition of spread options (Case~3) to vanilla and basket options (Case~2) only yields a significant improvement to the bounds for small strikes ($\le4$), because spread options with large strikes are not traded in the market.
(iv) In all cases, we observe that the no-arbitrage gap is significantly smaller for the (integer) strikes, where traded derivative prices are present.
In the bottom left panel of Figure~\ref{fig:exp1}, for example, the red line (Case~\textit{V+B+S+R}) almost touches the black line (reference price) for the upper bound at strikes 3 and 4, being only $0.47\%$ and $0.30\%$ higher than the respective reference ask prices, while for the lower bound the gap is very small at the same strikes, being $3.67\%$ and $1.82\%$ lower than the respective reference bid prices.
On the other hand, for the intermediate strikes between 3 and 4, for example when the strike is 3.4, the gap grows to $3.05\%$ for the upper bound and to $13.42\%$ for the lower bound.
This is due to the fact that all of the traded options in the synthetic financial market have integer strike prices. 
In the dual optimization problem (\ref{eqn:superreplicate}), one needs to interpolate traded options with integer strike prices in order to sub-replicate and super-replicate the call-on-max option with non-integer strike prices.
Therefore, whenever possible, practitioners should include in the sub- and super-replicating portfolios derivatives with the same strike price as that of the target derivative, in order to reduce the no-arbitrage gap.
(v) As observed from the bottom right panel of Figure~\ref{fig:exp1}, the upper bounds in Case~5 (\textit{V+R}) coincide with the upper bounds in Case~4, while the lower bounds in Case~5 coincide with Case~3 for strikes between 0 and 1.8, coincide with Case~4 for strikes between 2.6 and 10, and fall between Case~3 and Case~4 for strikes between 2 and 2.4. This shows that although the inclusion of more derivatives produces tighter bounds than fewer derivatives, derivatives with similar payoff structure provide more improvement than other derivatives. 
Therefore, since the computation of price bounds is faster when fewer derivatives are included in the sub- and super-replicating portfolios, practitioners should prioritize including derivatives with payoff structure similar to that of the target derivative when the computation time is limited.


\subsection{Experiment~2}
\label{ssec:exp2}

In this experiment, we consider a financial market with 60 assets ($d=60$).
We consider Setting~2 (\textit{i.e.} Assumption~\ref{asp:setting2}) where $\Omega=[0,100]^{60}$. 
We use Algorithms~\ref{alg:ecpalgo} and \ref{alg:accpalgo} to compute the model-free lower and upper price bounds for a call-on-min option on the first 50 out of the 60 assets, with the strike price ranging from 0~to~1 with an increment of 0.1. 
The purpose of this experiment is to demonstrate that Algorithms~\ref{alg:ecpalgo} and \ref{alg:accpalgo} work even when the number of assets is large. 

A total of 400 financial derivatives are traded in the market ($m=400$). 
These financial derivatives include: the 60 assets, 180 vanilla call options (\textit{V}), 3 basket call options (\textit{B}), 147 spread call options (\textit{S}), and 10 call-on-min options (\textit{R}). 
The bid and ask prices of the assets and derivatives are synthetically generated using the market models specified in Section~\ref{ecssec:exp2}. 
For simplicity, we only consider Cases \textit{V+B+S} and \textit{V+B+S+R} in this experiment.
The inputs of the two algorithms for this experiment are detailed in Section~\ref{ecssec:exp2}.

\begin{figure}[t]
\centering
\includegraphics[width=0.48\linewidth]{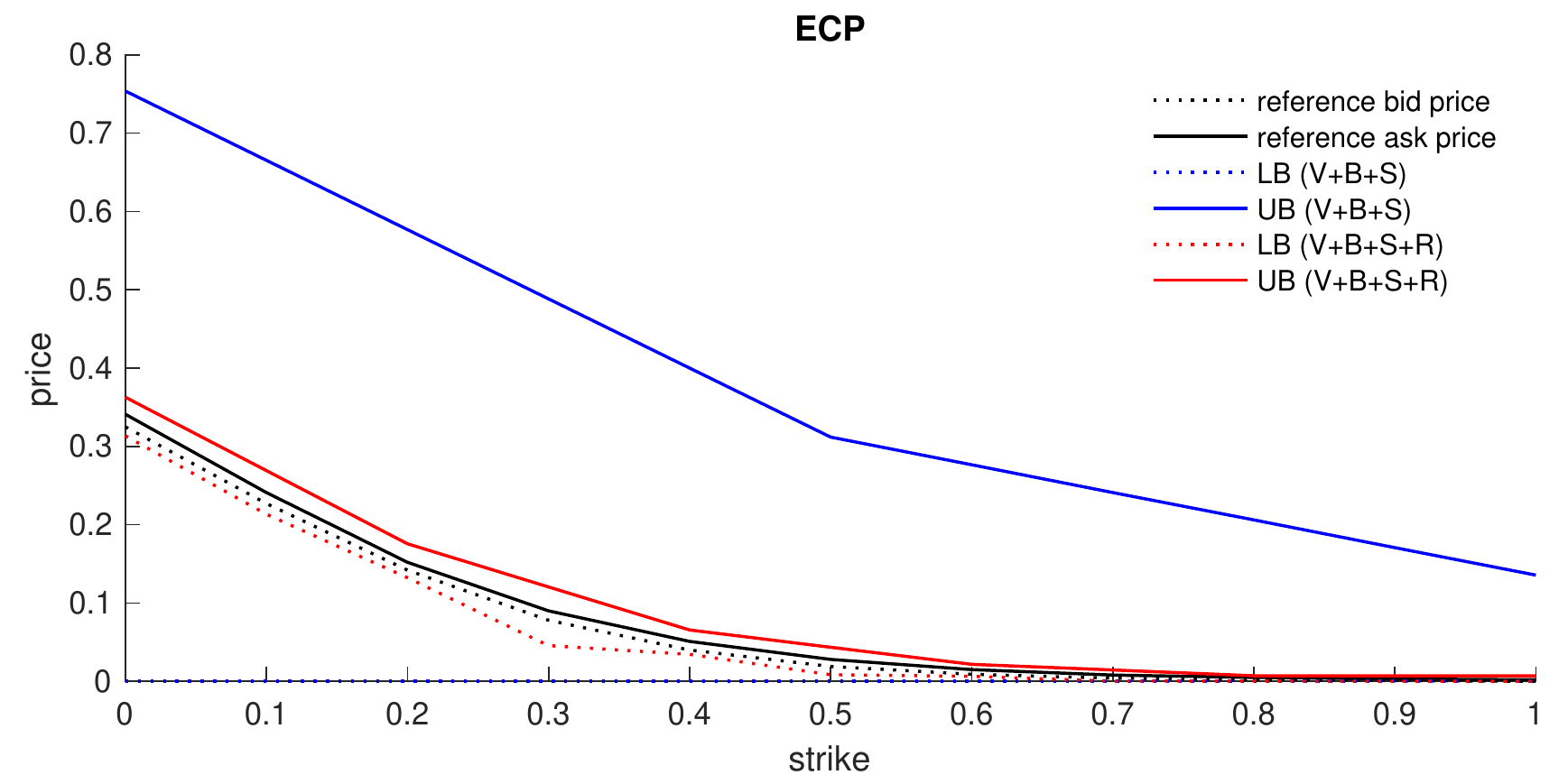}
~
\includegraphics[width=0.48\linewidth]{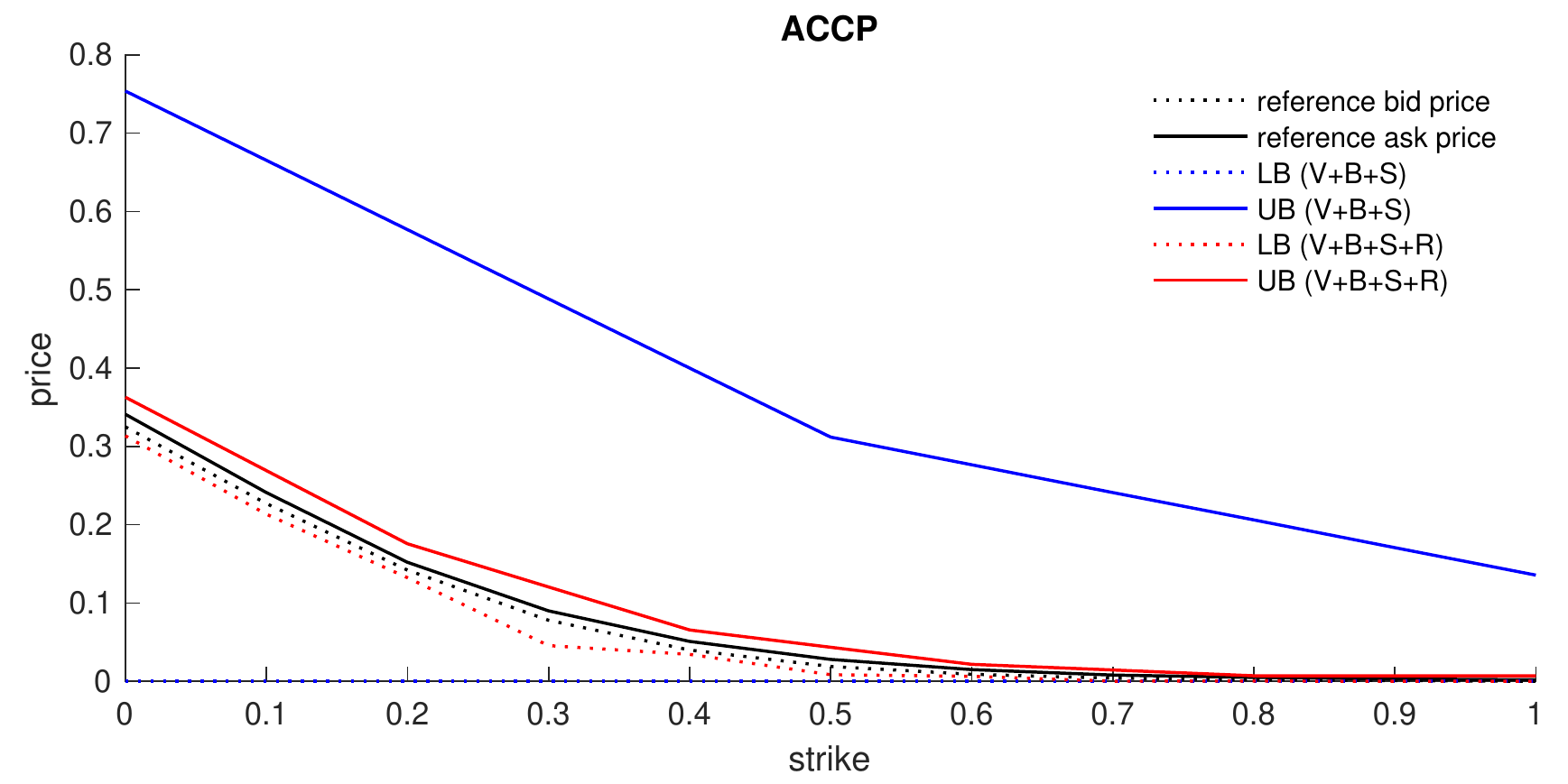}
\caption{Experiment 2 -- Model-free lower and upper price bounds for a call-on-min option with strikes between 0 and 1. }
\label{fig:exp2}
\end{figure}

Figure~\ref{fig:exp2} shows the computed lower and upper price bounds for the call-on-min option with different strikes, along with the reference bid and ask prices. 
Once again, the price bounds computed by the two algorithms are almost identical, and we have checked that all of the absolute differences between the bounds computed by the two algorithms are below $\varepsilon=0.001$.
The following observations ensue from this example, which are mostly in line with the observations from the previous one. 
(i) The price bounds in the two cases are distinct, and the addition of more information improves the bounds and reduces the no-arbitrage gap.
(ii) The addition of traded prices of call-on-min options results in a significant improvement of the bounds, since the payoffs used for sub- and super-replicating and the target payoff are of the same type.
(iii) However, in this high-dimensional example, we notice that the lower price bounds in Case~\textit{V+B+S} are identically zero, showing that the traded vanilla, basket, and spread options do not provide enough information for a non-trivial lower price bound of the call-on-min options, and that it is not possible to sub-replicate the payoff of a call-on-min option with these traded options. 
Therefore, we conclude once again that, whenever possible, practitioners should include in their sub- and super-replicating portfolios not only as many derivatives as possible, but also as many derivatives with similar payoff structure as possible. 

\begin{table}[t]
\caption{Experiment 2 -- Total number of LP and MILP problems solved by the two algorithms. }
\centering
\begin{tabular}{c|c|r|r}
Algorithm & Problem & \textit{V+B+S} & \textit{V+B+S+R} \\
\hline 
\hline 
\multirow{2}{*}{ECP} & LP & 4789 & 3339 \\ 
\cline{2-4} 
& MILP & 4789 & 3339 \\ 
\hline 
\multirow{2}{*}{ACCP} & LP & 1639 & 1714 \\ 
\cline{2-4} 
 & MILP & 1461 & 1574 \\ 
\end{tabular} 
\label{tab:exp2lpmilp}
\end{table}

Table~\ref{tab:exp2lpmilp} shows the total number of LP and MILP problems solved throughout this experiment by the two algorithms. 
The ACCP algorithm achieved convergence faster than the ECP algorithm in this experiment. 
Moreover, in the ACCP algorithm, the MILP problems were only approximately solved with relative gap tolerance $\zeta=0.8$, as explained in Remark~\ref{rmk:accpalg}. 
As a result, the ACCP algorithm was much faster than the ECP algorithm in this experiment.


\subsection{Experiment~3}
\label{ssec:exp3}

In this experiment, we want to demonstrate how the Fundamental Theorem can be combined with the numerical algorithms developed in order to detect arbitrage opportunities in the financial market.
We consider the case where the no-arbitrage assumption (\textit{i.e.} Assumption~\ref{asp:na}) does not hold and use Algorithm~\ref{alg:accpalgo} to detect the presence of arbitrage opportunities in the market. 
We can actually detect a very delicate form of arbitrage, since we consider a financial market with several single-asset options and two multi-asset options; the multi-asset options are priced within their own no-arbitrage intervals, \textit{i.e.} when considered separately from each, there is no arbitrage in the market.
However, when they are considered together an arbitrage opportunity arises and this is detected by the numerical algorithm.
This experiment is inspired by similar examples in \citet[Section~4]{Tavin_2015} and \citet[Section~5.2]{papapantoleon2020detection}; in their setting, the marginals of the pricing measure are given.

We consider Setting~2 (\textit{i.e.} Assumption~\ref{asp:setting2}) where $\Omega=[0,100]^5$, and consider the 5 assets $x_1$, $x_2$, $x_3$, $x_4$, $x_5$ and vanilla call options on the 5 assets with strikes $1,2,\ldots,10$ as the traded financial derivatives. 
In addition, we include a call-on-min option on the five assets with strike~1 and a put-on-min option on the five assets with strike~4. 
The bid and ask prices of the single-asset derivatives are synthetically generated using the method specified at the beginning of Section~\ref{sec:exp}. 

We set the bid and ask prices of the call-on-min option as 0.83 and 0.85, respectively. As for the put-on-min option, we set its bid and ask prices as 3.18 and 3.20, respectively.
Subsequently, we let $f=0$ and run Algorithm~\ref{alg:accpalgo} with $\varepsilon=0.001$, $\tau=1$, $\gamma=0.1$, $\zeta=0.8$, $\delta=0.7$, $\overline{c}=100$, $\overline{\BIy}=100\cdot\vecone$, $\underline{\phi}=\overline{\phi}=0$. 
When only the call-on-min option is considered as traded multi-asset option, Algorithm~\ref{alg:accpalgo} terminates without reaching Line~\ref{alglin:accpunbounded}, and the outputs satisfy $\phi(f)^{\TL\TB}>-0.001$, $\phi(f)^{\TU\TB}=0$. 
Similarly, when only the put-on-min option is considered as traded multi-asset option, Algorithm~\ref{alg:accpalgo} terminates without reaching Line~\ref{alglin:accpunbounded}, and the outputs again satisfy $\phi(f)^{\TL\TB}>-0.001$, $\phi(f)^{\TU\TB}=0$. 
These numerical results imply that there is no arbitrage opportunity in the market with the single-asset derivatives and the call-on-min option, as well as in the market with the single-asset derivatives and the put-on-min option. 
However, when the single-asset derivatives together with both the call-on-min option and the put-on-min option are considered as traded options, Algorithm~\ref{alg:accpalgo} reaches Line~\ref{alglin:accpunbounded} before termination, indicating the violation of Assumption~\ref{asp:na} and the presence of an arbitrage opportunity, as stated by Theorem~\ref{thm:accpalg}(iv). 
The detected arbitrage strategy is given by $c^\star$ and $\BIy^\star$, which specify a portfolio with non-negative payoff such that $c^\star+\pi(\BIy)<0$. 

In Experiment~4 (see Section~\ref{ssec:exp4}), a similar procedure for detecting arbitrage opportunities using Algorithm~\ref{alg:ecpalgo} is applied to Setting~1 (\textit{i.e.} Assumption~\ref{asp:setting1}) where bid and ask prices of traded derivatives are obtained from \textit{real market data}. 
This demonstrates the real-world applicability of the proposed algorithms for arbitrage detection. 


\subsection{Experiment~4}
\label{ssec:exp4}
In this experiment, we use \textit{real market prices} of European call and put options written on the Dow Jones Industrial Average (DJIA) index as well as European call and put options written on the 30 constituent stocks of the DJIA index. This type of market data has been considered by \citet{Hobson_Laurence_Wang_2005_1, dAspremont_ElGhaoui_2006}, and \citet{Pena_Vera_Zuluaga_2010, pena2012computing} for illustrations.

\subsubsection*{Data collection}
The following market prices (corresponding to the closing prices on April 5, 2021 at 16:00 EDT) were collected from MarketWatch\footnote{\url{http://marketwatch.com}} on April 6, 2021. 
\begin{itemize}
\item The prices of the 30 constituent stocks of the DJIA. 
\item The bid and ask prices of the call and put options written on the 30 constituent stocks of the DJIA with expiration date May 21, 2021. 
\item The bid and ask prices of the call and put options written on the SPDR Dow Jones Industrial Average ETF Trust (symbol: DIA), which is an exchange traded fund (ETF) that tracks the DJIA index. 
These DIA options are regarded as basket options written on the 30 constituent stocks with equal weights $w_{\TD\TI\TA}:=0.0658$. The weight $w_{\TD\TI\TA}$ is equal to $\frac{1}{100}$ of the inverse of the Dow divisor\footnote{See \textit{e.g.} \url{https://www.investopedia.com/terms/d/dowdivisor.asp}, accessed: April 30, 2021.} calculated based on the stock prices on April 5, 2021. 
\end{itemize}
\begin{remark}
Even though Experiments~1 and 2 have demonstrated that one should gather the market prices of as many derivatives as possible in order to obtain tight price bounds, exotic options such as, \textit{e.g.}, spread options, call-on-max options, and call-on-min options are usually only traded in over-the-counter (OTC) markets. Therefore, the prices of these exotic options are not publicly available and we are thus unable to collect real market data of this type. 
\end{remark}

\subsubsection*{Data preprocessing}
We apply a procedure using Algorithm~\ref{alg:ecpalgo} to detect whether arbitrage opportunities are present with the bid and ask prices of call and put options written on each of the stocks and on DIA, similar to Experiment~3. 
We found that the prices of DIA options are arbitrage-free, while arbitrage opportunities are present among the prices of options written on 5 of the 30 underlying stocks.
Before proceeding to the next step of the experiment, we adjust the bid and ask prices slightly to remove these arbitrage opportunities. 
In this process, the $l_1$-norm of the price adjustment is minimized to encourage sparsity in the same spirit as \citet{cohen2020detecting}. 
We refer the reader to Section~\ref{ecssec:exp4} in the online appendices for details of this arbitrage removal process. 
After this process, 27 out of the 4304 prices have been adjusted and the largest change (in absolute value) is \$0.38. This shows that only very few options were mis-priced and the market was close to being arbitrage-free. 

\subsubsection*{Experimental setting}
We consider Setting~1 (\textit{i.e.} Assumption~\ref{asp:setting1}) and let $\Omega=\R^{30}_+$. We rank the 30 constituent stocks of the DJIA index based on the market capitalization of the respective companies, \textit{i.e.}, in $(x_1,\ldots,x_{30})^\mathsf{T}\in\Omega$, $x_1$ corresponds to the stock price of the company with the highest market capitalization, and $x_{30}$ corresponds to the stock price of the company with the lowest market capitalization.
Our goal is to use Algorithm~\ref{alg:ecpalgo} to compute the model-free lower and upper price bounds of two basket call options with the following payoff functions:
\begin{align}
\begin{split}
f_1(\BIx)=&\left[\left(\sum_{i=6}^{30}w_{\TD\TI\TA}x_i\right)-\kappa\right]^+, \\ f_2(\BIx)=&\left[\left(\sum_{i=1}^{10}1.2w_{\TD\TI\TA}x_i\right)+\left(\sum_{i=11}^{20}w_{\TD\TI\TA}x_i\right)+\left(\sum_{i=21}^{30}0.8w_{\TD\TI\TA}x_i\right)-\kappa\right]^+,
\end{split}
\label{eqn:weightedboc}
\end{align}
where $\kappa$ is the strike price that is varied in this experiment.
Therefore, $f_1$ is the payoff of a basket call option written on DIA with the five largest companies in terms of market capitalization excluded (\textit{i.e.} a basket option written on a subset of the 30 constituent stocks of the DJIA index). $f_2$ is the payoff of a basket call option written on a weight-adjusted version of DIA, where the weights of the ten largest companies are increased by 20\% and the weights of the ten smallest companies are decreased by 20\%.

After preprocessing, we are left with the following 2152 traded financial derivatives ($m=2152$):
\begin{itemize}
\item 980 vanilla call options and 980 vanilla put options written on the 30 stocks.
\item 96 basket call options and 96 basket put options written on DIA.
\end{itemize}
We consider four different cases when computing the model-free lower and upper bounds. 
In the first three cases, denoted by \textit{V(25\%)}, \textit{V(50\%)}, and \textit{V(100\%)}, we randomly select 25\%, 50\%, and 100\% of the vanilla options, respectively. 
In the fourth case, denoted by \textit{V(100\%)+B}, we use all vanilla and basket options. 
The inputs of the ECP algorithm for this experiment are similar to those in Experiment~1.

\subsubsection*{Experimental results}
\begin{figure}[t]
\centering
\includegraphics[width=0.48\linewidth]{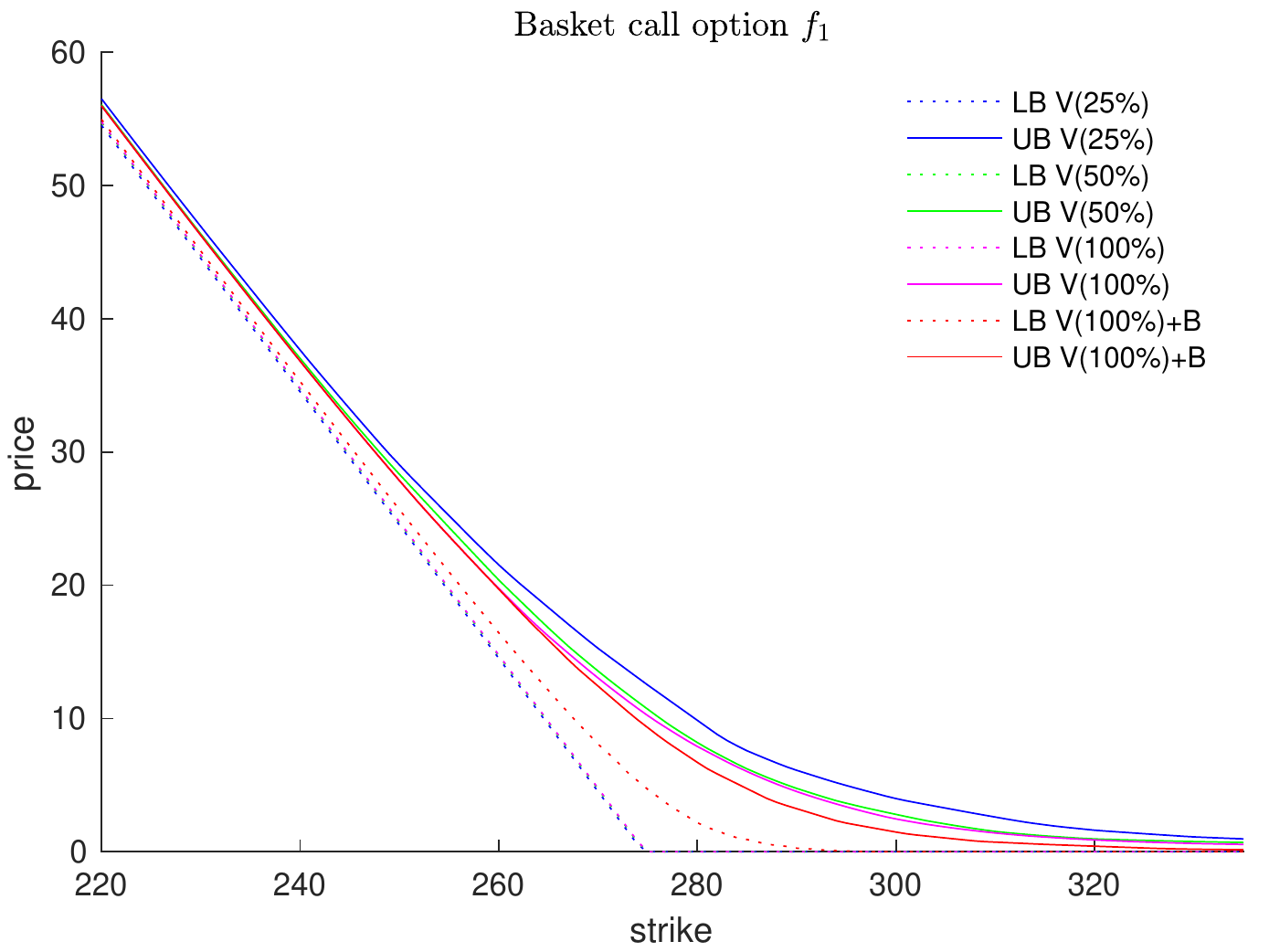}
~
\includegraphics[width=0.48\linewidth]{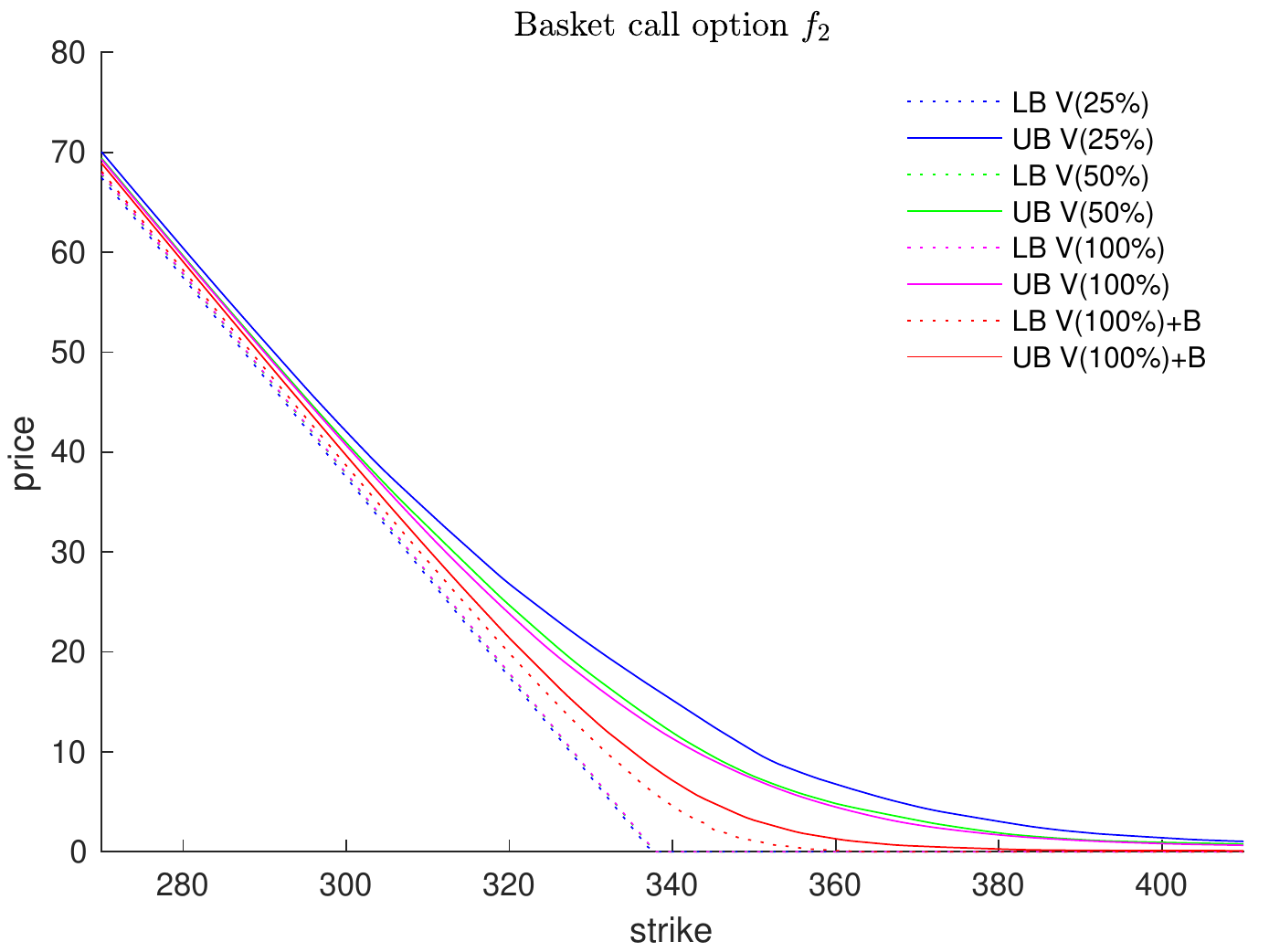}
\caption{Experiment 4 -- Model-free lower and upper price bounds for basket call options.}
\label{fig:exp4}
\end{figure}

Figure~\ref{fig:exp4} shows the computed model-free upper and lower price bounds of the two basket call options with payoff functions $f_1$ and $f_2$ defined in (\ref{eqn:weightedboc}).
When pricing the first basket call option $f_1$, the results from the first three cases show that when more vanilla options are considered, the uncertainty about the marginals of the pricing measure decreases and the model-free upper and lower price bounds are improved, where the improvements of the upper bounds are much more noticeable compared to the lower bounds. 
In Case~\textit{V(100\%)+B}, when other basket options are considered in addition to the vanilla options, the upper price bounds are further improved for large strikes, and the lower price bounds are improved substantially, especially for strikes between \$260 and \$280.
Compared to Case~\textit{V(100\%)}, the inclusion of basket options results in a maximum reduction of \$5.59 in the no-arbitrage gap when the strike is \$275, which reduces the gap by 54.7\%. 
The results when pricing the second basket call option $f_2$ are in line with the results from pricing $f_1$, and the improvements of the upper and lower price bounds in Case~\textit{V(100\%)+B} are more significant across all strikes. 
Compared to Case~\textit{V(100\%)}, the inclusion of basket options results in a maximum reduction of \$9.79 in the no-arbitrage gap when the strike is \$338, which corresponds to a reduction of 79.4\%. 
We conclude from these observations that when pricing a derivative using real market data, the inclusion of more traded derivatives in the sub- and super-replicating portfolios, especially derivatives whose payoff structures are similar to the target derivative, can produce substantial improvements to the arbitrage-free price bounds, and that the proposed algorithms are applicable to real market data.

\begin{APPENDICES}


\section{Proofs of the main results}
\label{sec:proofsmain}




\subsection{Proof of Theorem~\ref{thm:ecpalg}}
\label{ssec:proofsetting1main}

\proof{Proof of Theorem~\ref{thm:ecpalg}}
By the assumption about $\underline{\phi}$ in Remark~\ref{rmk:ecpalg}, there exist $c_0$ and $\BIy_0$ such that $c_0+\langle\BIy_0,\BIg\rangle\ge -f$ and $\underline{\phi}=-c_0-\pi(\BIy_0)$. 
For any $c$ and $\BIy$ such that $c+\langle\BIy,\BIg\rangle\ge f$, it holds that $c_0+c+\langle\BIy_0+\BIy,\BIg\rangle\ge0$ and thus $\langle\BIy_0+\BIy,\BIg\rangle-\pi(\BIy_0+\BIy)\ge-c_0-c-\pi(\BIy_0+\BIy)$.
By Assumption~\ref{asp:na}, $-c_0-c-\pi(\BIy_0+\BIy)\le0$, and so $-c_0-c\le\pi(\BIy_0+\BIy)\le\pi(\BIy_0)+\pi(\BIy)$,
thus $\underline{\phi}=-c_0-\pi(\BIy_0)\le c+\pi(\BIy)$. 
This implies that $\underline{\phi}\le\phi(f)$. 
Let $(\widehat{c},\widehat{\BIy}^+,\widehat{\BIy}^-)$ be an optimizer of (\ref{eqn:lsipdef}), which exists due to Theorem~\ref{thm:duality}(ii). 
We have $\widehat{c}+\pi(\widehat{\BIy}^+-\widehat{\BIy}^-)=\phi(f)>\underline{\phi}-\tau$.
Let $\sigma^{(r)}$ denote the system of linear inequalities $\sigma$ at iteration $r$. 
Proposition~\ref{prop:radconsalgo} states that $\inf_{\BIx\in\R^d_+}s_{\BIy}(\BIx)>-\infty$ if and only if $\BIy\in\R^m$ satisfies all constraints in $\widetilde{\sigma}$. It hence holds by Line~\ref{alglin:ecpinit2} that $(\widehat{c},\widehat{\BIy}^+,\widehat{\BIy}^-)$ satisfies all constraints in $\sigma^{(r)}$ for all $r$. 
Consequently, we have
\begin{align}
\begin{split}
\underline{\varphi}^{(r)}=&\inf_{(c,\BIy^+,\BIy^-)\text{ satisfies }\sigma^{(r)}}\big\{c+\pi(\BIy^+-\BIy^-)\big\}\le\widehat{c}+\pi(\widehat{\BIy}^+-\widehat{\BIy}^-)=\phi(f). 
\end{split}
\label{eqn:ecpalgolb}
\end{align}
$\underline{\varphi}^{(r)}$ is non-decreasing in $r$ because more constraints are added to $\sigma$. 
Moreover, by Proposition~\ref{prop:slackfunc}(iv) and the assumption on $\overline{\BIx}$, for any $\BIy\in\R^m$ that satisfies all constraints in $\widetilde{\sigma}$, it holds that $\inf_{\BIx\in\R^d_+}s_{\BIy}(\BIx)=\inf_{\veczero\le\BIx\le\overline{\BIx}}s_{\BIy}(\BIx)$. 
By Definition~\ref{def:slackfunc}, for all $r$ and for any $\BIx\in\Omega$, 
\begin{align}
\begin{split}
&c^{(r)}-{s}^{(r)}+\langle\BIy^{(r)},\BIg(\BIx)\rangle-f(\BIx)\\
=&c^{(r)}-c^{(r)}-\inf_{\veczero\le\BIx'\le\overline{\BIx}}\Big\{s_{\BIy^{(r)}}(\BIx')\Big\}+\langle\BIy^{(r)},\BIg(\BIx)\rangle-f(\BIx)\\
=&\langle\BIy^{(r)},\BIg(\BIx)\rangle-f(\BIx)-\inf_{\BIx'\in\R^d_+}s_{\BIy^{(r)}}(\BIx')\\
=&\langle\BIy^{(r)},\BIg(\BIx)\rangle-f(\BIx)-\inf_{\BIx'\in\R^d_+}\left\{\langle\BIy^{(r)},\BIg(\BIx')\rangle-f(\BIx')\right\}\\
\ge&0,
\end{split}
\label{eqn:ecpalgofeasshift}
\end{align}
and thus by Line~\ref{alglin:ecplp}, $\underline{\varphi}^{(r)}-{s}^{(r)}=c^{(r)}-{s}^{(r)}+\pi(\BIy^{(r)})\ge\phi(f)$. 
This and (\ref{eqn:ecpalgolb}) also show that $s^{(r)}\le0$.
We have proved statement~(i). 

If Algorithm~\ref{alg:ecpalgo} terminates, then by (\ref{eqn:ecpalgofeasshift}), it holds for all $\BIx\in\R^d_+$ that $c^\star+\langle\BIy^\star,\BIg(\BIx)\rangle-f(\BIx)=c^{(r-1)}-{s}^{(r-1)}+\langle\BIy^{(r-1)},\BIg(\BIx)\rangle-f(\BIx)\ge0$.
Therefore, $(c^\star,\BIy^\star)$ is feasible for (\ref{eqn:superreplicate}) and $\phi(f)^{\TL\TB}\le\phi(f)\le\phi(f)^{\TU\TB}$ follows directly from statement~(i). 
We have ${s}^{(r-1)}\ge-\varepsilon$ at termination, and thus $\phi(f)^{\TU\TB}-\phi(f)^{\TL\TB}\le\varepsilon$ and $(c^\star,\BIy^\star)$ is $\varepsilon$-optimal. 
We now show that Algorithm~\ref{alg:ecpalgo} terminates. 
Proposition~\ref{prop:slackfunc}(iv) states that there exists a partition $\FC$ of $\{\BIx\in\R^d:\veczero\le\BIx\le\overline{\BIx}\}$, such that each $C\in\FC$ is a polytope, $\bigcup_{C\in\FC}C=\{\BIx\in\R^d:\veczero\le\BIx\le\overline{\BIx}\}$, and that for all $\BIy\in\R^m$, $s_{\BIy}(\cdot)$ is an affine function when restricted to each $C\in\FC$.
Let\footnote{A convex subset $C'$ of a convex set $C\subseteq\R^d$ is called a face of $C$ if for all $\lambda\in(0,1)$ and $\BIx_1,\BIx_2\in C$, $\lambda\BIx_1+(1-\lambda)\BIx_2\in C'$ implies that $\BIx_1\in C',\BIx_2\in C'$, see \cite[Section~18]{rockafellar1970convex}.} $\FF:=\{F\ne\emptyset\text{ is a face of some }C\in\FC\}$. 
By Theorem~18.2 of \cite{rockafellar1970convex}, 
\begin{align}
\bigcup_{F\in\FF}\relint(F)=\bigcup_{C\in\FC}C=\{\BIx\in\R^d_+:\veczero\le\BIx\le\overline{\BIx}\}.
\label{eqn:proofecpfaces}
\end{align}
Moreover, by Theorem~19.1 of \cite{rockafellar1970convex}, $|\FF|<\infty$. Let $\BIx^{(r)}$ be a minimizer of the MILP problem in Line~\ref{alglin:ecpmilp}, that is, $c^{(r)}+s_{\BIy^{(r)}}(\BIx^{(r)})={s}^{(r)}$. 
We prove that either Algorithm~\ref{alg:ecpalgo} terminates, or for each $F\in\FF$, there exists at most one $r\in\N$ such that $\BIx^{(r)}\in\relint(F)$. 
Suppose, for the sake of contradiction, that Algorithm~\ref{alg:ecpalgo} does not terminate, and that there exists $r,l\in\N,r<l$, and $\BIx^{(r)},\BIx^{(l)}\in\relint(F)$ for some $F\in\FF$.
Since $\BIx^{(r)}\in X^{(r)}$, we have $c^{(l)}+s_{\BIy^{(l)}}(\BIx^{(r)})\ge0$ by Line~\ref{alglin:ecpfeascut}. We also have $c^{(l)}+s_{\BIy^{(l)}}(\BIx^{(l)})=s^{(l)}<0$, since otherwise Algorithm~\ref{alg:ecpalgo} will terminate at the $l$-th iteration. 
For every $\lambda\in\R$, let $\BIx_\lambda:=(1-\lambda)\BIx^{(r)}+\lambda\BIx^{(l)}$. 
Since $\BIx^{(r)},\BIx^{(l)}\in\relint(F)$, there exists $\widehat{\lambda}>1$ such that $\BIx_{\widehat{\lambda}}\in F\subset C$. 
Since $c^{(l)}+s_{\BIy^{(l)}}(\cdot)$ is an affine function when restricted to the set $C$, we have by $c^{(l)}+s_{\BIy^{(l)}}(\BIx^{(r)})\ge0$ and $c^{(l)}+s_{\BIy^{(l)}}(\BIx^{(l)})<0$ that
\begin{align*}
c^{(l)}+s_{\BIy^{(l)}}(\BIx_{\widehat{\lambda}})=&(1-{\widehat{\lambda}})\Big(c^{(l)}+s_{\BIy^{(l)}}(\BIx^{(r)})\Big)+{\widehat{\lambda}}\Big(c^{(l)}+s_{\BIy^{(l)}}(\BIx^{(l)})\Big)\\
\le&{\widehat{\lambda}}\Big(c^{(l)}+s_{\BIy^{(l)}}(\BIx^{(l)})\Big)<c^{(l)}+s_{\BIy^{(l)}}(\BIx^{(l)}),
\end{align*}
contradicting the fact that $\BIx^{(l)}$ is a minimizer of the MILP problem in Line~\ref{alglin:ecpmilp}. 
Since for each~$r$, $\BIx^{(r)}\in\relint(F)$ for some $F\in\FF$ as a consequence of  (\ref{eqn:proofecpfaces}), and since $|\FF|<\infty$, Algorithm~\ref{alg:ecpalgo} terminates eventually. 
The proof of statement~(ii) is now complete.

Finally, if Line~\ref{alglin:ecpunbounded} of Algorithm~\ref{alg:ecpalgo} is reached, then $\underline{\phi}>\underline{\varphi}^{(r-1)}-{s}^{(r-1)}$. By (\ref{eqn:ecpalgofeasshift}), $c^{(r-1)}-{s}^{(r-1)}+\langle\BIy^{(r-1)},\BIg\rangle\ge f$.
By the assumption about $\underline{\phi}$ in Remark~\ref{rmk:ecpalg}, there exist $c_0$ and $\BIy_0$ such that $c_0+\langle\BIy_0,\BIg\rangle\ge -f$ and $\underline{\phi}=-c_0-\pi(\BIy_0)$. 
Hence, we have $c_0+c^{(r-1)}-{s}^{(r-1)}+\langle\BIy_0+\BIy^{(r-1)},\BIg\rangle\ge0$,
and thus
\begin{align*}
\langle\BIy_0+\BIy^{(r-1)},\BIg\rangle-\pi(\BIy_0+\BIy^{(r-1)})&\ge -c_0-c^{(r-1)}+{s}^{(r-1)}-\pi(\BIy_0+\BIy^{(r-1)})\\
&\ge -c_0-\pi(\BIy_0)-c^{(r-1)}+{s}^{(r-1)}-\pi(\BIy^{(r-1)})\\
&=\underline{\phi}-\underline{\varphi}^{(r-1)}+{s}^{(r-1)}>0,
\end{align*}
which is a violation of Assumption~\ref{asp:na}. 
The proof is now complete. 
\endproof

\subsection{Proofs of Theorem~\ref{thm:accpalg} and Corollary~\ref{coro:ecpprimal}}
\label{ssec:proofsetting2}

\proof{Proof of Theorem~\ref{thm:accpalg}}
If Assumption~\ref{asp:na} holds, then by the same argument as in the proof of Theorem~\ref{thm:ecpalg}(i), $\underline{\phi}\le\phi(f)$. 
Hence, $\underline{\varphi}^{(0)}\le\phi(f)\le\overline{\varphi}^{(0)}$ and $c^{\star(0)}+\langle\BIy^{\star(0)},\BIg\rangle\ge f$ follow from our assumptions. 
For $r\ge1$, suppose that $\underline{\varphi}^{(r-1)}\le\phi(f)$. 
Then, $\underline{\varphi}^{(r)}\ne\underline{\varphi}^{(r-1)}$ only when Line~\ref{alglin:accplblp} is reached. 
This implies that $\sigma(\overline{c},\overline{\BIy},\underline{\varphi}^{(r-1)},{\varphi}^{(r)},X)=\emptyset$. 
By the assumption in Remark~\ref{rmk:accpalg}, there exists an optimizer $(\widehat{c},\widehat{\BIy}^+,\widehat{\BIy}^-)$ of (\ref{eqn:lsipdef})
that satisfies $|\widehat{c}|\le\overline{c}-1$, $\veczero\le\widehat{\BIy}^{+}\le\overline{\BIy}-\vecone$, $\veczero\le\widehat{\BIy}^{-}\le\overline{\BIy}-\vecone$.
In particular, $\widehat{c}+\langle\widehat{\BIy}^{+}-\widehat{\BIy}^{-},\BIg(\BIx)\rangle\ge f(\BIx)$ for all $\BIx\in X$ and
\begin{align}
\widehat{c}+\langle\widehat{\BIy}^{+},\overline{\Bpi}\rangle-\langle\widehat{\BIy}^{-},\underline{\Bpi}\rangle=\widehat{c}+\pi(\widehat{\BIy})=\phi(f).
\label{eqn:accpthmstep1}
\end{align}
By (\ref{eqn:accpthmstep1}) and by the assumption that $\underline{\varphi}^{(r-1)}\le\phi(f)$, it holds that $\sigma(\overline{c},\overline{\BIy},\underline{\varphi}^{(r-1)},{\varphi}^{(r)},X)=\emptyset$ implies $\varphi^{(r)}<\phi(f)$. Therefore, by Line~\ref{alglin:accplblp}, $\underline{\varphi}^{(r)}>\varphi^{(r)}>\underline{\varphi}^{(r-1)}$.
Since $(\widehat{c},\widehat{\BIy}^{+},\widehat{\BIy}^{-})$ also satisfies all constraints in the LP problem in Line~\ref{alglin:accplblp}, we have, again by (\ref{eqn:accpthmstep1}), that $\underline{\varphi}^{(r)}\le\phi(f)$. 
By induction, we have proved that $\underline{\varphi}^{(r)}$ is non-decreasing in $r$ and $\underline{\varphi}^{(r)}\le\phi(f)$ for all $r$. 

For $r\ge1$, $\overline{\varphi}^{(r)}\ne\overline{\varphi}^{(r-1)}$, $c^{\star(r)}\ne c^{\star(r-1)}$, or $\BIy^{\star(r)}\ne \BIy^{\star(r-1)}$ only if Line~\ref{alglin:accpobjcut} is reached. 
By Line~\ref{alglin:accpobjcutcheck} and Line~\ref{alglin:accpobjcut}, $\overline{\varphi}^{(r)}<\overline{\varphi}^{(r-1)}$. 
By the same reasoning as in the proof of Theorem~\ref{thm:ecpalg}(i) in equation~(\ref{eqn:ecpalgofeasshift}), we have $c^{\star(r)}+\langle\BIy^{\star(r)},\BIg\rangle\ge f$ and $\overline{\varphi}^{(r)}\ge\phi(f)$. 
We have thus proved statement~(i). 

If Assumption~\ref{asp:na} holds and Algorithm~\ref{alg:accpalgo} terminates, then $\phi(f)^{\TL\TB}\le\phi(f)\le\phi(f)^{\TU\TB}$ and the feasibility and $\varepsilon$-optimality of $(c^\star,\BIy^\star)$ follow directly from statement~(i) and Line~\ref{alglin:accploop}. 
Thus, we only need to show that Algorithm~\ref{alg:accpalgo} terminates. 
Notice that the Strong Slater Condition in Theorem~1 of \citet{betro2004accelerated} holds because one may take $(\widehat{c},\widehat{\BIy}^+,\widehat{\BIy}^-)$ defined earlier
and choose any $0<\eta<\frac{1}{2}$. Subsequently, 
one checks that 
\begin{align*}
&|\widehat{c}+\eta|\le\overline{c}-\eta,\quad\hspace{4em}\eta\vecone\le\widehat{\BIy}^{+}+\eta\vecone\le\overline{\BIy}-\eta\vecone,\\
&\eta\vecone\le\widehat{\BIy}^{-}+\eta\vecone\le\overline{\BIy}-\eta\vecone,\quad(\widehat{c}+\eta)+\langle(\widehat{\BIy}^{+}+\eta\vecone)-(\widehat{\BIy}^{-}+\eta\vecone),\BIg\rangle\ge f+\eta.
\end{align*}
Thus, $(\widehat{c}+\eta,\widehat{\BIy}^{+}+\eta\vecone,\widehat{\BIy}^{-}+\eta\vecone)$ satisfies the Strong Slater Condition.
Moreover, under Assumption~\ref{asp:setting2}, $\sup_{\BIx\in\Omega}\|\BIg(\BIx)\|<\infty$. 
Suppose, for the sake of contradiction, that Algorithm~\ref{alg:accpalgo} loops infinitely and does not terminate. 
Then, one can deduce that after finitely many iterations, Line~\ref{alglin:accplblp} is never reached, since each time Line~\ref{alglin:accplblp} is reached, $\overline{\varphi}^{(r)}-\underline{\varphi}^{(r)}\le\frac{3}{4}(\overline{\varphi}^{(r-1)}-\underline{\varphi}^{(r-1)})$ by Line~\ref{alglin:accpbisect1} and Line~\ref{alglin:accpbisect2}. 
Similarly, Line~\ref{alglin:accpobjcut} is never reached again after finitely many iterations since each time Line~\ref{alglin:accpobjcut} is reached $\overline{\varphi}^{(r)}-\underline{\varphi}^{(r)}\le\overline{\varphi}^{(r-1)}-\underline{\varphi}^{(r-1)}-\varepsilon$. 
The rest of the proof of statement~(ii) follows exactly as the proof of Theorem~1 in \cite{betro2004accelerated}. 

For statement~(iii), notice that since $\underline{\varphi}^{(0)}=\underline{\phi}-\tau\le\phi(f)-\tau<\phi(f)^{\TL\TB}$, Line~\ref{alglin:accplblp} is reached at least once before termination. 
Thus, $c^\dagger$ and $\BIy^\dagger$ are defined. Let $\BIy^{+\dagger}$ and $\BIy^{-\dagger}$ be defined in Line~\ref{alglin:accplblp}. 
Then, by Line~\ref{alglin:accplblp}, $(c^\dagger,\BIy^{+\dagger},\BIy^{-\dagger})$ is an optimal solution of the LP problem:
\begin{align}
\begin{split}
\text{minimize }\quad& c+\langle\BIy^{+},\overline{\Bpi}\rangle-\langle\BIy^{-},\underline{\Bpi}\rangle\\
\text{subject to }\quad&c+\langle\BIy^+-\BIy^-,\BIg(\BIx)\rangle \ge f(\BIx)\quad\forall\BIx\in X^\dagger,\\
&-\overline{c}\le c\le\overline{c},\;\veczero\le\BIy^{+}\le\overline{\BIy},\;\veczero\le\BIy^{-}\le\overline{\BIy}.
\end{split}
\label{eqn:accpprimallpwithbound}
\end{align}
Thus, since $\underline{\varphi}^{(r)}$ is updated whenever $(c^\dagger,\BIy^\dagger)$ are updated, we have $\phi(f)^{\TL\TB}=c^\dagger+\langle\BIy^{+\dagger},\overline{\Bpi}\rangle-\langle\BIy^{-\dagger},\underline{\Bpi}\rangle$, and $c^\dagger+\langle\BIy^{+\dagger}-\BIy^{-\dagger},\BIg(\BIx)\rangle\ge f(\BIx)$ for all $\BIx\in X^\dagger$.
Let $\widetilde{B}:=\left\{(c,\BIy^+,\BIy^-):-\overline{c}\le c\le\overline{c},\;\veczero\le\BIy^{+}\le\overline{\BIy},\;\veczero\le\BIy^{-}\le\overline{\BIy}\right\}\subset\R^{2m+1}$.
By the assumption of statement~(iii), $-\overline{c}<c^\dagger<\overline{c}$, $\veczero\le\BIy^{+\dagger}<\overline{\BIy}$, $\veczero\le\BIy^{-\dagger}<\overline{\BIy}$, and we claim that $(c^\dagger,\BIy^{+\dagger},\BIy^{-\dagger})$ is also optimal for the following LP problem:
\begin{align}
\begin{split}
\text{minimize }\quad& c+\langle\BIy^{+},\overline{\Bpi}\rangle-\langle\BIy^{-},\underline{\Bpi}\rangle\\
\text{subject to }\quad&c+\langle\BIy^+-\BIy^-,\BIg(\BIx)\rangle\ge f(\BIx)\quad\forall\BIx\in X^\dagger,\\
&\BIy^+\ge\veczero,\;\BIy^-\ge\veczero.
\end{split}
\label{eqn:accpprimallp}
\end{align}
Suppose, for the sake of contradiction, that (\ref{eqn:accpprimallp}) has optimal solution $(\widetilde{c},\widetilde{\BIy}^{+},\widetilde{\BIy}^{-})$ with $\kappa:=\widetilde{c}+\langle\widetilde{\BIy}^{+},\overline{\Bpi}\rangle-\langle\widetilde{\BIy}^{-},\underline{\Bpi}\rangle<\phi(f)^{\TL\TB}$. 
Then, since $\phi(f)^{\TL\TB}$ is the optimal value of (\ref{eqn:accpprimallpwithbound}), we have $(\widetilde{c},\widetilde{\BIy}^{+},\widetilde{\BIy}^{-} )\notin \widetilde{B}$, $\widetilde{\BIy}^{+}\ge\veczero$, $\widetilde{\BIy}^{-}\ge\veczero$. 
Let $c_{\lambda}:=\lambda c^\dagger+(1-\lambda)\widetilde{c}$, $\BIy^+_\lambda:=\lambda\BIy^{+\dagger}+(1-\lambda)\widetilde{\BIy}^+$, $\BIy^-_\lambda:=\lambda\BIy^{-\dagger}+(1-\lambda)\widetilde{\BIy}^-$. 
Then, there exists some $\lambda\in(0,1)$, such that $(c_\lambda,\BIy^+_\lambda,\BIy^-_\lambda)=\lambda(c^\dagger,\BIy^{+\dagger},\BIy^{-\dagger})+(1-\lambda)(\widetilde{c},\widetilde{\BIy}^+,\widetilde{\BIy}^-)\in \widetilde{B}$,
$c_\lambda+\langle\BIy^+_{\lambda}-\BIy^-_{\lambda},\BIg(\BIx)\rangle\ge f(\BIx)$ for all $\BIx\in X^\dagger$, 
and $c_\lambda+\langle\BIy^{+}_\lambda,\overline{\Bpi}\rangle-\langle\BIy^{-}_\lambda,\underline{\Bpi}\rangle=\lambda\phi(f)^{\TL\TB}+(1-\lambda)\kappa<\phi(f)^{\TL\TB}$, contradicting the optimality of $(c^\dagger,\BIy^{+\dagger},\BIy^{-\dagger})$ for (\ref{eqn:accpprimallpwithbound}). 
Therefore, $(c^\dagger,\BIy^{+\dagger},\BIy^{-\dagger})$ is also optimal for (\ref{eqn:accpprimallp}), whose corresponding dual LP problem is exactly (\ref{eqn:accpdual}). 
Then, an optimal solution $(\mu^\star_{\BIx})_{\BIx\in X^\dagger}$ of (\ref{eqn:accpdual}) exists, its corresponding finitely supported measure $\mu^\star$ is a probability measure which satisfies $\underline{\pi}_j\le\int_{\Omega}g_j \Td\mu^\star\le\overline{\pi}_j$ for $j=1,\ldots,m$, and thus $\mu^\star\in\CQ$. 
Moreover, due to the strong duality of LP problems, $\int_{\Omega}f\Td\mu^\star=\phi(f)^{\TL\TB}\ge\phi(f)-\varepsilon$ by statement~(ii), and $\mu^\star$ is $\varepsilon$-optimal for the right-hand side of (\ref{eqn:duality}) by Theorem~\ref{thm:duality}(iii). We have completed the proof of statement~(iii). 

The proof of statement~(iv) is exactly the same as the proof of Theorem~\ref{thm:ecpalg}(iii). The proof is now complete. 
\endproof

\proof{Proof of Corollary~\ref{coro:ecpprimal}}
The proof that Algorithm~\ref{alg:ecpalgo} terminates is identical to the proof of Theorem~\ref{thm:ecpalg}(ii). 
Hence, as in Theorem~\ref{thm:ecpalg}(ii), we have $\phi(f)^{\TL\TB}\ge\phi(f)-\varepsilon$. Since Line~\ref{alglin:ecpinitaux} is not used, and that $\phi(f)^{\TL\TB}\ge\underline{\phi}-\varepsilon>\underline{\phi}-\tau$, we have that $\phi(f)^{\TL\TB}$ is the optimal value of the following LP problem:
\begin{align*}
\text{minimize }\quad&c+\langle\BIy^+,\overline{\Bpi}\rangle-\langle\BIy^-,\underline{\Bpi}\rangle\\
\text{subject to }\quad&c+\langle\BIy^+-\BIy^-,\BIg(\BIx)\rangle\ge f(\BIx)\quad\forall\BIx\in X=\bigcup_{l=0}^{r-1}X^{(l)},
\end{align*}
whose dual LP problem is exactly (\ref{eqn:ecpprimallp}). Consequently, by the argument in the proof of Theorem~\ref{thm:accpalg}(iii), $\int_{\Omega}f\Td\mu^\star=\phi(f)^{\TL\TB}\ge\phi(f)-\varepsilon$, and $\mu^\star$ is $\varepsilon$-optimal for the right-hand side of (\ref{eqn:duality}) by Theorem~\ref{thm:duality}(iii). 
\endproof

\end{APPENDICES}

\ACKNOWLEDGMENT{We thank Daniel Bartl, Julien Guyon, and Steven Vanduffel for fruitful comments and discussions that initiated this project.
We are also grateful to Stephan Eckstein and Michael Kupper for fruitful discussions during the work on these topics. AN
gratefully acknowledges the financial support by his Nanyang Assistant Professorship Grant (NAP Grant) \textit{Machine Learning
based Algorithms in Finance and Insurance}. AP gratefully acknowledges the financial support from the
Hellenic Foundation for Research and Innovation Grant No. HFRI-FM17-2152.}

\bibliographystyle{informs2014} 
\input{ModelFree_arXiv_main.bbl} 

\ECSwitch


\ECHead{Online Appendices}

\renewcommand{\theequation}{\thesection.\arabic{equation}}
\numberwithin{equation}{section}

\section{Additional remarks about the duality result}
\label{ecsec:duality}

\begin{remark}
Theorem~\ref{thm:duality} is similar in spirit to the multi-marginal optimal transport problem, possibly under additional constraints, see \textit{e.g.} \citet{kellerer1984duality}, \citet{Rachev_Rueschendorf_1994}, \citet{Zaev_2015}, and \citet{bartl2017marginal}. 
However, in contrast to these articles, we do not assume that the marginals are known and fixed.
Moreover, notice the subtle differences in the attainment of supremum and infimum. 
In Theorem~\ref{thm:duality}, the infimum in (\ref{eqn:superreplicate}) corresponding to the super-replication portfolio is attained, whereas the supremum on the right-hand side of (\ref{eqn:duality}) corresponding to the ``worst-case'' probability measure is not necessarily attained. 
In the multi-marginal optimal transport problem, the supremum corresponding to the ``worst-case'' probability measure (or optimal coupling of the marginals) is attained due to compactness, which holds when the marginals of the pricing measure are fixed. 
\end{remark}

The following example demonstrates that the supremum on the right-hand side of (\ref{eqn:duality}) is not necessarily attained. 
\begin{example}
Let $\Omega=\R_+$, $m=1$, and, for $x\in\R_+$, set:
\begin{align*}
g_1(x) = x, \quad f(x)=\max\{x-1,0\}, \quad \underline{\pi}_1=0, \quad\text{and }\overline{\pi}_1=1.
\end{align*}
Clearly, 
\begin{align*}
&c+y_1g_1(x)-f(x)\ge0\quad\forall x\in\R_+\\
\Longleftrightarrow\quad&c+y_1x-\max\{x-1,0\}\ge0\quad\forall x\in\R_+\\
\Longrightarrow\quad&c\ge0,\;y_1\ge1.
\end{align*}
In the case where $c=0$, $y_1=1$, $c+y_1g_1(x)-f(x)\ge0$ $\forall x\in\R_+$ holds, and we have that $c^\star=0$, $y_1^\star=1$ is an optimizer of (\ref{eqn:superreplicate}).
Thus $\phi(f)=c^\star+y_1^\star\overline{\pi}_1=1$. On the other hand, if $\mu$ is a probability measure on $\R_+$, $\underline{\pi}_1\le\int_{\R_+}g_1\Td\mu\le\overline{\pi}_1$, and $\int_{\R_+} f\Td\mu=1$, then 
\begin{align*}
\int_{\R_+} \min\{x,1\}\Td\mu=\int_{\R_+} (g_1-f)\Td\mu=\int_{\R_+} g_1\Td\mu-\int_{\R_+} f\Td\mu\le\overline{\pi}_1-1=0.
\end{align*}
Since $\min\{x,1\}>0$ for $x>0$, it is implied that $\mu((0,+\infty))=0$, which is impossible. Hence, $\sup_{\mu\in\CQ}\int_{\Omega} f\Td\mu$ is not attained. 
\label{ex:primalnotattained}
\end{example}

The following proposition presents a specific setting in which Assumption~\ref{asp:na} holds, and its proof can be found in Section~\ref{sec:proofduality}.
\begin{proposition}
Let $\Omega\subseteq\R^d$ for some $d\in\N$, $g_1,\ldots,g_m$ be continuous functions on $\Omega$, and  assume there exists $\widehat{\mu}\in\CQ$ such that $\widehat{\mu}$ is equivalent to the Lebesgue measure on $\Omega$.
Then, Assumption~\ref{asp:na} holds. 
\label{prop:examplena}
\end{proposition}

\section{Additional notes about the numerical methods}
\label{ecsec:numerical}

\subsection{Notes about CPWA functions}
\label{ecssec:cpwa}

\begin{example}\label{ex:CPWA-payoffs}
Many popular payoff functions in finance belong to the class of CPWA functions.
The list below contains those payoff functions used in the present work, alongside their CPWA representation.
\begin{enumerate}[label=(\roman*)]
	\item Setting $K=1, \xi_1=1, I_1=2, \BIa_{1,1}=\BIe_i, \BIa_{1,2}=\veczero$ and $b_{1,1}=-\kappa, b_{1,2}=0$, we get the payoff of a call option on the $i$-th asset with strike $\kappa$, via
		\[
			h(\BIx) = \max \{ x_i-\kappa,0 \} = (x_i-\kappa)^+.
		\] 
	\item Using the previous setting and replacing $\BIa_{1,1}$ with $\BIa_{1,1}=\BIw=(w_1,\ldots,w_d)^{\mathsf{T}}\in\R^d_+$, we get the payoff of a basket call option on the $d$ assets with weights $\BIw$ and strike $\kappa$, via
		\[
			h(\BIx) = \max\Big\{ \sum_{i=1}^dw_ix_i-\kappa,0 \Big\} = \Big( \sum_{i=1}^dw_ix_i-\kappa \Big)^+.
		\]
	\item Using the previous setting, and replacing $\BIa_{1,1}$ with $\BIa_{1,1} = \sum_{i\in \CI} \BIe_i - \sum_{j\in \CI'} \BIe_j$, for $\CI,\CI'\subset \{1,\dots,d\}$ and $\CI\cap\CI'=\emptyset$, we get the payoff of a spread call option with strike $\kappa$, via
		\[
			h(\BIx) = \max\Big\{ \sum_{i\in \CI} x_i - \sum_{j\in \CI'} x_j-\kappa, 0 \Big\}  
					= \Big(  \sum_{i\in \CI} x_i - \sum_{j\in \CI'} x_j-\kappa \Big)^+.
		\] 
	\item Setting $K=1, \xi_1=1,I_1=d+1, \BIa_{1,i}=\BIe_i, b_{1,i}=-\kappa$ for all $i\in\{1,\dots,d\}$, $\BIa_{1,d+1}=\veczero$ and $b_{1,d+1}=0$, we get the payoff of a call-on-max option on the $d$ assets with strike $\kappa$, via
		\begin{align*}
			h(\BIx) &= \max\big\{ x_1-\kappa, \dots, x_d-\kappa, 0 \big\} 
					 = \max\big\{ \max\big\{x_1, \dots, x_d\big\}- \kappa, 0 \big\} \\ 
					&= \big( x_1 \vee \dots \vee x_d - \kappa \big)^+.
		\end{align*}
	\item Setting $K=2, \xi_1=1, \xi_2=-1, I_1=d+1, I_2=d, \BIa_{1,i}=\BIa_{2,i}=-\BIe_i, b_{1,i}=b_{2,i}=\kappa$ for all $i\in\{1,\dots,d\}$, and $\BIa_{1,d+1}=\veczero$, $b_{1,d+1}=0$, we get the payoff of a call-on-min option on the $d$ assets with strike $\kappa$, via
		\begin{align*}
			h(\BIx) &= \max\big\{ \kappa-x_1, \dots, \kappa-x_d, 0 \big\} - \max\big\{ \kappa-x_1, \dots, \kappa-x_d \big\}\\
					&= \max\big\{ \min\big\{\kappa-x_1, \dots, \kappa-x_d\big\}, 0 \big\} \\ 
					&= \big( x_1 \wedge \dots \wedge x_d - \kappa \big)^+.
		\end{align*}
	\item Finally, setting $K=1, \xi_1=1, I_1=d+1, \BIa_{1,i}=\BIe_i, b_{1,i}=-\kappa_i$ for all $i\in\{1,\dots,d\}$, and $\BIa_{1,d+1}=\veczero, b_{1,d+1}=0$, we get the payoff of a best-of-call option on the $d$ assets with strikes $\kappa_1,\dots,\kappa_d$, via
		\begin{align*}
			h(\BIx) &= \max\big\{ x_1-\kappa_1, \dots, x_d-\kappa_d, 0 \big\} \\
					&= \max\big\{ (x_1-\kappa_1)^+, \dots, (x_d-\kappa_d)^+ \big\} \\ 
					&= (x_1-\kappa_1)^+ \vee (x_2-\kappa_2)^+ \vee \dots \vee (x_d-\kappa_d)^+.
		\end{align*}
\end{enumerate}	
Let us point out that these representations are not unique.
Moreover, by replacing the vectors $\BIe_i$ or $-\BIe_i$ with suitable vectors in examples~(iii), (iv), (v), and (vi) above, one can create weighted versions of the aforementioned payoffs. 
These can be interpreted as options written on a number of indices. 
\end{example}


The following properties of CPWA functions can be deduced from Definition~\ref{def:cpwa}.
\begin{lemma}
The following are properties of CPWA functions.
\begin{enumerate}
\item[(i)] A finite linear combination of CPWA functions is again a CPWA function;
\item[(ii)] If $h:\R^d_+\to\R$ is a CPWA function, then it admits the following local representation:
\begin{align}
h(\BIx)=\begin{cases}
\langle\BIa_{1},\BIx\rangle+b_{1}, & \text{if }\BIx\in\Omega_1,\\
\qquad\vdots & \quad\vdots\\
\langle\BIa_{J},\BIx\rangle+b_{J}, & \text{if }\BIx\in\Omega_J,
\end{cases}
\label{eqn:cpwalocal}
\end{align}
where for $j=1,\ldots,J$, $\BIa_{j}\in\R^d$, $b_{j}\in\R$, $\Omega_j$ is a polyhedron, and $\bigcup_{j=1}^J\Omega_j=\R^d_+$. 
In addition, for $j\ne j'$, $\inter(\Omega_j)\cap\inter(\Omega_{j'})=\emptyset$. 
If $\BIx\in\Omega_j\cap\Omega_{j'}$, then $\langle\BIa_{j},\BIx\rangle+b_{j}=\langle\BIa_{j'},\BIx\rangle+b_{j'}$, and thus the above function is well-defined. 
\item[(iii)] In the local representation (\ref{eqn:cpwalocal}), it can be assumed without loss of generality that $\inter(\Omega_j)\ne\emptyset$ for $j=1,\ldots,J$. 
\end{enumerate}
\label{lem:cpwaproperties}
\end{lemma}
\proof{Proof of Lemma~\ref{lem:cpwaproperties}}
See Section~\ref{sec:proofalgo}. 
\endproof

The following proposition establishes the connection between the lower-boundedness of a CPWA function and the non-negativity of its associated radial function. 
\begin{proposition}
Let $h:\R^d_+\to\R$ be a CPWA function, and let $\widetilde{h}:\R^d_+\to\R$ be its radial function. 
Then, 
\begin{enumerate}
\item[(i)] $\inf_{\BIx\in\R^d_+}h(\BIx)>-\infty$ if and only if $\widetilde{h}(\BIz)\ge0$ for all $\BIz\in\R^d_+$;
\item[(ii)] if $\inf_{\BIx\in\R^d_+}h(\BIx)>-\infty$, there exists $\BIx^\star\in\R^d_+$ such that $h(\BIx^\star)=\inf_{\BIx\in\R^d_+}h(\BIx)$. 
\end{enumerate}
\label{prop:cpwaboundedness}
\end{proposition}
\proof{Proof of Proposition~\ref{prop:cpwaboundedness}}
See Section~\ref{sec:proofalgo}.
\endproof

\subsection{Notes about the exterior cutting plane (ECP) algorithm}
\label{ecssec:ecp}

A crucial step in the ECP algorithm (Algorithm~\ref{alg:ecpalgo}) is to minimize the slack function ${s}_{\BIy}(\BIx)$ over $\BIx\in\R^d_+$. 
In order to do so, one needs to first guarantee that $\inf_{\BIx\in\R^d_+}{s}_{\BIy}(\BIx)>-\infty$. Proposition~\ref{prop:slackfunc} below shows that this is equivalent to guaranteeing that the slack function $\widetilde{s}_{\BIy}(\cdot)$ of $s_{\BIy}(\cdot)$ satisfies $\widetilde{s}_{\BIy}(\BIz)\ge0$ for all $\BIz\in\R^d_+$. Moreover, Proposition~\ref{prop:slackfunc} states that, in order to minimize ${s}_{\BIy}(\BIx)$ over $\BIx\in\R^d_+$, one only needs to minimize over a compact set $\{\BIx\in\R^d:\veczero\le\BIx\le\overline{\BIx}\}$ for some $\overline{\BIx}>\veczero$\footnote{$\BIx=(x_1,\ldots,x_d)^{\mathsf T}>\veczero  \Leftrightarrow  x_1>0,\ldots,x_d>0.$}. 

\pagebreak
\begin{proposition}
Under Assumption~\ref{asp:setting1}, the following statements hold.
\begin{enumerate}
\item[(i)] For any fixed $\BIy\in\R^m$, the slack function $s_{\BIy}:\R^d_+\to\R$ is a CPWA function. 
		For any fixed $\BIx\in\R^d_+$, $s_{\BIy}(\BIx)$ is linear when regarded as a function of $\BIy$. 
		In particular, $s_{\BIy}(\BIx)\ge -c$ corresponds to a linear inequality constraint on $(c,\BIy)$ for every fixed $\BIx\in\R^d_+$. 
\item[(ii)] There exist $K\in\N$, $\BIw_k\in\R^d$, $z_k\in\R$, $I_k\in\N$ for $k=1,\ldots,K$, $\BIa_{k,i}\in\R^d$, $b_{k,i}\in\R$, for $i=1,\ldots,I_k,k=1,\ldots,K$, such that
	\begin{align}
	\begin{split}
		s_{\BIy}(\BIx)=\sum_{k=1}^K(\langle\BIy,\BIw_k\rangle+z_k)\max\left\{\langle\BIa_{k,i},\BIx\rangle+b_{k,i}:1\le i\le I_k\right\}. 
	\end{split}
	\label{eqn:lsipconspoint}
	\end{align}
	Let the radial function of $s_{\BIy}$ be denoted by $\widetilde{s}_{\BIy}$. 
	There exist $\widetilde{K}\in\N$, $\widetilde{\BIw}_k\in\R^d$, $\widetilde{z}_k\in\R$, $\widetilde{I}_k\in\N$ for $k=1,\ldots,K$, $\widetilde{\BIa}_{k,i}\in\R^d$, for $i=1,\ldots,\widetilde{I}_k,k=1,\ldots,\widetilde{K}$, such that
	\begin{align}
	\begin{split}
		\widetilde{s}_{\BIy}(\BIz)=\sum_{k=1}^{\widetilde{K}}(\langle\BIy,\widetilde{\BIw}_k\rangle+\widetilde{z}_k)\max\left\{\langle\widetilde{\BIa}_{k,i},\BIz\rangle:1\le i\le \widetilde{I}_k\right\}.
	\end{split}
	\label{eqn:lsipconsdir}
	\end{align}
\item[(iii)] The slack function $s_{\BIy}(\BIx)$ is bounded from below as a function of $\BIx$ if and only if $\widetilde{s}_{\BIy}(\BIz)\ge0$ for all $\BIz\in\R^d_+$. Moreover, in this case $\inf_{\BIx\in\R^d_+}s_{\BIy}(\BIx)$ is attained at some $\BIx^\star$.
\item[(iv)] There exists $\overline{\BIx}\in\R^d$ with $\overline{\BIx}>\veczero$ such that for every $\BIy\in\R^m$ that satisfies $\inf_{\BIx\in\R^d_+}s_{\BIy}(\BIx)>-\infty$, it holds that $\inf_{\veczero\le\BIx\le\overline{\BIx}}s_{\BIy}(\BIx)=\inf_{\BIx\in\R^d_+}s_{\BIy}(\BIx)$. Moreover, there exists a collection $\FC$ of polytopes such that $\bigcup_{C\in\FC}C=\{\BIx\in\R^d:\veczero\le\BIx\le\overline{\BIx}\}$ and that for every $\BIy\in\R^m$, $s_{\BIy}(\cdot)$ is an affine function when restricted to each $C\in\FC$. 
\end{enumerate}
\label{prop:slackfunc}
\end{proposition}

\proof{Proof of Proposition~\ref{prop:slackfunc}}
See Section~\ref{sec:proofalgo}. 
\endproof


The next proposition translates the condition: $\widetilde{s}_{\BIy}(\BIz)\ge0$ for all $\BIz\in\R^d_+$ into sufficient and necessary linear constraints (\ref{eqn:lpconsdir}) on $\BIy$ with auxiliary variables. 
\begin{proposition}
Let $\widetilde{s}_{\BIy}(\BIz)$ be given by (\ref{eqn:lsipconsdir}).
Let $\widetilde{\CI}=\{(i_k)_{k=1:\widetilde{K}}:1\le i_k\le\widetilde{I}_k, k=1,\ldots,\widetilde{K}\}$. 
For $(i_k)\in\widetilde{\CI}$, define $A_{(i_k)}:=\Big\{\widetilde{\BIa}_{k,i_k}-\widetilde{\BIa}_{k,i}:i\in\{1,\ldots,\widetilde{I}_k\}\setminus\{i_k\},k=1,\ldots,\widetilde{K}\Big\}\setminus\{\veczero\}$.
Let $\dual(C)$ denote the dual cone of a set $C\subset\R^d$. 
Then, the following statements hold.
\begin{enumerate}
\item[(i)] 
For each $(i_k)\in\widetilde{\CI}$, it holds that 
\begin{align*}
\widetilde{s}_{\BIy}(\BIz)=\sum_{k=1}^{\widetilde{K}}(\langle\BIy,\widetilde{\BIw}_k\rangle+\widetilde{z}_k)\langle\widetilde{\BIa}_{k,i_k},\BIz\rangle,
\end{align*}
for all $\BIz\in\dual(A_{(i_k)})\cap\R^d_+$ and $\BIy\in\R^m$. 
We have $\bigcup_{(i_k)\in\widetilde{\CI}}\dual(A_{(i_k)})\cap\R^d_+=\R^d_+$. If $(i_k)\in\widetilde{\CI}$, $(i'_k)\in\widetilde{\CI}$, $(i_k)\ne(i'_k)$, then
\begin{align*}
\inter(\dual(A_{(i_k)})\cap\R^d_+)\cap\inter(\dual(A_{(i'_k)})\cap\R^d_+)=\emptyset.
\end{align*}
Moreover, by possibly removing some elements from $\widetilde{\CI}$, one can assume without loss of generality that $\inter(\dual(A_{(i_k)})\cap\R^d_+)\ne\emptyset$ for all $(i_k)\in\widetilde{\CI}$. 
\item[(ii)] For each $(i_k)\in\widetilde{\CI}$, $\widetilde{s}_{\BIy}(\BIz)\ge0$ for all $\BIz\in\dual(A_{(i_k)})\cap\R^d_+$ if and only if for each $\BIv\in A_{(i_k)}$ there exists $\eta^{(i_k)}_\BIv\ge0$ such that 
\begin{align}
\sum_{k=1}^{\widetilde{K}}(\langle\BIy,\widetilde{\BIw}_k\rangle+\widetilde{z}_k)\widetilde{\BIa}_{k,i_k}\ge \sum_{\BIv\in A_{(i_k)}}\eta^{(i_k)}_{\BIv} \BIv.
\label{eqn:lpconsdir}
\end{align}
\item[(iii)] For each $(i_k)\in\widetilde{\CI}$, $\inter(\dual(A_{(i_k)})\cap\R^d_+)=\emptyset$ is equivalent to the following condition:
\begin{align}
\begin{split}
\text{there exists}\quad&\eta^{(i_k)}_{\BIv}\ge 0\text{ for each }\BIv\in A_{(i_k)},\\
\text{such that}\quad&\sum_{\BIv\in A_{(i_k)}}\eta^{(i_k)}_{\BIv}=1,\\
&\sum_{\BIv\in A_{(i_k)}}\eta^{(i_k)}_{\BIv} \BIv\le\veczero,
\end{split}
\label{eqn:lpconsdircheck}
\end{align}
which can be checked by a linear programming (LP) solver.
\end{enumerate}
\label{prop:radcons}
\end{proposition}

\proof{Proof of Proposition~\ref{prop:radcons}}
See Section~\ref{sec:proofalgo}. 
\endproof

\setcounter{algocf}{-1}
\begin{algorithm}[t]
\KwIn{$\widetilde{s}_{\BIy}(\BIz)=\sum_{k=1}^{\widetilde{K}}(\langle\BIy,\widetilde{\BIw}_k\rangle+\widetilde{z}_k)\max\left\{\langle\widetilde{\BIa}_{k,i},\BIz\rangle:1\le i\le \widetilde{I}_k\right\}$ for all $\BIz\in\R^d_+$ and $\BIy\in\R^m$}
\KwOut{System of linear inequalities on $\BIy$ with auxiliary variables to guarantee $\widetilde{s}_{\BIy}(\BIz)\ge0$ for all $\BIz\in\R^d_+$, denoted by $\widetilde{\sigma}$}
\nl Initialize $\widetilde{\sigma}$ to be empty. \\
\nl \label{alglin:radloop}\For{every $(i_k)\in\widetilde{\CI}$}{
\nl Compute $A_{(i_k)}\leftarrow\{\widetilde{\BIa}_{k,i_k}-\widetilde{\BIa}_{k,i}:1\le i\le\widetilde{I}_k,i\ne i_k,1\le k\le\widetilde{K}\}\setminus\{\veczero\}$. \label{alglin:radvectors}\\
\nl \label{alglin:radlp}\If{the condition (\ref{eqn:lpconsdircheck}) fails to hold}{
\nl Expand $\widetilde{\sigma}$ to include additional non-negative auxiliary variables $(\eta^{(i_k)}_{\BIv})_{\BIv\in A_{(i_k)}}$. \label{alglin:radaux}\\
\nl Append the linear inequality~(\ref{eqn:lpconsdir}) to $\widetilde{\sigma}$. \label{alglin:radagg}
}
}
\nl \Return $\widetilde{\sigma}$.
    \caption{{\bf Generation of Radial Constraints}}
    \label{alg:radcons}
\end{algorithm}

Based on Proposition~\ref{prop:radcons}, Algorithm~\ref{alg:radcons} generates the necessary and sufficient constraints on $\BIy$ to guarantee that $\widetilde{s}_{\BIy}(\BIz)\ge0$ for all $\BIz\in\R^d_+$, which are jointly referred to as the \textit{radial constraints}. 
The output of Algorithm~\ref{alg:radcons} is a system of linear inequalities of the form (\ref{eqn:lpconsdir}) and is denoted by $\widetilde{\sigma}$. 
Apart from $\BIy$, these linear inequalities also contain auxiliary variables. 
In Line~\ref{alglin:radaux}, the system of linear inequalities is ``expanded'' to include additional auxiliary variables. 
In Line~\ref{alglin:radagg}, a new linear inequality is added to the current inequality system. 
This is sometimes referred to as \textit{aggregating} the constraint to $\widetilde{\sigma}$ in the literature. 
Concretely, if the linear inequality constraints are stored in a matrix, then Line~\ref{alglin:radaux} corresponds to adding extra columns to the matrix to increase the number of decision variables, and Line~\ref{alglin:radagg} corresponds to adding an extra row to the matrix as an additional constraint. 
The correctness of Algorithm~\ref{alg:radcons} is shown by the following proposition.

\begin{proposition}
Let $s_{\BIy}(\BIz)$ be the slack function given by (\ref{eqn:lsipconspoint}) and let $\widetilde{s}_{\BIy}(\BIz)$ be its radial function given by (\ref{eqn:lsipconsdir}). Let $\widetilde{\sigma}$ be the system of linear inequalities returned by Algorithm~\ref{alg:radcons}. Then, for any $\BIy\in\R^m$, $\inf_{\BIz\in\R^d_+}s_{\BIy}(\BIz)>-\infty$ if and only if $\BIy$ satisfies all linear inequalities in $\widetilde{\sigma}$. 
\label{prop:radconsalgo}
\end{proposition}
\proof{Proof of Proposition~\ref{prop:radconsalgo}}
See Section~\ref{sec:proofalgo}. 
\endproof

\begin{remark}
\begin{enumerate}
\item[(i)]
	In Algorithm~\ref{alg:radcons}, auxiliary variables are introduced because each $\cone(A_{(i_k)})$ is represented by the conic hull of its extreme directions $A_{(i_k)}$ (\textit{i.e.} a $V$-polyhedron). 
	Such a cone also admits an equivalent representation by the intersection of finitely many half-spaces (\textit{i.e.} an $H$-polyhedron), that is, 
	\begin{align*}
	\cone(A_{(i_k)})=\big\{\BIv\in\R^d_+:\langle\BIv,\BIu\rangle\ge 0,\forall\BIu\in\mathrm{exdir}(\dual(\cone(A_{(i_k)})))\big\},
	\end{align*}
	where $\mathrm{exdir}(\cdot)$ denotes the set of extreme directions of a polyhedron. 
	Using this representation, there is no need to introduce auxiliary variables. 
	Unfortunately, to obtain such a representation, one needs to enumerate all the extreme directions of $\dual(\cone(A_{(i_k)}))$, which is usually computationally very costly (see, for example, \citet{le1992note} and Section~5.2 of \citet{ciripoibensolve}).
\item[(ii)]
	The computational complexity of Algorithm~\ref{alg:radcons} may grow fast as $\widetilde{K}$ grows. 
	In certain cases, the growth can be exponential. 
	The proposed approach is only practical when the total number of times Line~\ref{alglin:radagg} is executed is not too large. 
	Fortunately, this is usually the case with the problems we study, since vanilla options and basket options with non-negative weights do not increase this number. 
\end{enumerate}
\end{remark}

After imposing constraints on $\BIy$ to guarantee that $\inf_{\BIx\in\R^d_+}s_{\BIy}(\BIx)>-\infty$, Algorithm~\ref{alg:ecpalgo} solves the global optimization problem $\min_{\veczero\le\BIx\le\overline{\BIx}}s_{\BIy}(\BIx)$ by formulating it into a mixed-integer linear programming (MILP) problem. This is detailed in the following lemma. 

\begin{lemma}
Let $h:\R^d_+\to\R$ be a CPWA function represented as \eqref{eqn:cpwadef}. 
Assume that $\inf_{\BIx\in\R^d_+}h(\BIx)>-\infty$. 
Let $\overline{\BIx}=(\overline{x}_1,\ldots,\overline{x}_d)^{\mathsf T}>\veczero$. 
For $i=1,\ldots,I_k$, $k=1,\ldots,K$, let $M_{k,i}$ be the constant given by
\begin{align}
M_{k,i}:=\max\{\langle\BIa_{k,i'}-\BIa_{k,i},\BIx\rangle+b_{k,i'}-b_{k,i}:\veczero\le\BIx\le\overline{\BIx},1\le i'\le I_k,i'\ne i\}.
\label{eqn:milpM}
\end{align}
Then, the global minimization problem
\begin{align}
\begin{split}
\text{minimize }\quad& h(\BIx)=\sum_{k=1}^{K}\xi_k\max\left\{\langle\BIa_{k,i},\BIx\rangle+b_{k,i}:1\le i\le I_k\right\}\\
\text{subject to }\quad&\veczero\le\BIx\le\overline{\BIx}
\end{split}
\label{eqn:milpglobalmin}
\end{align}
is equivalent to the following mixed-integer linear programming (MILP) problem:
\begin{align}
\begin{split}
\text{minimize }\quad&\sum_{k=1,\ldots,K,\xi_k=1}\lambda_k+\sum_{k=1,\ldots,K,\xi_k=-1}-\zeta_k\\
\text{subject to }\quad
&\text{for }k=1,\ldots,K \text{ and }\xi_k=1:\\
&\quad\quad\quad\begin{cases}
\lambda_k\in\R,\\
\langle\BIa_{k,i},\BIx\rangle+b_{k,i}\le\lambda_k,&\forall 1\le i\le I_k,
\end{cases}\\
\quad& \text{for }k=1,\ldots,K \text{ and }\xi_k=-1:\\
&\quad\quad\quad\begin{cases}
\zeta_k\in\R,\\
\delta_{k,i}\in\R,&\forall 1\le i\le I_k,\\
\iota_{k,i}\in\{0,1\} &\forall 1\le i\le I_k,\\
\langle\BIa_{k,i},\BIx\rangle+b_{k,i}+\delta_{k,i}=\zeta_k &\forall 1\le i\le I_k,\\
0\le\delta_{k,i}\le M_{k,i}(1-\iota_{k,i}) &\forall 1\le i\le I_k,
\end{cases}\\
&\quad\quad\quad \sum_{i=1}^{I_k}\iota_{k,i}=1,\\
&\veczero\le\BIx\le\overline{\BIx}.
\end{split}
\label{eqn:milp}
\end{align}
\label{lem:milp}
\end{lemma}
\proof{Proof of Lemma~\ref{lem:milp}}
See Section~\ref{sec:proofalgo}. 
\endproof

\begin{remark}[The branch-and-bound (BnB) algorithm]
Consider the MILP problem: $\text{minimize }p(\Bgamma)\text{ subject to }\Bgamma\in \Gamma$, where $\Gamma$ denotes the feasible set.
The BnB algorithm generates a sequence of solutions $(\Bgamma_{t})\subset \Gamma$ that satisfy the equality/inequality constraints and integer constraints with non-increasing objective values, that is, for $\overline{p}_t:=p(\Bgamma_t)$, we have
\begin{align*}
\overline{p}_1\ge\overline{p}_2\ge\ldots\ge p^\star,
\end{align*}
where $p^\star$ denotes the optimal value of the problem. 
These solutions are usually referred to as integer feasible solutions. 
At the same time, the BnB algorithm produces a sequence of relaxations of the original problem with non-decreasing objective values to serve as lower bounds of the original problem, that is,
\begin{align*}
\underline{p}_1\le\underline{p}_2\le\ldots\le p^\star.
\end{align*}
The BnB algorithm terminates on iteration $t$ when the following condition is satisfied:
\begin{align}
\frac{\overline{p}_t-\underline{p}_t}{|\overline{p}_t|}\le \zeta,
\label{eqn:bnbtermin}
\end{align}
where $\zeta\in(0,1)$ is referred to as the relative gap tolerance. 
At termination, $\Bgamma_t$ is returned as the approximate optimizer, $\overline{p}_t$ is returned as the approximate optimal value, and $\underline{p}_t$ is returned as the best objective bound. 
Gurobi Optimization provides a basic overview of mixed-integer programming\footnote{\url{https://www.gurobi.com/resource/mip-basics/}, accessed: April 18, 2021.}. 
\label{rmk:bnb}
\end{remark}

\begin{remark}[Inputs of Algorithm~\ref{alg:ecpalgo}]
The following list explains the various inputs of Algorithm~\ref{alg:ecpalgo} and how to set them.
\begin{itemize}
\item $\overline{\Bpi}$, $\underline{\Bpi}$, $(g_j)_{j=1:m}$, $f$ are given by problem~(\ref{eqn:lsipdef}). 
\item $X^{(0)}$ is the initial set of points in $\R^d_+$ used to generate the initial constraints on $c$ and $\BIy$. 
	It is possible to set $X^{(0)}=\emptyset$ and the LP problem in Line~\ref{alglin:ecplp} remains bounded due to the additional constraint introduced in Line~\ref{alglin:ecpinitlb}. 
	The return value $X$ (or a subset of $X$) can be reused as $X^{(0)}$ for another run of Algorithm~\ref{alg:ecpalgo} with a similar problem to achieve significant speed-up. 
\item $\overline{\BIx}$ is the upper bound used in the MILP problem~(\ref{eqn:milp}). 
	For the validity of the algorithm, $\overline{\BIx}$ should satisfy the conditions in Proposition~\ref{prop:slackfunc}(iv). 
\item $\underline{\phi}$ specifies an initial lower bound of $\phi(f)$. 
	In order to guarantee $\underline{\phi}\le\phi(f)$, one finds $c_0$ and $\BIy_0$ such that $c_0+\langle\BIy_0,\BIg\rangle\ge -f$ and set $\underline{\phi}=-c_0-\pi(\BIy_0)$. The existence of such $c_0$ and $\BIy_0$ is guaranteed by Assumption~\ref{asp:setting1}(iii) and Theorem~\ref{thm:duality}(ii). 
	Under Assumption~\ref{asp:na}, one can show that $\underline{\phi}\le\phi(f)$ always holds. 
	This is detailed in Theorem~\ref{thm:ecpalg}. 
\item $\varepsilon$ is a positive number indicating the numerical accuracy of the algorithm. 
	This is detailed in Theorem~\ref{thm:ecpalg}(ii). 
\item $\tau$ can be set as any positive number to provide a strict lower bound on $\phi(f)$. 
\item $\delta\in(0,1]$ controls the number of additional constraints generated in each iteration of the algorithm in Line~\ref{alglin:ecpintfeas}, each of which corresponds to a sub-optimal integer feasible solution found by the BnB algorithm before the optimizer (see Remark~\ref{rmk:bnb}). 
	When $\delta$ is set to 1, $X^{(r)}$ will be a singleton and contains only the $\BIx^\star$ that corresponds to an optimizer of problem~(\ref{eqn:milp}). 
	When $\delta$ is set close to 0, $X^{(r)}$ corresponds to most of the integer feasible solutions with $\BIx$ such that $c^{(r)}+s_{\BIy^{(r)}}(\BIx)<0$ and thus more constraints are included in the next LP problem in Line~\ref{alglin:ecplp}. 
	$\delta$ can usually be set small, since we empirically observe that it often results in a reduction of the total number of iterations. 
\end{itemize}
\label{rmk:ecpalg}
\end{remark}

\begin{remark}
Algorithm~\ref{alg:ecpalgo} is inspired by the Conceptual Algorithm~11.4.1 in \citet{goberna1998linear}. 
There are a few notable differences, as detailed below. 
\begin{itemize}
\item In Line~\ref{alglin:ecpinitaux}, the initial constraints involve auxiliary variables from the radial constraints generated by Algorithm~\ref{alg:radcons}.
\item Line~\ref{alglin:ecpinitlb} introduces a lower bound on $\phi(f)$ to guarantee that the LP problem in Line~\ref{alglin:ecplp} is bounded. 
\item Line~\ref{alglin:ecpintfeas} uses sub-optimal integer feasible solutions of the MILP problem~(\ref{eqn:milp}), allowing for more than one cuts to be generated in each iteration of Algorithm~\ref{alg:ecpalgo}, which significantly speeds up the algorithm. 
\item In contrast to the Conceptual Algorithm~11.4.1 in \cite{goberna1998linear} which returns an infeasible solution in general, Algorithm~\ref{alg:ecpalgo} returns a feasible solution due to Line~\ref{alglin:ecpfeas}. 
	This is discussed in detail in Theorem~\ref{thm:ecpalg}.  
\end{itemize}
\label{rmk:ecpalgdiff}
\end{remark}


\subsection{Notes about the accelerated central cutting plane (ACCP) algorithm}
\label{ecssec:accp}

\begin{remark}
It is not possible to compute a Chebyshev of a $P$-polytope resulting from a projection (\textit{e.g.} a system of linear inequalities containing auxiliary variables) by solving an LP problem. 
Thus, Algorithm~\ref{alg:accpalgo} is not applicable in Setting~1 (\textit{i.e.} Assumption~\ref{asp:setting1}). 
\citet{zhen2018computing} have made a detailed exposition of various approaches to approximately compute the Chebyshev center of such a polytope. 
However, the computational cost of these approaches is too high to be considered for our problem. 
\end{remark}

\begin{remark}[Inputs of Algorithm~\ref{alg:accpalgo}]
The following list explains the various inputs of Algorithm~\ref{alg:accpalgo} and how to set them. 
\begin{itemize}
\item $\overline{\Bpi}$, $\underline{\Bpi}$, $(g_j)_{j=1:m}$, $f$ are given by problem~(\ref{eqn:lsipdef}). 
\item $X^{(0)}$ is the same as in Algorithm~\ref{alg:ecpalgo} and Remark~\ref{rmk:ecpalg}. 
\item $\overline{\BIx}$ is specified in Assumption~\ref{asp:setting2}. 
\item $\underline{\phi}$ specifies an initial lower bound of $\phi(f)$ and is obtained in the same way as in Remark~\ref{rmk:ecpalg}. 
\item $\overline{c}>1$, $\overline{\BIy}>\vecone$ specify a bounding box on $(c,\BIy^+,\BIy^-)$ as in (\ref{eqn:accppolytope}). We specify $\overline{c}$ and $\overline{\BIy}$ to be large enough such that 
there exists an optimizer $(\widehat{c},\widehat{\BIy}^+,\widehat{\BIy}^-)$ of (\ref{eqn:lsipdef}) that satisfies $|\widehat{c}|\le\overline{c}-1$, $\widehat{\BIy}^+\le\overline{\BIy}-\vecone$, and $\widehat{\BIy}^-\le\overline{\BIy}-\vecone$.
This is to guarantee the validity of Algorithm~\ref{alg:accpalgo}.
\item $(c^{\star(0)},\BIy^{\star(0)})$ specifies an initial feasible point of (\ref{eqn:superreplicate}),
and $\overline{\phi}=c^{\star(0)}+\pi(\BIy^{\star(0)})$ is an initial upper bound of $\phi(f)$. 
Moreover, we specify $c^{\star(0)},\BIy^{\star(0)}$ such that $\big|c^{\star(0)}\big|\le\overline{c}-1$ and $-\overline{\BIy}+\vecone\le\BIy^{\star(0)}\le\overline{\BIy}-\vecone$.
\item $\varepsilon$ is a positive number specifying the numerical accuracy of the algorithm. This is detailed in Theorem~\ref{thm:accpalg}(ii) and~(iii). 
\item $\tau$ can be set as any number greater than $\varepsilon$ to provide a strict lower bound on $\phi(f)$. 
\item $\gamma\in[0,1)$ controls the removal of active feasibility cuts in Line~\ref{alglin:accpremove}. 
	When $\gamma$ is set to 0, then Line~\ref{alglin:accpremove} is never reached and all feasibility cuts are retained. 
	If $\gamma$ is set close to 1, then more feasibility cuts are removed, making solving the LP problem in Line~\ref{alglin:accpcclp} and Line~\ref{alglin:accplblp} faster. However, when $\gamma$ is close to 1, Algorithm~\ref{alg:accpalgo} may take more iterations to converge.
\item $\zeta\in(0,1)$ controls the termination of the BnB algorithm, hence also the solution quality of the MILP problem in Line~\ref{alglin:accpmilp}. 
	Suppose that the true optimal value of the MILP problem is $s^{(r)}$. 
	As in Remark~\ref{rmk:bnb}, in the BnB algorithm, $\overline{s}^{(r)}$ is the approximate optimal value, and $\underline{s}^{(r)}$ is the best objective bound at termination. 
	By the termination condition in (\ref{eqn:bnbtermin}), $\frac{\overline{s}^{(r)}-\underline{s}^{(r)}}{|\overline{s}^{(r)}|}\le \zeta$. 
	If $\overline{s}^{(r)}\ge0$, then $\overline{s}^{(r)}\ge{s}^{(r)}\ge\underline{s}^{(r)}\ge0$. 
	If $\overline{s}^{(r)}<0$, then $(1+\zeta)\overline{s}^{(r)}\le\underline{s}^{(r)}\le {s}^{(r)}\le\overline{s}^{(r)}<0$, hence ${s}^{(r)}\le\overline{s}^{(r)}\le\frac{1}{1+\zeta}{s}^{(r)}$.
\item $\delta$ is the same as in Algorithm~\ref{alg:ecpalgo} and Remark~\ref{rmk:ecpalg}.
\end{itemize}
\label{rmk:accpalg}
\end{remark}

\begin{remark}
Algorithm~\ref{alg:accpalgo} is partially based on the accelerated central cutting plane algorithm in Section~2 of \citet{betro2004accelerated}, with the following differences:
\begin{itemize}
\item Instead of starting with $\underline{\varphi}^{(0)}=-\infty$, Algorithm~\ref{alg:accpalgo} starts with a given lower bound $\underline{\varphi}^{(0)}>-\infty$ similar to Algorithm~\ref{alg:ecpalgo}. 
\item In Algorithm~\ref{alg:accpalgo}, whenever the current polytope $\sigma(\overline{c},\overline{\BIy},\underline{\varphi}^{(r)},{\varphi}^{(r)},X)$ is found to be empty, instead of setting $\underline{\varphi}^{(r)}\leftarrow\varphi^{(r)}$, an LP problem is solved in Line~\ref{alglin:accplblp} to update the lower bound $\underline{\varphi}^{(r)}$. 
\item Similar to Algorithm~\ref{alg:ecpalgo}, in Line~\ref{alglin:accpintfeas} of Algorithm~\ref{alg:accpalgo}, sub-optimal integer feasible solutions found by the MILP solver are used to generate additional feasibility cuts to speed up the algorithm. 
\item In Line~\ref{alglin:accpobjcut} of Algorithm~\ref{alg:accpalgo}, a feasible point of the LSIP problem is generated from a possibly infeasible one by setting $c^{\star(r)}\leftarrow c^{(r)}-\underline{s}^{(r)}$. 
	This can result in an objective cut being made even when the current Chebyshev center is infeasible. 
\end{itemize}
\label{rmk:accpalgdiff}
\end{remark}


\subsection{Details about the implementation of the algorithms}
\label{ecssec:expdetails}

All algorithms and experiments are implemented using the MATLAB language. 
We adopt the MATLAB interface of Gurobi to solve linear programs and mixed-integer linear programs within the algorithms. 
The MATLAB code used in this work is available on GitHub\footnote{\url{https://github.com/qikunxiang/ModelFreePriceBounds}}. 

The list below briefly discusses some numerical considerations in the implementation of the algorithms. 
The details of these considerations can be found in the MATLAB code. 
\begin{itemize}
\item All LP problems in Algorithms~\ref{alg:ecpalgo} and \ref{alg:accpalgo} are solved using the interior point algorithm for better numerical accuracy. 
	Even though the dual simplex algorithm is usually more efficient, we observed empirically that it is numerically less accurate compared to the interior point algorithm. 
\item The efficiency and the numerical stability of the LP solver are sensitive to the numerical condition of the coefficients of the model, which is associated with the ratio between the largest and smallest coefficients in absolute values. 
	Thus, we round the inputs of the model when constructing coefficient matrices of LP problems. 
	Concretely, in Line~\ref{alglin:ecpintfeas} of Algorithm~\ref{alg:ecpalgo} and in Line~\ref{alglin:accpintfeas} of Algorithm~\ref{alg:accpalgo}, each component of a solution $\BIx$ is rounded to a multiple of $10^{-4}$ to control the numerical condition of the coefficient matrix in the LP problem.
\item When formulating the MILP problem in (\ref{eqn:milp}), we discard the terms in (\ref{eqn:milpglobalmin}) in which $\BIa_{k,i}$ and $b_{k,i}$ are very close to zero (\textit{i.e.} $\le 10^{-6}$) to reduce the size of the resulting MILP problem. 
\item In order to prevent Algorithms~\ref{alg:ecpalgo} and \ref{alg:accpalgo} from spending too much time on solving the MILP problem (\ref{eqn:milp}), we set an upper limit on the number of nodes explored by the BnB algorithm. 
This is fine for Algorithm~\ref{alg:ecpalgo}, and for Algorithm~\ref{alg:accpalgo} when $\overline{s}^{(r)}<0$ in Line~\ref{alglin:accpmilp}. 
However, in Line~\ref{alglin:accpmilp} of Algorithm~\ref{alg:accpalgo}, when the node limit of the BnB algorithm is reached and $\overline{s}^{(r)}\ge0$, the algorithm might get stuck when Line~\ref{alglin:accpobjcut} is not subsequently reached. 
To resolve this issue, we heuristically modify $\underline{s}^{(r)}$ to guarantee that Line~\ref{alglin:accpobjcut} is subsequently reached. 
We have observed that the BnB algorithm is almost always able to find the optimal integer feasible solution early on while spending most of the execution time refining the lower bound. 
Moreover, we have empirically observed that the described scenario only occurred when $\phi(f)^{\TU\TB}-\phi(f)^{\TL\TB}$ is very small. 
Therefore, this heuristic adjustment usually does not affect the validity of Algorithm~\ref{alg:accpalgo}. 
\end{itemize}

\section{Additional notes about the numerical experiments}
\label{ecsec:experiments}

\subsection{Notes about Experiment~1}
\label{ecssec:exp1}
In Experiment~1, the bid and ask prices of the derivatives are synthetically generated using market models specified as follows. 
\begin{itemize}
\item The marginal distributions of the asset prices at terminal time are log-normal distributions truncated to $[0,100]$. Two groups of marginal distributions are used. Let $\mu_{j,i}$ and $\sigma^2_{j,i}$ denote the $\mu$ and $\sigma^2$ parameters of the $i$-th marginal in the $j$-th group, for $i=1,\ldots,5$, $j=1,2$. 
In the first group, we set 
\begin{itemize}
\item $\mu_{1,1}=0.5$, $\sigma^2_{1,1}=0.2$, 
\item $\mu_{1,2}=1$, $\sigma^2_{1,2}=0.4$, 
\item $\mu_{1,3}=1$, $\sigma^2_{1,3}=0.2$, 
\item $\mu_{1,4}=0.5$, $\sigma^2_{1,4}=0.4$, 
\item $\mu_{1,5}=0.5$, $\sigma^2_{1,5}=0.2$.
\end{itemize}
In the second group, we set 
\begin{itemize}
\item $\mu_{2,1}=0.5$, $\sigma^2_{2,1}=0.21$, 
\item $\mu_{2,2}=1$, $\sigma^2_{2,2}=0.42$, 
\item $\mu_{2,3}=1$, $\sigma^2_{2,3}=0.21$, 
\item $\mu_{2,4}=0.5$, $\sigma^2_{2,4}=0.42$, 
\item $\mu_{2,5}=0.5$, $\sigma^2_{2,5}=0.21$.
\end{itemize}
\item The dependence model corresponds to a $t$-copula with a correlation matrix which has a 2-factor model structure:
\begin{align*}
	\BC=\BLambda\BD\BLambda^{\mathsf T}+\BPsi,
\end{align*}
where $\BLambda\in\R^{d\times 2}$, $\BD\in\R^{2\times 2}$ is a diagonal matrix, and $\BPsi\in\R^{d\times d}$ is a diagonal matrix. 
Two dependence models with identical correlation matrices and with degrees of freedom equal to 3 and 4 are used. 
\end{itemize}
Effectively, we have now created four market models, each made up of a combination of marginal distributions and a dependence model. 
We have kept the four market models similar in order to simulate realistic bid--ask spreads that are consistent with typical observations from real markets. 

The inputs of Algorithms~\ref{alg:ecpalgo} and \ref{alg:accpalgo} in Experiment~1 are specified as follows. 
In Algorithm~\ref{alg:ecpalgo}, we set $\varepsilon=0.001,\tau=0.1,\delta=0.7$. We set $\underline{\phi}$ according to Remark~\ref{rmk:ecpalg}. 
When computing $\phi(f)$, we set $\underline{\phi}=0$, which is given by $0\ge-f$. 
When computing $\phi(-f)$, we set $\underline{\phi}=-\overline{\pi}_2-\overline{\pi}_3-\overline{\pi}_4$, which is given by $x_2+x_3+x_4\ge f$. 
In Algorithm~\ref{alg:accpalgo}, we set $\varepsilon=0.001,\tau=0.1,\gamma=0.1,\zeta=0.01,\delta=0.7,\overline{c}=100,\overline{\BIy}=100\cdot\vecone$. 
$\underline{\phi}$ is set to be the same as in Algorithm~\ref{alg:ecpalgo} and $\overline{\phi}$ is set according to Remark~\ref{rmk:accpalg}. 
When computing $\phi(f)$, we set $\overline{\phi}=\overline{\pi}_2+\overline{\pi}_3+\overline{\pi}_4$ and $c^{\star(0)}+\langle\BIy^{\star(0)},\BIg\rangle=x_2+x_3+x_4\ge f$.
When computing $\phi(-f)$, we set $\underline{\phi}=0$ and $c^{\star(0)}+\langle\BIy^{\star(0)},\BIg\rangle=0\ge-f$.

\subsection{Notes about Experiment~2}
\label{ecssec:exp2}
In Experiment~2, the bid and ask prices of the derivatives are synthetically generated using market models similar to those used in Experiment~1, as detailed below. 
\begin{itemize}
\item The marginal distributions of the asset prices at terminal time are log-normal distributions truncated to $[0,100]$. Two groups of marginal distributions are used. The $\mu$ parameters of the marginals in the first group are uniformly randomly sampled from $[-0.3,0.1]$, and the $\sigma^2$ parameters of the marginals in the first group are uniformly randomly sampled from $[0.2,0.8]$. 
In the second group, the $\mu$ parameters of the marginals are set to be identical to the first group, while the $\sigma^2$ parameters are increased by a random amount uniformly sampled from $[0,0.1]$. 
\item The dependence model is a $t$-copula with a correlation matrix which has a 3-factor model structure. Two dependence models with slightly different correlation matrices and with degrees of freedom equal to 3 and 20 are used. 
\end{itemize}
Same as in Experiment~1, we have created four market models, each made up of a combination of marginal distributions and a dependence model. 

The inputs of Algorithms~\ref{alg:ecpalgo} and \ref{alg:accpalgo} in Experiment~2 are specified as follows. 
In Algorithm~\ref{alg:ecpalgo}, we set $\varepsilon=0.001,\tau=1,\delta=0.7$. 
We set $\underline{\phi}$ according to Remark~\ref{rmk:ecpalg}. 
When computing $\phi(f)$, we set $\underline{\phi}=0$, which is given by $0\ge-f$. 
When computing $\phi(-f)$, we set $\underline{\phi}=-\min\{\overline{\pi}_1,\ldots,\overline{\pi}_{50}\}$, which follows from $x_1\ge f,\ldots,x_{50}\ge f$. 
In Algorithm~\ref{alg:accpalgo}, we set $\varepsilon=0.001,\tau=1,\gamma=0.1,\zeta=0.8,\delta=0.7,\overline{c}=100,\overline{\BIy}=100\cdot\vecone$. 
$\underline{\phi}$ is set to be the same as in Algorithm~\ref{alg:ecpalgo} and $\overline{\phi}$ is set according to Remark~\ref{rmk:accpalg}. 
When computing $\phi(f)$, we set $\overline{\phi}=\min\{\overline{\pi}_1,\ldots,\overline{\pi}_{50}\}$ and $c^{\star(0)}+\langle\BIy^{\star(0)},\BIg\rangle=x_{i_{\min}}\ge f$, where $i_{\min}:=\argmin_{1\le i\le 50}\overline{\pi}_i$. 
When computing $\phi(-f)$, we set $\underline{\phi}=0$ and $c^{\star(0)}+\langle\BIy^{\star(0)},\BIg\rangle=0\ge-f$.

\begin{figure}[t]
\centering
\includegraphics[width=0.30\linewidth]{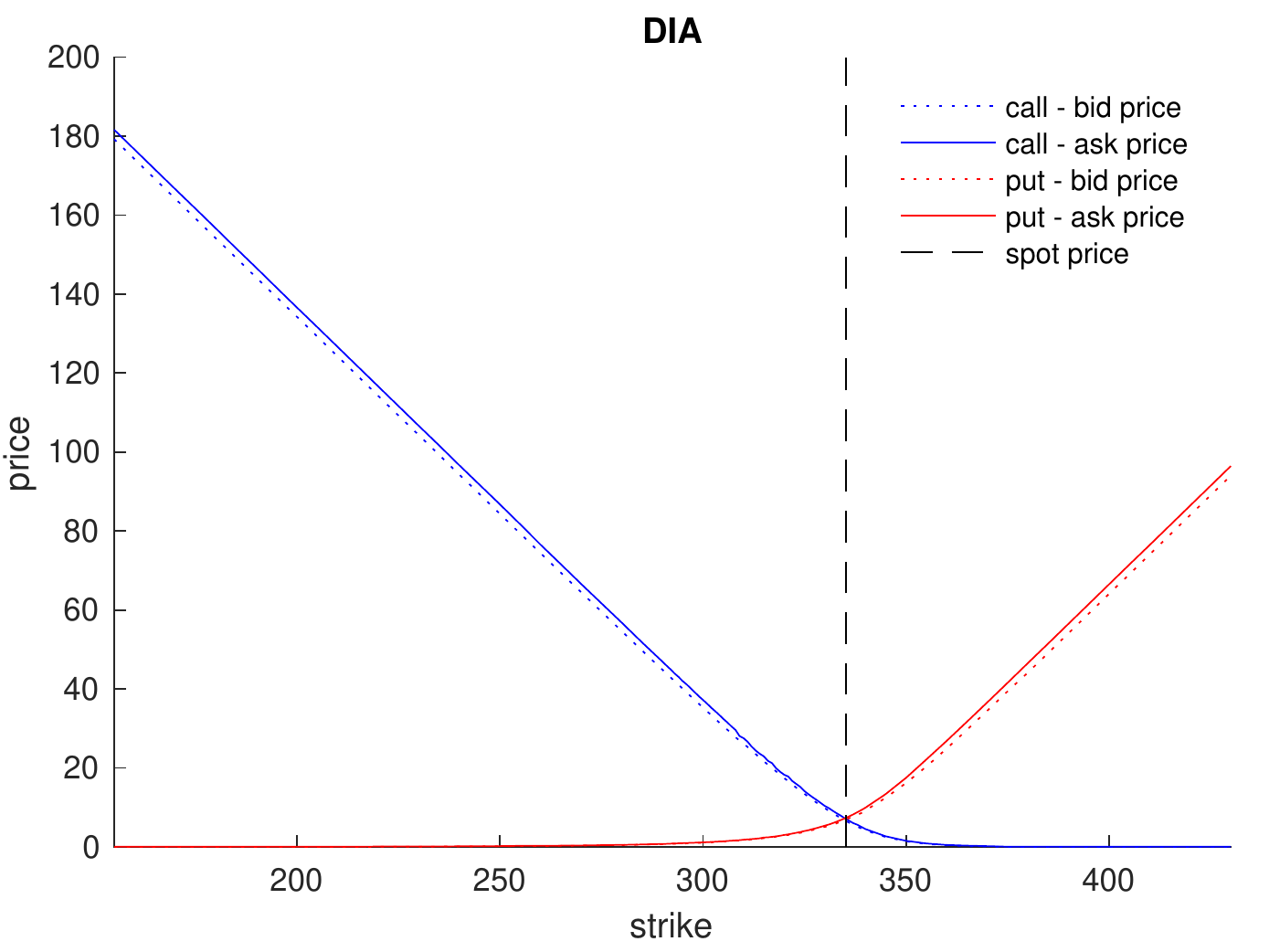}
~
\includegraphics[width=0.30\linewidth]{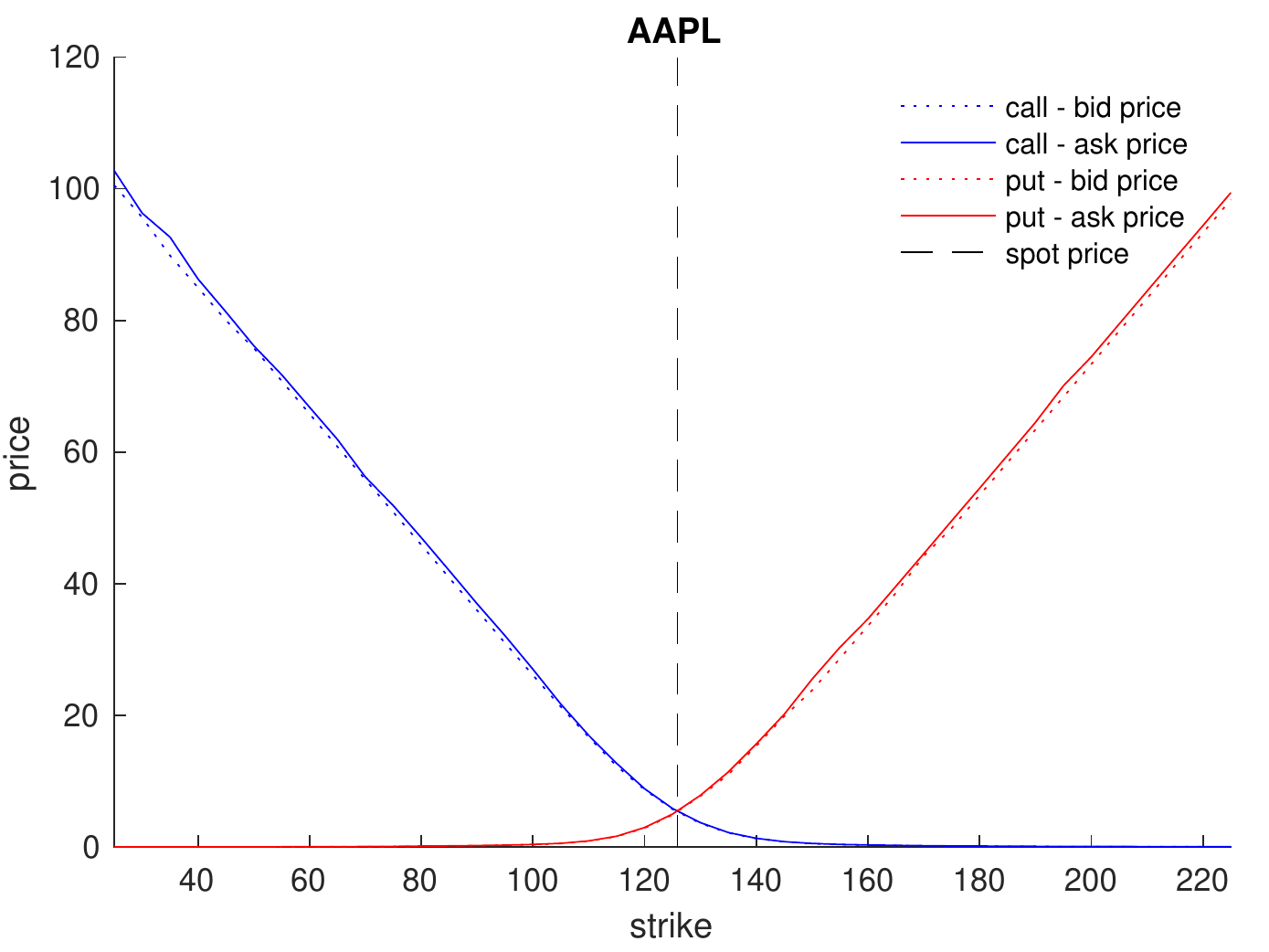}
~
\includegraphics[width=0.30\linewidth]{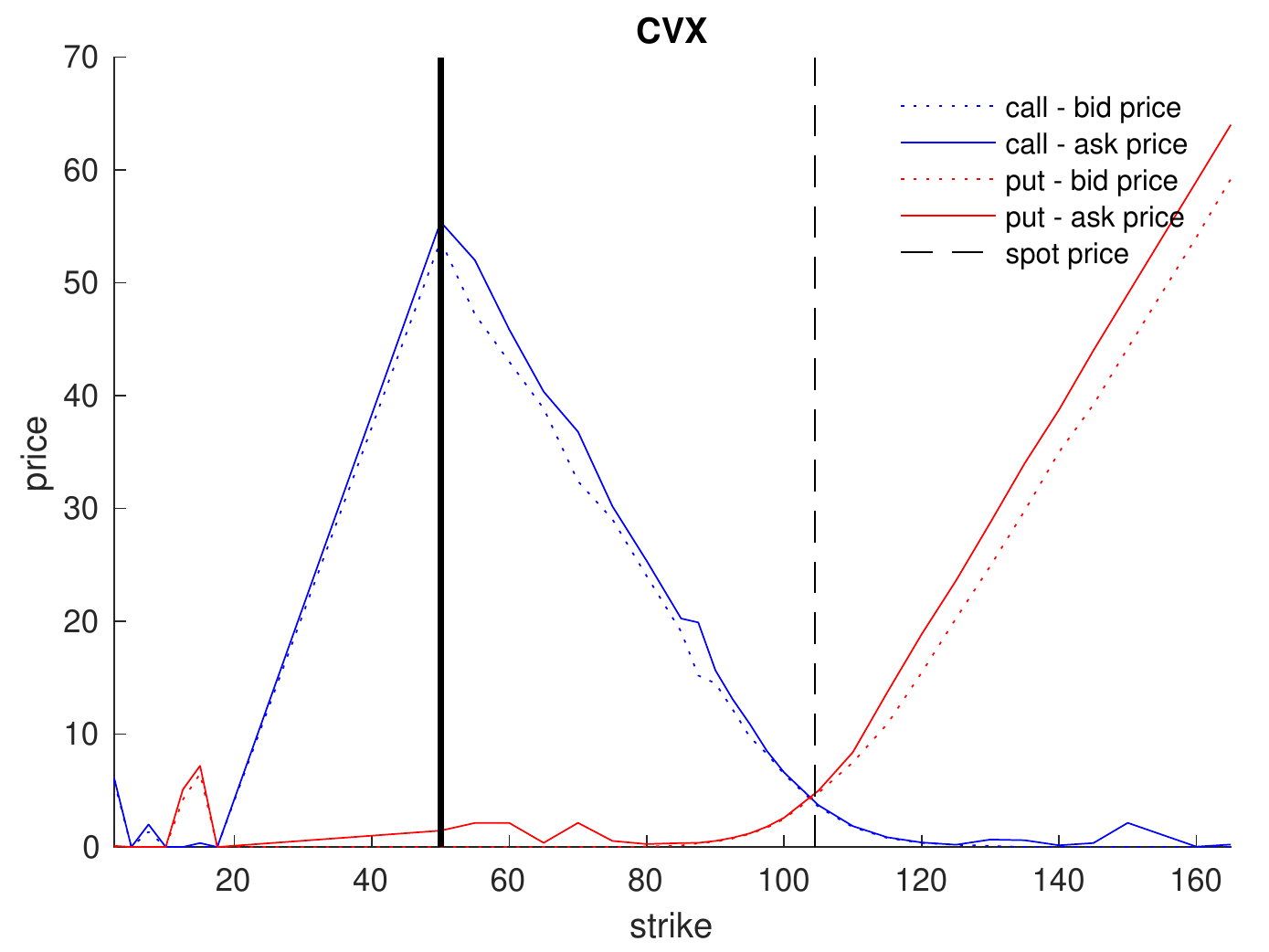}
\caption{Prices of options written on the SPDR Dow Jones Industrial Average ETF Trust (DIA), {Apple Inc.} (AAPL), and Chevron Corporation (CVX) used in Experiment~4. Prices of options written on CVX with strike prices below \$50 (indicated by the black solid vertical line) are anomalous and are removed from the experiment.}
\label{fig:pricedata}
\end{figure}

\subsection{Notes about Experiment~4}
\label{ecssec:exp4}

While preprocessing the price data, we noticed that the prices of options written on CVX with strike prices below \$50 are clearly anomalous (see the third panel of Figure~\ref{fig:pricedata}). Thus, we remove these options from Experiment~4. 

As mentioned in Section~\ref{ssec:exp4}, we adopt an arbitrage removal process to adjust the bid and ask prices of the options. 
The specific procedure is detailed as follows.
\begin{itemize}
\item We treat the DIA ETF as an asset and treat the options written on DIA as vanilla call and put options.
\item For each asset (including the DIA ETF), we repeat the following steps.
\begin{itemize}
\item Let $0<\kappa_1<\cdots<\kappa_m$ denote the strike prices of the call and put options written on this asset. For $j=1,\ldots,m$, let $\underline{\pi}^{\mathrm{call}}_j$ and $\overline{\pi}^{\mathrm{call}}_j$ denote the bid and ask prices of the call option with strike $\kappa_j$, respectively. Similarly, let $\underline{\pi}^{\mathrm{put}}_j$ and $\overline{\pi}^{\mathrm{put}}_j$ denote the bid and ask prices of the put option with strike $\kappa_j$, respectively. 
\item Fix $\overline{x}>\kappa_m$ to be large enough\footnote{In the MATLAB implementation we take $\overline{x}:=2\kappa_m$.}. Then, let $\big(\widehat{v}^{\mathrm{call}-}_j$, $\widehat{v}^{\mathrm{call}+}_j$, $\widehat{v}^{\mathrm{put}-}_j$, $\widehat{v}^{\mathrm{put}+}_j\big)_{j=1:m}$ be an approximate minimizer of the following minimization problem:
\begin{align}
\begin{split}
\min\Big\{&\textstyle\sum_{j=1}^m\big(v^{\mathrm{call}-}_j+v^{\mathrm{call}+}_j+v^{\mathrm{put}-}_j+v^{\mathrm{put}+}_j\big):\\
&\quad\exists\mu\in\CP([0,\overline{x}]),\;\{0,\kappa_1,\ldots,\kappa_m,\overline{x}\}\subseteq\mathrm{supp}(\mu),\\
&\quad\textstyle\underline{\pi}^{\mathrm{call}}_j-v^{\mathrm{call}-}_j\le\int_{[0,\overline{x}]}(x-\kappa_j)^+\,\mu(\Td x)\le\overline{\pi}^{\mathrm{call}}_j+v^{\mathrm{call}+}_j,\\
&\quad\textstyle\underline{\pi}^{\mathrm{put}}_j-v^{\mathrm{put}-}_j\le\int_{[0,\overline{x}]}(\kappa_j-x)^+\,\mu(\Td x)\le\overline{\pi}^{\mathrm{put}}_j+v^{\mathrm{put}+}_j,\\
&\quad v^{\mathrm{call}-}_j\ge0,v^{\mathrm{call}+}_j\ge0,v^{\mathrm{put}-}_j\ge0,v^{\mathrm{put}+}_j\ge0\;\forall 1\le j\le m\Big\}.
\end{split}
\label{eqn:arbitrage-removal}
\end{align}
\item For $j=1,\ldots,m$, adjust the bid and ask prices of the call option with strike $\kappa_j$ to $\underline{\pi}^{\mathrm{call}}_j-\widehat{v}^{\mathrm{call}-}_j$ and $\overline{\pi}^{\mathrm{call}}_j+\widehat{v}^{\mathrm{call}+}_j$, respectively. Adjust the bid and ask prices of the put option with strike $\kappa_j$ to $\underline{\pi}^{\mathrm{put}}_j-\widehat{v}^{\mathrm{put}-}_j$ and $\overline{\pi}^{\mathrm{put}}_j+\widehat{v}^{\mathrm{put}+}_j$, respectively. 
\end{itemize}
\end{itemize}

The following proposition shows that the above procedure produces arbitrage-free option prices. 
\begin{proposition}
Let $m\in\N$ and $0<\kappa_1<\cdots<\kappa_m<\overline{x}$. 
For $j=1,\ldots,m$, let $\underline{\pi}^{\CALL}_j<\overline{\pi}^{\CALL}_j$, $\underline{\pi}^{\PUT}_j<\overline{\pi}^{\PUT}_j$ be fixed, and define
\begin{align*}
\pi\big(y^{\CALL}_1,\ldots,y^{\CALL}_m,y^{\PUT}_1,\ldots,y^{\PUT}_m\big):=\sum_{j=1}^m&(y^{\CALL}_j \vee 0)\overline{\pi}^{\CALL}_j-(-y^{\CALL}_j \vee 0)\underline{\pi}^{\CALL}_j\\
&+(y^{\PUT}_j \vee 0)\overline{\pi}^{\PUT}_j-(-y^{\PUT}_j \vee 0)\underline{\pi}^{\PUT}_j
\end{align*}
for all $\big(y^{\CALL}_j\in\R\big)_{j=1:m},\big(y^{\PUT}_j\in\R\big)_{j=1:m}$. 
Suppose there exists $\mu\in\CP([0,\overline{x}])$ with $\{0,\kappa_1,\ldots,\kappa_m,\overline{x}\}\subseteq\mathrm{supp}(\mu)$, such that for $j=1,\ldots,m$,
\begin{align}
\begin{split}
&\underline{\pi}^{\mathrm{call}}_j\le\int_{[0,\overline{x}]}(x-\kappa_j)^+\,\mu(\Td x)\le\overline{\pi}^{\mathrm{call}}_j,\\
&\underline{\pi}^{\mathrm{put}}_j\le\int_{[0,\overline{x}]}(\kappa_j-x)^+\,\mu(\Td x)\le\overline{\pi}^{\mathrm{put}}_j.
\end{split}
\label{eqn:no-arbitrage-measure}
\end{align}
Then, for any $y^{\CALL}_1,\ldots,y^{\CALL}_m,y^{\PUT}_1,\ldots,y^{\PUT}_m$ such that 
\begin{align}
\left[\sum_{j=1}^m y^{\CALL}_j(x-\kappa_j)^++y^{\PUT}_j(\kappa_j-x)^+\right]-\pi\big(y^{\CALL}_1,\ldots,y^{\CALL}_m,y^{\PUT}_1,\ldots,y^{\PUT}_m\big)\ge0\quad \forall x\in\R_+,
\label{eqn:no-arbitrage-nonneg}
\end{align}
it holds that 
\begin{align}
\left[\sum_{j=1}^m y^{\CALL}_j(x-\kappa_j)^++y^{\PUT}_j(\kappa_j-x)^+\right]-\pi\big(y^{\CALL}_1,\ldots,y^{\CALL}_m,y^{\PUT}_1,\ldots,y^{\PUT}_m\big)=0\quad \forall x\in\R_+.
\label{eqn:no-arbitrage-zero}
\end{align}
\label{prop:no-arbitrage}
\end{proposition}

\proof{Proof of Proposition~\ref{prop:no-arbitrage}}
See Section~\ref{ecsec:proofarbitrageremoval}. 
\endproof

In order to numerically solve the minimization problem~(\ref{eqn:arbitrage-removal}), we need the following proposition. 
\begin{proposition}
Let $m\in\N$ and $0<\kappa_1<\cdots<\kappa_m<\overline{x}$. 
For $j=1,\ldots,J$ where $J\in\N$, let $g_j:[0,\overline{x}]\to\R$ be a continuous function which is piece-wise affine on the intervals $[0,\kappa_1]$, $[\kappa_1,\kappa_2]$, $\ldots$, $[\kappa_{m-1},\kappa_m]$, $[\kappa_m,\overline{x}]$. 
For $j=1,\ldots,J$, let $\pi_j\in\R$ be fixed. 
Suppose there exists $\mu\in\CP([0,\overline{x}])$ with $\{0,\kappa_1,\ldots,\kappa_m,\overline{x}\}\subseteq\mathrm{supp}(\mu)$, such that $\int_{[0,\overline{x}]}g_j\,\Td\mu=\pi_j$ for $j=1,\ldots,J$. 
Then, there exists a probability measure $\widehat{\mu}\in\CP([0,\overline{x}])$ with $\mathrm{supp}(\widehat{\mu})=\{0,\kappa_1,\ldots,\kappa_m,\overline{x}\}$ such that $\int_{[0,\overline{x}]}g_j\,\Td\widehat{\mu}=\pi_j$ for $j=1,\ldots,J$. 
\label{prop:arbitrage-discrete}
\end{proposition}

\proof{Proof of Proposition~\ref{prop:arbitrage-discrete}}
See Section~\ref{ecsec:proofarbitrageremoval}. 
\endproof

Notice that the payoff functions of the call and put options are all continuous and piece-wise affine on the intervals $[0,\kappa_1]$, $[\kappa_1,\kappa_2]$, $\ldots$, $[\kappa_{m-1},\kappa_m]$, $[\kappa_m,\overline{x}]$. Therefore, Proposition~\ref{prop:arbitrage-discrete} states that for fixed $\big(v^{\mathrm{call}-}_j\ge0,v^{\mathrm{call}+}_j\ge0,v^{\mathrm{put}-}_j\ge0,v^{\mathrm{put}+}_j\ge0\big)_{j=1:m}$, there exists $\mu\in\CP([0,\overline{x}])$ with $\{0,\kappa_1,\ldots,\kappa_m,\overline{x}\}\subseteq\mathrm{supp}(\mu)$, such that for $j=1,\ldots,J$,
\begin{align*}
\underline{\pi}^{\mathrm{call}}_j-v^{\mathrm{call}-}_j\le\int_{[0,\overline{x}]}(x-\kappa_j)^+\,\mu(\Td x)\le\overline{\pi}^{\mathrm{call}}_j+v^{\mathrm{call}+}_j,\\
\underline{\pi}^{\mathrm{put}}_j-v^{\mathrm{put}-}_j\le\int_{[0,\overline{x}]}(\kappa_j-x)^+\,\mu(\Td x)\le\overline{\pi}^{\mathrm{put}}_j+v^{\mathrm{put}+}_j,
\end{align*}
if and only if there exists a probability measure $\widehat{\mu}\in\CP([0,\overline{x}])$ with $\mathrm{supp}(\widehat{\mu})=\{0,\kappa_1,\ldots,\kappa_m,\overline{x}\}$ such that for $j=1,\ldots,J$,
\begin{align*}
\underline{\pi}^{\mathrm{call}}_j-v^{\mathrm{call}-}_j\le\int_{[0,\overline{x}]}(x-\kappa_j)^+\,\widehat{\mu}(\Td x)\le\overline{\pi}^{\mathrm{call}}_j+v^{\mathrm{call}+}_j,\\
\underline{\pi}^{\mathrm{put}}_j-v^{\mathrm{put}-}_j\le\int_{[0,\overline{x}]}(\kappa_j-x)^+\,\widehat{\mu}(\Td x)\le\overline{\pi}^{\mathrm{put}}_j+v^{\mathrm{put}+}_j.
\end{align*}
Thus, the minimization problem~(\ref{eqn:arbitrage-removal}) can be reformulated as follows:
\begin{align}
\begin{split}
\min\Big\{&\textstyle\sum_{j=1}^m\big(v^{\mathrm{call}-}_j+v^{\mathrm{call}+}_j+v^{\mathrm{put}-}_j+v^{\mathrm{put}+}_j\big):\\
&\quad\exists p_0>0,p_1>0,\ldots,p_m>0,p_{m+1}>0,\;\textstyle\sum_{i=0}^{m+1}p_i=1,\\
&\quad\textstyle\underline{\pi}^{\mathrm{call}}_j-v^{\mathrm{call}-}_j\le\left[\sum_{i=1}^m p_i(\kappa_i-\kappa_j)^+\right]+p_{m+1}(\overline{x}-\kappa_j)\le\overline{\pi}^{\mathrm{call}}_j+v^{\mathrm{call}+}_j,\\
&\quad\textstyle\underline{\pi}^{\mathrm{put}}_j-v^{\mathrm{put}-}_j\le p_0\kappa_j+\left[\sum_{i=1}^m p_i(\kappa_j-\kappa_i)^+\right]\le\overline{\pi}^{\mathrm{put}}_j+v^{\mathrm{put}+}_j,\\
&\quad v^{\mathrm{call}-}_j\ge0,v^{\mathrm{call}+}_j\ge0,v^{\mathrm{put}-}_j\ge0,v^{\mathrm{put}+}_j\ge0\;\forall 1\le j\le m\Big\}.
\end{split}
\label{eqn:arbitrage-removal-reform}
\end{align}
Now, problem~(\ref{eqn:arbitrage-removal-reform}) can be formulated into a linear programming problem (the constraints $p_0>0$, $p_1>0$, $\ldots,p_m>0$, $p_{m+1}>0$ need to be approximated by $p_0\ge\eta$, $p_1\ge\eta$, $\ldots$, $p_m\ge\eta$, $p_{m+1}\ge\eta$ for $\eta>0$ very small) which can be solved efficiently. 


\section{Proofs of the Fundamental Theorem and the superhedging duality}
\label{sec:proofduality}

This section is dedicated to the proof of the fundamental theorem (Theorem~\ref{thm:equivalence}) and the superhedging duality (Theorem~\ref{thm:duality}).
Let us first introduce the following notions and definitions.
For any vector $\BIx$ and any $i<j$, let $[\BIx]_{i:j}$ denote the vector formed with the $i$-th through the $j$-th entries of $\BIx$, \textit{i.e.} $[\BIx]_{i:j}:=(x_i,x_{i+1},\ldots,x_{j})^{\mathsf{T}}$. 
For a subset $A$ of a Euclidean space, let $\cone(A)$ denote the conic hull of $A$. 
Let $\Delta$ denote the following polytope,
\begin{align}
\Delta:=[\underline{\pi}_1,\overline{\pi}_1]\times\dots\times[\underline{\pi}_m,\overline{\pi}_m]\subset\R^m.
\label{eqn:Deltadef}
\end{align}
One verifies from definitions \eqref{eqn:Deltadef} and \eqref{eqn:pricepi} that, for every $\BIy=(y_1,\ldots,y_m)^{\mathsf T}\in\R^m$
\begin{align*}
\pi(\BIy)=\max_{\BIq\in\Delta}\langle\BIy,\BIq\rangle,
\end{align*}
thus $\pi(\cdot)$ is sublinear.
Let us define the sets $\widetilde{\CP}\subseteq\CP(\Omega)$ and $\Gamma\subseteq\R^m$ as follows:
\begin{align*}
\widetilde{\CP}:=&\left\{\mu\in\CP(\Omega):\int_{\Omega}|f|+\textstyle{\sum_{j=1}^m|g_j|}\Td\mu<\infty\right\},
\\
\Gamma:=&\left\{\int_{\Omega}\BIg \Td\mu:\mu\in\widetilde{\CP}\right\},
\end{align*}
where $\int_{\Omega}\BIg \Td\mu$ denotes $\left(\int_{\Omega} g_1 \Td\mu,\ldots,\int_{\Omega} g_m \Td\mu\right)^{\mathsf T}\in\R^m$.
Moreover, let us denote the set of payoff functions that can be super-replicated with no initial cost by
\begin{align}
\CC:=\left\{\langle\BIy,\BIg\rangle-\pi(\BIy):\BIy\in\R^m\right\}-\CL^0_+,
\label{eqn:lemma1Cdef}
\end{align}
where $\CL^0_+:=\big\{h:\Omega\to\R:h\text{ is Borel measurable}, h(\omega)\ge 0\;\forall \omega\in\Omega\big\}$ denotes the set of all Borel measurable non-negative functions.

Let us first prove the following two lemmata which are intermediate results. 

\begin{lemma}
Under Assumption~\ref{asp:na}, if $(f_n)_{n\ge1}\subset\CC$ and $f_n\to f$ point-wise, then $f\in\CC$. 
\label{lem:limitclosure}
\end{lemma}

\proof{Proof of Lemma~\ref{lem:limitclosure} (adapted from Theorem~2.2 of \cite{bouchard2015arbitrage})}
Suppose $(f_n)_{n\ge 1}\subset\CC$ and $f_n\to f$ point-wise. 
Then, there exist $\BIy_n\in\R^m,h_n\in\CL^0_+$ such that $f_n=\langle\BIy_n,\BIg\rangle-\pi(\BIy_n)-h_n$, for $n\ge 1$. 
The goal is to find a subsequence $(\BIy_{n_k})_{k\ge1}$ such that $\lim_{k\to\infty}\langle\BIy_{n_k},\BIg\rangle-\pi(\BIy_{n_k})=\langle\BIy,\BIg\rangle-\pi(\BIy)$ point-wise for some $\BIy\in\R^m$. 
Given such a subsequence, for any $\omega\in\Omega$,
\begin{align*}
\langle\BIy,\BIg(\omega)\rangle-\pi(\BIy)-f(\omega)=&\lim_{k\to\infty}\langle\BIy_{n_k},\BIg(\omega)\rangle-\pi(\BIy_{n_k})-f_{n_k}(\omega)=\lim_{k\to\infty}h_{n_k}(\omega)
\ge0,
\end{align*}
hence $f\in\CC$. 

Let us consider the following two separate cases. 

\emph{Case 1: $\liminf_{n\to\infty}\|\BIy_n\|<\infty$.} 
In this case, there exists a subsequence $(\BIy_{n_k})_{k\ge1}$ that converges to some $\BIy\in\R^m$ by compactness, and $\lim_{k\to\infty}\langle\BIy_{n_k},\BIg\rangle-\pi(\BIy_{n_k})=\langle\BIy,\BIg\rangle-\pi(\BIy)$ follows from the continuity of $\pi(\cdot)$. 

\emph{Case 2: $\liminf_{n\to\infty}\|\BIy_n\|=\infty$.} 
Define $\widetilde{\BIy}_{n}:=\frac{\BIy_n}{1+\|\BIy_n\|}$. 
Then $\lim_{n\to\infty}\|\widetilde{\BIy}_n\|=1$ and there exists a subsequence of $(\widetilde{\BIy}_n)_{n\ge1}$, denoted again by $(\widetilde{\BIy}_n)_{n\ge1}$ for simplicity, such that $\widetilde{\BIy}_n\to\widetilde{\BIy}\in\R^m$ as $n\to\infty$. 
Moreover, for all $\omega\in\Omega$, 
\begin{align*}
\langle\widetilde{\BIy},\BIg(\omega)\rangle-\pi(\widetilde{\BIy})=&\lim_{n\to\infty}\langle\widetilde{\BIy}_n,\BIg(\omega)\rangle-\pi(\widetilde{\BIy}_n)\\
=& \lim_{n\to\infty}\frac{1}{1+\|\BIy_n\|}\left[\langle\BIy_n,\BIg(\omega)\rangle-\pi(\BIy_n)\right]\\
\ge& \lim_{n\to\infty}\frac{f_n(\omega)}{1+\|\BIy_n\|}\\
=& \lim_{n\to\infty}\frac{f_n(\omega)-f(\omega)}{1+\|\BIy_n\|}+\frac{f(\omega)}{1+\|\BIy_n\|}
=0,
\end{align*}
where the second equality follows from the positive homogeneity of $\pi(\cdot)$.
By Assumption~\ref{asp:na}, we have 
\begin{align}
\langle\widetilde{\BIy},\BIg\rangle-\pi(\widetilde{\BIy})=0.
\label{eqn:lemma1case2}
\end{align}

The rest of the proof follows from induction, with the following induction hypothesis:
\begin{align}
\begin{split}
&[\BIy_n]_{j}=0\;\forall j\in I,n\ge 1,\text{ where }I\subset\{1,\ldots,m\} \text{ and }|I|=l\\[0.8em]
\Longrightarrow\qquad & (\BIy_n)_{n\ge1}\text{ admits a subsequence }(\BIy_{n_k})_{k\ge1}\text{ such that }\exists\BIy\in\R^m \text{ satisfying} \\
&\lim_{k\to\infty}\langle\BIy_{n_k},\BIg\rangle-\pi(\BIy_{n_k})=\langle\BIy,\BIg\rangle-\pi(\BIy).
\end{split}
\label{eqn:lemma1inductionhypo}
\end{align}
The claim in (\ref{eqn:lemma1inductionhypo}) holds trivially when $l=m$. The goal is to show that (\ref{eqn:lemma1inductionhypo}) holds when $l=0$. 

Now, let $\widetilde{I}\subset\{1,\ldots,m\}$, $|\widetilde{I}|=l-1$ for some $1\le l\le m$, and suppose that $[\BIy_n]_j=0$ for all $j\in\widetilde{I}$ and ${n\ge1}$. 
If $\liminf_{n\to\infty}\|\BIy_n\|<\infty$, then the claim holds by Case~1 above. 
Therefore, let us assume that $\liminf_{n\to\infty}\|\BIy_n\|=\infty$. 
Let us first restrict $(\BIy_n)_{n\ge1}$ to a subsequence $(\BIy_{n_k})_{k\ge1}$ such that $\mathrm{sign}(\BIy_{n_k})=\mathrm{sign}(\BIy_{n_1})$ for all $k\ge1$, where for $\BIx\in\R^m$, $\mathrm{sign}(\BIx):=(\mathrm{sign}([\BIx]_1),\ldots,\mathrm{sign}([\BIx]_m))^{\mathsf T}\in\{-1,1\}^m$ and for $x\in\R$,
\begin{align*}
\mathrm{sign}(x):=\begin{cases}
1, & \text{if }x\ge 0,\\
-1, & \text{otherwise}.
\end{cases}
\end{align*}
Such a subsequence exists since $m$ is finite. 
Let this subsequence be denoted by $(\BIy_n)_{n\ge 1}$ again. 
By Case~2 and by setting $\widetilde{\BIy}_{n}:=\frac{\BIy_n}{1+\|\BIy_n\|}$, we get a subsequence (still denoted by $(\widetilde{\BIy}_n)_{n\ge1}$) such that $\widetilde{\BIy}_n\to\widetilde{\BIy}$ as $n\to\infty$. 
Note that we have $\mathrm{sign}(\BIy_{n})=\mathrm{sign}(\widetilde{\BIy}_{n})$ for all $n$. 
Moreover, (\ref{eqn:lemma1case2}) holds and for all $n\ge1$
\begin{align}
\mathrm{sign}([\BIy_{n}]_j)=\mathrm{sign}([\widetilde{\BIy}_{n}]_j)=\mathrm{sign}([\widetilde{\BIy}]_j)\quad\forall j\in\{1,\ldots,m\}\text{ such that }[\widetilde{\BIy}]_j\ne0.
\label{eqn:lemma1sign}
\end{align}
For every $n\ge1$, let $\lambda_n=\min\{\gamma_{n,1},\ldots,\gamma_{n,m}\}$, where
\begin{align*}
\gamma_{n,j}=\begin{cases}
\frac{[\BIy_n]_j}{[\widetilde{\BIy}]_j}, &\text{if }[\widetilde{\BIy}]_j\ne0,\\
+\infty, &\text{otherwise},
\end{cases}
\end{align*}
for $j=1,\ldots,m$. 
By (\ref{eqn:lemma1sign}), $\gamma_{n,j}\ge0$ for all $n$ and $j$. 
We also have $\lambda_n<\infty$ for all $n\ge1$ since $\|\widetilde{\BIy}\|=1$. 
Let $\bar{\BIy}_n:=\BIy_n-\lambda_n\widetilde{\BIy}$, and observe that for all $j\in\{1,\ldots,m\}\setminus\widetilde{I}$,
\begin{align}
[\bar{\BIy}_n]_j=[\BIy_n]_j-\lambda_n[\widetilde{\BIy}]_j\begin{cases}
\ge0, & \text{if }[\BIy_n]_j\ge0,\\
\le0, & \text{if }[\BIy_n]_j<0,\\
=0, & \text{if }\lambda_n=\gamma_{n,j}.
\end{cases}
\label{eqn:lemma1sign3}
\end{align}
It follows for all $n\ge1$ that
\begin{align}
\mathrm{sign}([\BIy_{n}]_j)=\mathrm{sign}([\bar{\BIy}_n]_j) \quad \forall j\in\{1,\ldots,m\} \text{ such that } [\bar{\BIy}_n]_j\ne0. 
\label{eqn:lemma1sign2}
\end{align}
Combining (\ref{eqn:lemma1sign}) and (\ref{eqn:lemma1sign2}), it follows that if $\pi(\BIy_n)=\langle\BIy_n,\BIq\rangle$ for $\BIq\in\Delta$, then $\pi(\bar{\BIy}_n)=\langle\bar{\BIy}_n,\BIq\rangle$ and $\pi(\widetilde{\BIy})=\langle\widetilde{\BIy},\BIq\rangle$. 
Hence, $\pi(\bar{\BIy}_n)=\pi(\BIy_n)-\lambda_n\pi(\widetilde{\BIy})$. 
Combining (\ref{eqn:lemma1sign3}) and the assumption that $[\bar{\BIy}_n]_j=0$ for all $j\in\widetilde{I}$, it holds that $\bar{\BIy}_n$ has at least $l$ zero-entries.
Again, by extracting a subsequence, we can assume without loss of generality that the additional zero-entries in $\bar{\BIy}_n$ all occur at the same position, that is, there exists $j'\in\{1,\ldots,m\}\setminus\widetilde{I}$ such that $[\bar{\BIy}_n]_{j'}=0$ for all $n\ge1$. 
Now since $|\widetilde{I}\cup \{j'\}|=l$, by the induction hypothesis (\ref{eqn:lemma1inductionhypo}), $(\bar{\BIy}_{n})_{n\ge1}$ admits a subsequence $(\bar{\BIy}_{n_k})_{k\ge1}$ such that $\exists\BIy\in\R^m$, $\lim_{k\to\infty}\langle\bar{\BIy}_{n_k},\BIg\rangle-\pi(\bar{\BIy}_{n_k})=\langle\BIy,\BIg\rangle-\pi(\BIy)$.
By (\ref{eqn:lemma1case2}), it holds for $n\ge1$ that
\begin{align*}
\langle\BIy_n,\BIg\rangle-\pi(\BIy_n)=&\langle\BIy_n,\BIg\rangle-\pi(\BIy_n)-\lambda_n(\langle\widetilde{\BIy},\BIg\rangle-\pi(\widetilde{\BIy}))=\langle\bar{\BIy}_n,\BIg\rangle-\pi(\bar{\BIy}_n),
\end{align*}
hence $\lim_{k\to\infty}\langle{\BIy}_{n_k},\BIg\rangle-\pi({\BIy}_{n_k})=\langle\BIy,\BIg\rangle-\pi(\BIy)$. 
We have thus proved the induction hypothesis for the case where $|I|=l-1$. 
By induction, the existence of a subsequence $(\BIy_{n_k})_{k\ge1}$ such that $\lim_{k\to\infty}\langle\BIy_{n_k},\BIg\rangle-\pi(\BIy_{n_k})=\langle\BIy,\BIg\rangle-\pi(\BIy)$ for some $\BIy\in\R^m$ is established. 
The proof is now complete. 
\endproof

\begin{lemma}
For any $\nu\in\CP(\Omega)$, let $\widetilde{\CP}_{\nu}:=\left\{\mu\in\CP(\Omega):\nu\ll\mu,\int_{\Omega}|f|+\sum_{j=1}^m|g_j|\Td\mu<\infty\right\}$, and let $\Gamma_{\nu}:=\left\{\int_{\Omega}\BIg \Td\mu:\mu\in\widetilde{\CP}_{\nu}\right\}\subset\R^m$. 
Under Assumption~\ref{asp:na}, $\Delta\cap\relint(\Gamma_{\nu})\ne\emptyset$. 
In particular, $\Delta\cap\relint(\Gamma)\ne\emptyset$. 
\label{lem:relint}
\end{lemma}

\proof{Proof of Lemma~\ref{lem:relint} (partially adapted from Lemma~3.3 of \cite{bouchard2015arbitrage})}
Fix a $\nu\in\CP(\Omega)$. 
Suppose, for the sake of contradiction, that $\Delta\cap\relint(\Gamma_{\nu})=\emptyset$. 
Since both $\Delta$ and $\Gamma_\nu$ are convex, and $\Delta$ is polyhedral, by Theorem~20.2 of \citet{rockafellar1970convex}, there exists a hyperplane (namely, there exist $\BIy\in\R^n$ with $\BIy\ne\veczero$ and $\alpha\in\R$)
\begin{align*}
H:=\{\BIw\in\R^n:\langle\BIy,\BIw\rangle=\alpha\}
\end{align*}
that separates $\Delta$ and $\Gamma_\nu$ properly and that $\Gamma_\nu\nsubseteq H$. 
Suppose without loss of generality that $\Delta$ is contained in the half-space $H_-:=\{\BIw:\langle\BIy,\BIw\rangle\le\alpha\}$, then
\begin{align*}
\pi(\BIy)=\sup_{\BIq\in\Delta}\langle\BIy,\BIq\rangle\le\alpha.
\end{align*}
On the other hand, 
\begin{align}
\int_{\Omega}\langle\BIy,\BIg\rangle \Td\mu\ge \alpha\ge\pi(\BIy)\quad\forall \mu\in\widetilde{\CP}_{\nu}.
\label{eqn:lemma2contradiction}
\end{align}
This implies that $\int_{\Omega}\langle\BIy,\BIg\rangle-\pi(\BIy) \Td\mu\ge0\;\forall\mu\in\widetilde{\CP}_{\nu}$. 
Suppose, for the sake of contradiction, that there exists $\omega\in \Omega$ such that $\langle\BIy,\BIg(\omega)\rangle-\pi(\BIy)=\beta<0$. 
Let $\widetilde{\nu}=(\nu+\delta_\omega)/2$, where $\delta_\omega$ is the Dirac measure at $\omega$. 
By Theorem~VII.57 of \citet{dellacherie1982probabilities}, there exists $\widehat{\nu}\in\CP(\Omega)$ such that $\widehat{\nu}\sim\widetilde{\nu}$ and $\int_{\Omega}|f|+\sum_{j=1}^m|g_j|\Td\widehat{\nu}<\infty$. 
Therefore, $\widehat{\nu}\in\widetilde{\CP}_{\nu}$ and $\widehat{\nu}(\{\omega\})>0$. 
Now for $\varepsilon>0$, define $\eta_{\varepsilon}$ by
\begin{align*}
\frac{\Td\eta_{\varepsilon}}{\Td\widehat{\nu}}=\frac{\INDI_{\{\omega\}}+\varepsilon}{\widehat{\nu}(\{\omega\})+\varepsilon}.
\end{align*}
It is clear that $\eta_{\varepsilon}\sim\widehat{\nu}$ and $\eta_{\varepsilon}\in\widetilde{\CP}_{\nu}$ for all $\varepsilon>0$. 
Moreover, 
\begin{align*}
\int_{\Omega}\langle\BIy,\BIg\rangle-\pi(\BIy) \Td\eta_{\varepsilon}=\frac{\beta\widehat{\nu}(\{\omega\})+\varepsilon\int_{\Omega}\langle\BIy,\BIg\rangle-\pi(\BIy) \Td\widehat{\nu}}{\widehat{\nu}(\{\omega\})+\varepsilon}\to \beta<0
\end{align*}
when $\varepsilon\to0$. 
This implies that there exist $\varepsilon>0$ and $\eta_{\varepsilon}\in\widetilde{\CP}_{\nu}$ such that $\int_{\Omega}\langle\BIy,\BIg\rangle-\pi(\BIy) \Td\eta_{\varepsilon}<0$, which is a contradiction to (\ref{eqn:lemma2contradiction}). 
Therefore, $\langle\BIy,\BIg(\omega)\rangle-\pi(\BIy)\ge0$ for all $\omega\in\Omega$ and it follows from Assumption~\ref{asp:na} that $\langle\BIy,\BIg\rangle-\pi(\BIy)=0$. 
Hence, by (\ref{eqn:lemma2contradiction}), we have for all $\mu\in\widetilde{\CP}_{\nu}$ that
\begin{align*}
\alpha\le\int_{\Omega}\langle\BIy,\BIg\rangle \Td\mu=\pi(\BIy)\le\alpha,
\end{align*}
and we conclude that $\int_{\Omega}\langle\BIy,\BIg\rangle \Td\mu=\alpha$ for all $\mu\in\widetilde{\CP}_{\nu}$ and thus $\Gamma_\nu\subseteq H$, which is a contradiction to $\Gamma_\nu\nsubseteq H$. 
The last statement of Lemma~\ref{lem:relint} follows from the fact that $\widetilde{\CP}=\bigcup_{\nu\in\CP(\Omega)}\widetilde{\CP}_\nu$, and that $\Gamma=\bigcup_{\nu\in\CP(\Omega)}\Gamma_\nu$. 
The proof is now complete. 
\endproof

\proof{Proof of Theorem~\ref{thm:equivalence} (adapted from Theorem~3.1 of \cite{bouchard2015arbitrage})}
Suppose that Assumption~\ref{asp:na} holds. 
We have by Lemma~\ref{lem:relint} that for any $\nu\in\CP(\Omega)$, there exists $\mu\in\CP(\Omega)$ such that $\nu\ll\mu$ and $\int_{\Omega}\BIg \Td\mu\in\Delta$, hence (ii) holds. 

Conversely, suppose that Assumption~\ref{asp:na} does not hold. 
Then, there exists $\BIy\in\R^m$ and $\omega\in\Omega$ such that $\langle\BIy,\BIg\rangle-\pi(\BIy)\ge 0$ and $\langle\BIy,\BIg(\omega)\rangle-\pi(\BIy)>0$.
With this $\BIy$ and any $\mu\in\CQ$, it holds that
\begin{align*}
\int_{\Omega}\langle\BIy,\BIg\rangle \Td\mu\le \sup_{\BIq\in\Delta}\langle\BIy,\BIq\rangle=\pi(\BIy),
\end{align*}
implying $\int_{\Omega}\langle\BIy,\BIg\rangle-\pi(\BIy) \Td\mu=0$.
Hence, $\mu(\{\omega\})=0$ and $\delta_\omega\centernot{\ll}\mu$ for all $\mu\in\CQ$, implying that (ii) does not hold. 
The proof is now complete.  
\endproof

\proof{Proof of Theorem~\ref{thm:duality} (partially adapted from Theorem~3.4 of  \cite{bouchard2015arbitrage})}

Let us first prove statement~(i). 
Suppose that $\phi(f)=-\infty$, then for all $n\ge 1$, there exists $\BIy_n\in\R^m$ such that $-n+\langle\BIy_n,\BIg\rangle-\pi(\BIy_n)\ge f$, which implies that 
\begin{align*}
\langle\BIy_n,\BIg\rangle-\pi(\BIy_n)\ge f+n\ge \min\{f+n,1\},
\end{align*}
and thus $\min\{f+n,1\}\in\CC$ where $\CC$ is the set of super-replicable payoff functions defined in (\ref{eqn:lemma1Cdef}). 
By Lemma~\ref{lem:limitclosure}, $\lim_{n\to\infty}\min\{f+n,1\}=1\in\CC$, which is a contradiction to Assumption~\ref{asp:na}. 

Next, let us prove statement~(ii) by showing the existence of $\BIy\in\R^m$ such that $\phi(f)+\langle\BIy,\BIg\rangle-\pi(\BIy)\ge f$. If $\phi(f)=+\infty$, the claim holds trivially. 
If $\phi(f)$ is finite, then for all $n\ge 1$, there exists $\BIy_n\in\R^m$ such that $\phi(f)+\frac{1}{n}+\langle\BIy_n,\BIg\rangle-\pi(\BIy_n)\ge f$, which implies that $f-\phi(f)-\frac{1}{n}\in\CC$. 
By Lemma~\ref{lem:limitclosure}, $f-\phi(f)=\lim_{n\to\infty}f-\phi(f)-\frac{1}{n}\in\CC$, hence there exists $\BIy\in\R^m$ such that $\phi(f)+\langle\BIy,\BIg\rangle-\pi(\BIy)\ge f$. 

Finally, let us prove statement~(iii), that is, the duality~(\ref{eqn:duality}). 
Assume first that $f$ is bounded from above. 
This part of the proof is based on the semi-infinite generalization of Fenchel's Duality Theorem by 
\citet{borwein1992partially}. 
Let $\CM(\Omega)$ denote the set of finite signed Borel measures on $\Omega$. 
Let 
\begin{align*}
\widetilde{\CM}:=&\bigg\{\mu\in\CM(\Omega):\int_{\Omega} \Big(\textstyle{\sum_{j=1}^m}|g_j|\Big)\Td|\mu|<\infty\bigg\},\\
\widetilde{\CM}^*:=&\bigg\{h:\Omega\to\R \text{ is Borel measurable}:|h|\le \alpha\Big(\textstyle{\sum_{j=1}^m}|g_j|+1\Big) \text{ for some }\alpha>0\bigg\},
\end{align*}
where $|\mu|$ denotes the total variation of the signed measure $\mu$. 
For $\mu\in\widetilde{\CM}$ and $h\in\widetilde{\CM}^*$, define
\begin{align*}
\langle h,\mu\rangle:=\int_{\Omega} h\Td\mu<\infty.
\end{align*}
$\widetilde{\CM}$ is a locally convex topological vector space equipped with the $\sigma(\widetilde{\CM},\widetilde{\CM}^*)$ topology. 
Let $\widetilde{\CM}_+:=\{\mu\in\widetilde{\CM}:\mu\text{ is a positive measure}\}$. 
For $\mu\in\widetilde{\CM}$, let  
\begin{align}
F(\mu):=\begin{cases}
\int_{\Omega} (-f)\Td\mu, & \text{if }\mu\in\widetilde{\CM}_+,\\
\infty, & \text{if }\mu\notin\widetilde{\CM}_+,
\end{cases}
\label{eqn:dualityFdef}
\end{align}
which is convex and proper since $-f$ is bounded from below. 
For $h\in\widetilde{\CM}^*$, let
\begin{align}
\begin{split}
F^*(h):=&\sup_{\mu\in\widetilde{\CM}}\left\{\langle h,\mu\rangle-F(\mu)\right\}\\
=&\sup_{\mu\in\widetilde{\CM}_+}\left\{\int_{\Omega} h+f \Td\mu\right\}\\
=&\begin{cases}
0 & \text{if }-h(\omega)\ge f(\omega)\;\forall\omega\in\Omega,\\
\infty & \text{if }\exists\omega\in\Omega,-h(\omega)<f(\omega),
\end{cases}
\end{split}
\label{eqn:dualityFconj}
\end{align}
where the last equality holds because $\gamma\delta_{\omega}\in\widetilde{\CM}_+$ for all $\omega\in\Omega$ and $\gamma>0$. 
Let $G:\widetilde{\CM}\to\R^{m+1}$ be defined as follows,
\begin{align}
G(\mu):=\left(\int_{\Omega} 1\Td\mu,\int_{\Omega} g_1\Td\mu,\ldots,\int_{\Omega} g_m\Td\mu\right)^{\mathsf T}.
\label{eqn:dualityGdef}
\end{align}
$G$ is linear and $\sigma(\widetilde{\CM},\widetilde{\CM}^*)$-continuous by definition. 
Let $G^{\mathsf T}$ be the adjoint of $G$, defined as
\begin{align*}
G^{\mathsf T}:\R^{m+1}&\to \widetilde{\CM}^*\\
\Blambda&\mapsto G^{\mathsf T}\Blambda=[\Blambda]_1+\langle[\Blambda]_{2:m+1},\BIg\rangle.
\end{align*}
Let us denote $[\Blambda]_1\in\R$ by $c$ and denote $[\Blambda]_{2:m+1}\in\R^m$ by $\BIy$. 
Thus, 
\begin{align}
G^{\mathsf T}\Blambda=c+\langle\BIy,\BIg\rangle.
\label{eqn:dualityadjoint}
\end{align}
Let $\widehat{\Delta}:=\{1\}\times\Delta\subset\R^{m+1}$. 
Let $\chi_{\widehat{\Delta}}$ be the characteristic function of $\widehat{\Delta}$, that is,
\begin{align}
\chi_{\widehat{\Delta}}(\BIq):=\begin{cases}
0, & \text{if }\BIq\in\widehat{\Delta},\\
\infty, & \text{if }\BIq\notin\widehat{\Delta}. 
\end{cases}
\label{eqn:dualitychidef}
\end{align}
Let $\chi_{\widehat{\Delta}}^*$ be the convex conjugate of $\chi_{\widehat{\Delta}}$, which is given by
\begin{align}
\begin{split}
\chi_{\widehat{\Delta}}^*(\Blambda)=&\sup_{\BIq\in\R^{m+1}}\left\{\langle\Blambda,\BIq\rangle-\chi_{\widehat{\Delta}}(\BIq)\right\}\\
=&[\Blambda]_1+\pi([\Blambda]_{2:m+1})\\
=&c+\pi(\BIy).
\end{split}
\label{eqn:dualitychiconj}
\end{align}
Let $\phi^D(f)$ denote the right-hand side of (\ref{eqn:duality}). 
By definitions (\ref{eqn:dualityFdef}), (\ref{eqn:dualityGdef}) and (\ref{eqn:dualitychidef}), $\phi^D(f)$ can be equivalently expressed as follows,
\begin{align*}
\begin{split}
\phi^D(f):=&\sup_{\mu\in\CQ}\int_{\Omega} f\Td\mu\\
=&-\inf_{\mu\in\widetilde{\CM}_+}\left\{\int_{\Omega} (-f)\Td\mu+\chi_{\widehat{\Delta}}(G(\mu))\right\}\\
=&-\inf_{\mu\in\widetilde{\CM}}\left\{F(\mu)+\chi_{\widehat{\Delta}}(G(\mu))\right\}.
\end{split}
\end{align*}
It holds that $\chi_{\widehat{\Delta}}$ is convex and proper in $\R^{m+1}$ and $F$ is convex and proper in $\widetilde{\CM}$. 
By Theorem~4.2 of \citet{borwein1992partially}, the following duality holds under the additional constraint qualification that $\relint(G(\dom F))\cap\dom\chi_{\widehat{\Delta}}\ne\emptyset$ (notice that $\chi_{\widehat{\Delta}}$ is a polyhedral function since its epigraph is a polyhedral set),
\begin{align*}
\begin{split}
-\phi^D(f)=&\inf_{\mu\in\widetilde{\CM}}\left\{F(\mu)+\chi_{\widehat{\Delta}}(G(\mu))\right\}\\
=&\sup_{\Blambda\in\R^{m+1}}\left\{-F^*(-G^{\mathsf T}\Blambda)-\chi_{\widehat{\Delta}}^*(\Blambda)\right\}\\
=&\sup\left\{-c-\pi(\BIy):c\in\R,\BIy\in\R^m,c+\langle\BIy,\BIg\rangle\ge f\right\}\\
=&-\inf\left\{c+\pi(\BIy):c\in\R,\BIy\in\R^m,c+\langle\BIy,\BIg\rangle\ge f\right\}
=-\phi(f),
\end{split}
\end{align*}
where the third equality follows from definitions (\ref{eqn:dualityFconj}), (\ref{eqn:dualityadjoint}) and (\ref{eqn:dualitychiconj}). 
Hence, we only need to show that the constraint qualification $\relint(G(\dom F))\cap\dom\chi_{\widehat{\Delta}}\ne\emptyset$ holds. 
It holds that
\begin{align*}
\begin{split}
\relint(G(\dom F))\cap\dom \chi_{\widehat{\Delta}}=&\relint(G(\dom F))\cap(\{1\}\times\R^m) \cap(\R\times\Delta)\\
=&\relint(G(\dom F)\cap(\{1\}\times\R^m)) \cap(\R\times\Delta).
\end{split}
\end{align*}
Since 
\begin{align*}
\dom F=&\left\{\mu\in\widetilde{\CM}_+:\int_{\Omega} (-f)\Td\mu<\infty\right\}\\
\supseteq&\left\{\mu\in\CM(\Omega):\mu\ge0,\int_{\Omega} |f|+\textstyle{\sum_{j=1}^m}|g_j|\Td\mu<\infty\right\},
\end{align*}
we have
\begin{align*}
G(\dom F)\cap(\{1\}\times\R^m) &\supseteq G\left(\left\{\mu\in\CM(\Omega):\mu\ge0,\mu(\Omega)=1,\int_{\Omega} |f|+\textstyle{\sum_{j=1}^m}|g_j|\Td\mu<\infty\right\}\right)\\
&= G(\widetilde{\CP})
= \{1\}\times \Gamma.
\end{align*}
Since $\relint(\{1\}\times {\Gamma}) \cap(\R\times\Delta)=\{1\}\times(\relint(\Gamma)\cap\Delta)$, which is non-empty by Lemma~\ref{lem:relint}, we have that $\relint(G(\dom F))\cap\dom \chi_{\widehat{\Delta}}\supseteq\relint(\{1\}\times {\Gamma}) \cap(\R\times\Delta)\ne\emptyset$ and the constraint qualification holds. 

Now suppose that $f$ is unbounded from above. 
Then it holds that for all $n\ge1$,
\begin{align}
\phi(\min\{f,n\})=\sup_{\mu\in\CQ}\int_{\Omega} \min\{f,n\}\Td\mu.
\label{eqn:thm2finalstep1}
\end{align}
On the one hand, 
\begin{align}
\begin{split}
\limsup_{n\to\infty}\sup_{\mu\in\CQ}\int_{\Omega} \min\{f,n\}\Td\mu\le&\sup_{\mu\in\CQ}\int_{\Omega} f\Td\mu.
\end{split}
\label{eqn:thm2finalstep2}
\end{align}
On the other hand, observe that $\limsup_{n\to\infty}\phi(\min\{f,n\})\le\phi(f)$. Let $\limsup_{n\to\infty}\phi(\min\{f,n\})=:\alpha\in\R\cup\{\infty\}$. If $\alpha=\infty$, then by (\ref{eqn:thm2finalstep2}) and (\ref{eqn:thm2finalstep1}),
\begin{align*}
\sup_{\mu\in\CQ}\int_{\Omega} f\Td\mu\ge\limsup_{n\to\infty}\sup_{\mu\in\CQ}\int_{\Omega} \min\{f,n\}\Td\mu=\limsup_{n\to\infty}\phi(\min\{f,n\})=\infty,
\end{align*}
and both sides of (\ref{eqn:duality}) are equal to $\infty$. 
Hence, we focus on the case where $\alpha<\infty$. 
Note that there exists $\BIy_n\in\R^m$ such that $\alpha+\langle\BIy_n,\BIg\rangle-\pi(\BIy_n)\ge\min\{f,n\}$ by statement~(ii). 
It follows that $\min\{f,n\}-\alpha\in\CC$ and hence $f-\alpha=\lim_{n\to\infty}\min\{f,n\}-\alpha\in\CC$ by Lemma~\ref{lem:limitclosure}. 
Therefore, there exists $\BIy\in\R^m$ such that $\alpha+\langle\BIy,\BIg\rangle-\pi(\BIy)\ge f$, and it hence holds that $\alpha\ge\phi(f)$.
Consequently, we have by (\ref{eqn:thm2finalstep2}) and (\ref{eqn:thm2finalstep1}) that
\begin{align*}
\phi(f)\le\alpha=\limsup_{n\to\infty}\phi(\min\{f,n\})=\limsup_{n\to\infty}\sup_{\mu\in\CQ}\int_{\Omega} \min\{f,n\}\Td\mu\le\sup_{\mu\in\CQ}\int_{\Omega} f\Td\mu.
\end{align*}
The reverse direction $\phi(f)\ge\sup_{\mu\in\CQ}\int_{\Omega} f\Td\mu$ is easier to verify. 
Suppose that $\phi(f)<\infty$, since otherwise $\phi(f)\ge\sup_{\mu\in\CQ}\int_{\Omega} f\Td\mu$ holds trivially. 
For any $c\in\R,\BIy\in\R^m$ such that $c+\langle\BIy,\BIg\rangle\ge f$, and any $\mu\in\CQ$, the following holds:
\begin{align*}
c+\pi(\BIy)=c+\sup_{\BIq\in\Delta}\langle\BIy,\BIq\rangle\ge \int_{\Omega}c+\langle\BIy,\BIg\rangle \Td\mu\ge\int_{\Omega}f \Td\mu.
\end{align*}
Taking the infimum over all $c\in\R,\BIy\in\R^m$ such that $c+\langle\BIy,\BIg\rangle\ge f$ and taking the supremum over all $\mu\in\CQ$ gives $\phi(f)\ge\sup_{\mu\in\CQ}\int_{\Omega} f\Td\mu$. 
It follows that (\ref{eqn:duality}) holds. 
The proof is now complete. 
\endproof

\proof{Proof of Proposition~\ref{prop:examplena}}
Suppose that Assumption~\ref{asp:na} does not hold, and there exist $\BIy\in\R^m$ and $\widehat{\omega}\in\Omega$ such that $\langle\BIy,\BIg\rangle-\pi(\BIy)\ge0$ and $\langle\BIy,\BIg(\widehat{\omega})\rangle-\pi(\BIy)=\alpha>0$. 
By the continuity of $\BIg$, there exists an open set $E\subset\Omega$ such that $\widehat{\omega}\in E$ and $\langle\BIy,\BIg(\omega)\rangle-\pi(\BIy)>\frac{\alpha}{2}$ for all $\omega\in E$. 
Since $\widehat{\mu}$ is equivalent to the Lebesgue measure on $\Omega$, we have $\widehat{\mu}(E)>0$, and thus
\begin{align*}
\int \langle\BIy,\BIg\rangle-\pi(\BIy)\Td\widehat{\mu}\ge\int \frac{\alpha}{2}\INDI_{E}\Td\widehat{\mu}=\frac{\alpha}{2}\widehat{\mu}(E)>0.
\end{align*}
This implies that $\int \langle\BIy,\BIg\rangle \Td\widehat{\mu}>\pi(\BIy)$. 
However, since $\widehat{\mu}\in\CQ$, we also have $\int \BIg \Td\widehat{\mu}\in\Delta$, therefore $\int \langle\BIy,\BIg\rangle \Td\widehat{\mu}\le\sup_{\BIq\in\Delta}\langle\BIy,\BIq\rangle=\pi(\BIy)$, which is a contradiction. 
The proof is now complete. 
\endproof

\section{Proofs of auxiliary results related to the numerical methods}
\label{sec:proofalgo}

\proof{Proof of Lemma~\ref{lem:cpwaproperties}}
Property (i) is immediately clear from Definition~\ref{def:cpwa}. 
To prove the property~(ii), suppose that $h$ has the form in (\ref{eqn:cpwadef}). 
Notice that each term in (\ref{eqn:cpwadef}) can be written as follows,
\begin{align*}
\xi_k\max\left\{\langle\BIa_{k,i},\BIx\rangle+b_{k,i}:1\le i\le I_k\right\}=\begin{cases}
\langle\xi_k\BIa_{k,1},\BIx\rangle+\xi_kb_{k,1}, & \text{if }\BIx\in \Omega_{k,1},\\
\qquad\vdots & \quad\vdots\\
\langle\xi_k\BIa_{k,I_k},\BIx\rangle+\xi_kb_{k,I_k}, & \text{if }\BIx\in \Omega_{k,I_k},
\end{cases}
\end{align*}
where $\Omega_{k,i}:=\{\BIx\in\R^d_+:\langle\BIa_{k,i},\BIx\rangle+b_{k,i}\ge\langle\BIa_{k,i'},\BIx\rangle+b_{k,i'}\;\forall 1\le i'\le I_k\}$ is a polyhedron (if non-empty) in $\R^d_+$ and $\bigcup_{i=1}^{I_k}\Omega_{k,i}=\R^d_+$. 
Moreover, $\inter(\Omega_{k,i})=\{\BIx\in\R^d_+:\langle\BIa_{k,i},\BIx\rangle+b_{k,i}>\langle\BIa_{k,i'},\BIx\rangle+b_{k,i'}\;\forall 1\le i'\le I_k\}$. 
Thus, for $i\ne i'$, $\inter(\Omega_{k,i})\cap\inter(\Omega_{k,i'})=\emptyset$. 
If $\BIx\in\Omega_{k,i}\cap\Omega_{k,i'}$, then by the definition of $\Omega_{k,i}$ and $\Omega_{k,i'}$ we have $\langle\xi_k\BIa_{k,i},\BIx\rangle+\xi_kb_{k,i}=\langle\xi_k\BIa_{k,i'},\BIx\rangle+\xi_kb_{k,i'}$. 
Let $\CI:=\{(i_k)_{k=1:K}:1\le i_k\le I_k\}$. 
$h$ can be decomposed into $|\CI|=\prod_{k=1}^KI_k$ cases and expressed by the following local representation (note that for some cases the corresponding set can be empty):
\begin{align}
h(\BIx)=\begin{cases}
\left\langle\sum_{k=1}^K\xi_k\BIa_{k,1},\BIx\right\rangle+\sum_{k=1}^K\xi_kb_{k,1}, & \text{if }\BIx\in\bigcap_{k=1}^K\Omega_{k,1},\\
\qquad\qquad\vdots & \qquad\vdots\\
\left\langle\sum_{k=1}^K\xi_k\BIa_{k,i_k},\BIx\right\rangle+\sum_{k=1}^K\xi_kb_{k,i_k}, & \text{if }\BIx\in\bigcap_{k=1}^K\Omega_{k,i_k},\\
\qquad\qquad\vdots & \qquad\vdots\\
\left\langle\sum_{k=1}^K\xi_k\BIa_{k,I_k},\BIx\right\rangle+\sum_{k=1}^K\xi_kb_{k,I_k}, & \text{if }\BIx\in\bigcap_{k=1}^K\Omega_{k,I_k}.
\end{cases}
\label{eqn:cpwaexplicitlocalrep}
\end{align}
Let $\Omega_{(i_k)}:=\bigcap_{k=1}^K\Omega_{k,i_k}$. 
It is straightforward to verify that for each $(i_k)\in\CI$, $\Omega_{(i_k)}$ is a polyhedron (if non-empty), and $\bigcup_{(i_k)\in\CI}\Omega_{(i_k)}=\R^d_+$. 
Let $(i_k)\in\CI,(i'_k)\in\CI$. 
If $i_j\ne i'_j$ for some $j\in\{1,\ldots,K\}$, then $\inter(\Omega_{(i_k)})\cap\inter(\Omega_{(i'_k)})\subset\inter(\Omega_{j,i_j})\cap\inter(\Omega_{j,i'_j})=\emptyset$. 
If $\BIx\in\Omega_{(i_k)}\cap\Omega_{(i'_k)}$, then
\begin{align*}
\left\langle\xi_k\BIa_{k,i_k},\BIx\right\rangle+\xi_kb_{k,i_k}=\left\langle\xi_k\BIa_{k,i'_k},\BIx\right\rangle+\xi_kb_{k,i'_k}
\end{align*}
for all $k=1,\ldots,K$, thus
\begin{align*}
\left\langle\sum_{k=1}^K\xi_k\BIa_{k,i_k},\BIx\right\rangle+\sum_{k=1}^K\xi_kb_{k,i_k}=\left\langle\sum_{k=1}^K\xi_k\BIa_{k,i'_k},\BIx\right\rangle+\sum_{k=1}^K\xi_kb_{k,i'_k}.
\end{align*}
We have thus completed the proof of property~(ii). 

To prove property (iii), we divide the collection $\CO:=\{\Omega_j:j=1,\ldots,J\}$ of all polyhedra in the representation (\ref{eqn:cpwadef}) into $\CN:=\{\Omega_j\in\CO:\inter(\Omega_j)\ne\emptyset\}$ and $\CO\setminus\CN$. 
We claim that $\bigcup_{O\in\CN}O=\R^d_+$. 
Suppose that $\R^d_+\setminus\left(\bigcup_{O\in\CN}O\right)\ne\emptyset$. $\bigcup_{O\in\CN}O$ is closed since it is a finite union of closed sets. 
Consequently, by property~(ii), $\R^d_+\setminus\left(\bigcup_{O\in\CN}O\right)$ is a non-empty open set contained in $\bigcup_{O\in\CO\setminus\CN}O$, which contradicts the fact that $\bigcup_{O\in\CO\setminus\CN}O$ has empty interior.
Hence, in the local representation (\ref{eqn:cpwalocal}), it is sufficient to only enumerate over those polyhedra with non-empty interior, and we have proved property~(iii). 
The proof is now complete. 
\endproof

\proof{Proof of Proposition~\ref{prop:cpwaboundedness}}
Let us first prove statement~(i). 
Suppose that $h$ has the representation (\ref{eqn:cpwadef}). 
Notice that for fixed $\BIx_0\in\R^d_+$ and $\BIz\in\R^d_+$, $\gamma_k(t):=\max\left\{\langle\BIa_{k,i},\BIx_0+t\BIz\rangle+b_{k,i}:1\le i\le I_k\right\}$ is a CPWA function in $t$ defined on $\R_+$. 
Thus, there exists $\overline{t}_k>0$ such that $\gamma_k$ is an affine function on $[\overline{t}_k,\infty)$.
Hence, for all $t\ge\overline{t}_k$,
\begin{align*}
\gamma_k(t)=&\max\left\{\langle\BIa_{k,i},\BIx_0+\overline{t}_k\BIz\rangle+b_{k,i}:1\le i\le I_k\right\}+(t-\overline{t}_k)\max\left\{\langle\BIa_{k,i},\BIz\rangle:1\le i\le I_k\right\}.
\end{align*}
Therefore, for $\overline{t}:=\max\{\overline{t}_1,\ldots,\overline{t}_K\}$ and all $t\ge\overline{t}$,
\begin{align}
\begin{split}
h(\BIx_0+t\BIz)=&\sum_{k=1}^{K}\xi_k\max\left\{\langle\BIa_{k,i},\BIx_0+\overline{t}\BIz\rangle+b_{k,i}:1\le i\le I_k\right\}\\
&+(t-\overline{t})\sum_{k=1}^{K}\xi_k\max\left\{\langle\BIa_{k,i},\BIz\rangle:1\le i\le I_k\right\}\\
=&\sum_{k=1}^{K}\xi_k\max\left\{\langle\BIa_{k,i},\BIx_0+\overline{t}\BIz\rangle+b_{k,i}:1\le i\le I_k\right\}+(t-\overline{t})\widetilde{h}(\BIz).
\end{split}
\label{eqn:cpwaonray}
\end{align}

Let us suppose first that $\inf_{\BIx\in\R^d_+}h(\BIx)>-\infty$, that is, there exists $\alpha\in\R$ such that $h(\BIx)\ge \alpha$ for all $\BIx\in\R^d_+$. 
Thus, $h(\BIx_0+t\BIz)\ge \alpha$ for all $t\ge 0$, and it holds by (\ref{eqn:cpwaonray}) that $\widetilde{h}(\BIz)\ge 0$. 
The same reasoning applies for every $\BIz\in\R^d_+$, and thus $\widetilde{h}(\BIz)\ge 0$ for all $\BIz\in\R^d_+$. 

Conversely, suppose that $\inf_{\BIx\in\R^d_+}h(\BIx)=-\infty$. 
By Lemma~\ref{lem:cpwaproperties}(ii), $h$ admits the following local representation,
\begin{align}
h(\BIx)=\begin{cases}
\langle\BIa_{1},\BIx\rangle+b_{1}, & \text{if }\BIx\in\Omega_1,\\
\qquad\vdots & \quad\vdots\\
\langle\BIa_{J},\BIx\rangle+b_{J}, & \text{if }\BIx\in\Omega_J,
\end{cases}
\label{eqn:radfunclocalrep}
\end{align}
hence there exists a $j\in\{1,\ldots,J\}$ such that $\inf\{\langle\BIa_{j},\BIx\rangle:\BIx\in\Omega_j\}=-\infty$. 
By well-known results about polyhedra (see \textit{e.g.} Theorem~19.1 of \citet{rockafellar1970convex}), the polyhedron $\Omega_j$ can be represented by its finite number of extreme points and extreme directions, that is,
\begin{align}
\Omega_j=\left\{\sum_{\BIv\in V}\lambda_{\BIv}\BIv+\sum_{\BIz\in D}\zeta_{\BIz}\BIz:\lambda_{\BIv}\ge0\;\forall\BIv\in V,\sum_{\BIv\in V}\lambda_{\BIv}=1,\zeta_{\BIz}\ge0\;\forall\BIz\in D\right\},
\label{eqn:polyhrep}
\end{align}
where $V\subset\R^d_+$ is the finite set of extreme points and $D\subset\R^d_+$ is the finite set of extreme directions.
Given this representation, we have
\begin{align*}
\inf\{\langle\BIa_{j}&,\BIx\rangle:\BIx\in\Omega_j\}\\
&=\inf\left\{\sum_{\BIv\in V}\lambda_{\BIv}\langle\BIa_j,\BIv\rangle+\sum_{\BIz\in D}\zeta_{\BIz}\langle\BIa_j,\BIz\rangle:\lambda_{\BIv}\ge0\;\forall\BIv\in V,\sum_{\BIv\in V}\lambda_{\BIv}=1,\zeta_{\BIz}\ge0\;\forall\BIz\in D\right\},
\end{align*}
and thus,
\begin{align}
\inf\{\langle\BIa_{j},\BIx\rangle:\BIx\in\Omega_j\}=-\infty\quad\Longrightarrow\quad\exists\BIz_0\in D \text{ such that }\langle\BIa_j,\BIz_0\rangle<0. 
\label{eqn:radfuncintstep}
\end{align}
By the representation (\ref{eqn:polyhrep}), for fixed $\BIx_0\in\Omega_j$ and $\BIz_0$ given by (\ref{eqn:radfuncintstep}), $\BIx_0+t\BIz_0\in\Omega_j$ for all $t\ge 0$ as $\BIz_0$ is an extreme direction. Thus, by (\ref{eqn:radfunclocalrep}), for all $t\ge0$,
\begin{align*}
h(\BIx_0+t\BIz_0)=\langle\BIa_{j},\BIx_0\rangle+b_{j}+t\langle\BIa_{j},\BIz_0\rangle.
\end{align*}
Since $\langle\BIa_j,\BIz_0\rangle<0$, we have $\lim_{t\to\infty}h(\BIx_0+t\BIz_0)=-\infty$. 
By (\ref{eqn:cpwaonray}), we can conclude that $\widetilde{h}(\BIz_0)<0$. 
This settles statement~(i).

To prove statement~(ii), suppose that $\inf_{\BIx\in\R^d_+}h(\BIx)>-\infty$. 
Again by Lemma~\ref{lem:cpwaproperties}(ii), $h$ admits the local representation in (\ref{eqn:radfunclocalrep}), hence $\inf_{\BIx\in\Omega_j}h(\BIx)>-\infty$ for $j=1,\ldots,J$. 
For each $j$, $\Omega_j$ can be decomposed as in (\ref{eqn:polyhrep}) into the set of extreme points $V$ and the set of extreme directions $D$. 
By (\ref{eqn:polyhrep}) and the fact that $\inf_{\BIx\in\Omega_j}h(\BIx)>-\infty$, $\langle\BIa_j,\BIz\rangle\ge0$ for all $\BIz\in D$. 
Hence, there exists $\BIx_j^\star\in V$ such that $h(\BIx_j^\star)=\inf_{\BIx\in\Omega_j}h(\BIx)$. 
Therefore, there exists $\BIx^\star\in\R^d_+$ such that $h(\BIx^\star)=\inf_{\BIx\in\R^d_+}h(\BIx)$. 
This settles statement~(ii).
The proof is now complete. 
\endproof


\proof{Proof of Proposition~\ref{prop:slackfunc}}
Statement~(i) follows directly from Assumption~\ref{asp:setting1}, Definition~\ref{def:slackfunc}, and Lemma~\ref{lem:cpwaproperties}(i). 

To prove statement~(ii), let us explicitly express $g_j$ for $j=1,\ldots,m$ and $f$ as CPWA functions as in (\ref{eqn:cpwadef}): for $j=1,\ldots,m$,
\begin{align}
g_j(\BIx)=&\sum_{k=1}^{K^{g_j}}\xi^{g_j}_k\max\left\{\langle\BIa^{g_j}_{k,i},\BIx\rangle+b^{g_j}_{k,i}:1\le i\le I^{g_j}_k\right\},
\label{eqn:cpwagdef}
\end{align}
where $K^{g_j}\in\N$, $I^{g_j}_k\in\N$ for $k=1,\ldots,K^{g_j}$, $\BIa^{g_j}_{k,i}\in\R^d$, $b^{g_j}_{k,i}\in\R$, $\xi^{g_j}_k\in\{-1,1\}$, for $i=1,\ldots,I^{g_j}_k,k=1,\ldots,K^{g_j}$,
\begin{align}
f(\BIx)=&\sum_{k=1}^{K^{f}}\xi^{f}_k\max\left\{\langle\BIa^{f}_{k,i},\BIx\rangle+b^{f}_{k,i}:1\le i\le I^{f}_k\right\},
\label{eqn:cpwafdef}
\end{align}
where $K^{f}\in\N$, $I^{f}_k\in\N$ for $k=1,\ldots,K^{f}$, $\BIa^{f}_{k,i}\in\R^d$, $b^{f}_{k,i}\in\R$, $\xi^{f}_k\in\{-1,1\}$, for $i=1,\ldots,I^{f}_k,k=1,\ldots,K^{f}$. 
By Definition~\ref{def:slackfunc} and (\ref{eqn:cpwagdef}), (\ref{eqn:cpwafdef}),
\begin{align*}
s_{\BIy}(\BIx)=&\sum_{j=1}^my_jg_j-f\\
=&\sum_{j=1}^m\sum_{k=1}^{K^{g_j}}\xi^{g_j}_ky_j\max\left\{\langle\BIa^{g_j}_{k,i},\BIx\rangle+b^{g_j}_{k,i}:1\le i\le I^{g_j}_k\right\}\\
&-\sum_{k=1}^{K^{f}}\xi^{f}_k\max\left\{\langle\BIa^{f}_{k,i},\BIx\rangle+b^{f}_{k,i}:1\le i\le I^{f}_k\right\},
\end{align*}
which has the desired form in (\ref{eqn:lsipconspoint}). 
Notice that one may extract positive coefficients out of $\max\{\cdot\}$ and combine terms to simplify this representation. 
The representation (\ref{eqn:lsipconsdir}) follows directly from (\ref{eqn:lsipconspoint}) and Definition~\ref{def:cpwa}. 
We remark that the number of terms in (\ref{eqn:lsipconsdir}) could be different from the number of terms in (\ref{eqn:lsipconspoint}) due to possible simplification. 
We have now proved statement~(ii). 

Statement~(iii) is a direct consequence of Proposition~\ref{prop:cpwaboundedness}. 

Finally, to prove statement~(iv), let us fix a $\BIy$ such that $\inf_{\BIx\in\R^d_+}s_{\BIy}(\BIx)>-\infty$ and express $s_{\BIy}(\BIx)$ by (\ref{eqn:lsipconspoint}). 
Let $\CI:=\{(i_k)_{k=1:K}:1\le i_k\le I_k\}$. 
Same as in the proof of Lemma~\ref{lem:cpwaproperties}(ii), $s_{\BIy}(\BIx)$ has a local representation as in (\ref{eqn:cpwaexplicitlocalrep}), with $\xi_k$ replaced by $\langle\BIy,\BIw_k\rangle+z_k$, over the finite partition $\widetilde{\FC}:=\left\{\bigcap_{k=1}^K\Omega_{k,i_k}:(i_k)\in\CI\right\}$ of $\R^d_+$, where $\Omega_{k,i}:=\{\BIx\in\R^d_+:\langle\BIa_{k,i},\BIx\rangle+b_{k,i}\ge\langle\BIa_{k,i'},\BIx\rangle+b_{k,i'}\;\forall 1\le i'\le I_k\}$. 
Notice that the partition $\widetilde{\FC}$ does not depend on the value of $\BIy$, and each $C\in\widetilde{\FC}$ is a polyhedron. Let $X^\star:=\bigcup_{C\in\widetilde{\FC}}V(C)$, where $V(C)$ denotes the finite set of extreme points of the polyhedron $C$ given in (\ref{eqn:polyhrep}). 
Thus, $|X^\star|<\infty$, and by (\ref{eqn:polyhrep}), for each $C\in\widetilde{\FC}$, there exists $\BIx^\star\in V(C)$ (depending on $\BIy$) such that $s_{\BIy}(\BIx^\star)=\inf_{\BIx\in C}s_{\BIy}(\BIx)$. 
Let us choose $\overline{\BIx}\in\R^d$ with $\overline{\BIx}>\veczero$ such that $X^\star\subset\{\BIx\in\R^d:\veczero\le\BIx\le\overline{\BIx}\}$. Subsequently, for every $\BIy\in\R^m$ such that $\inf_{\BIx\in\R^d_+}s_{\BIy}(\BIx)>-\infty$, we have $\inf_{\veczero\le\BIx\le\overline{\BIx}}s_{\BIy}(\BIx)=\inf_{\BIx\in X^\star}s_{\BIy}(\BIx)=\inf_{\BIx\in\R^d_+}s_{\BIy}(\BIx)$. Let $\FC$ be the finite partition of $\{\BIx\in\R^d:\veczero\le\BIx\le\overline{\BIx}\}$ defined by
\begin{align*}
\FC:=\Big\{C\cap\{\BIx\in\R^d:\veczero\le\BIx\le\overline{\BIx}\}:C\in\widetilde{\FC}\Big\}.
\end{align*}
Notice that for every $\BIy\in\R^m$, the function $s_{\BIy}(\cdot)$ is affine when restricted to each $C\in\FC$, and thus $\FC$ satisfies all required properties. 
The proof is now complete. 
\endproof

\proof{Proof of Proposition~\ref{prop:radcons}}
Let us first prove statement~(i). Since 
$$\widetilde{s}_{\BIy}(\BIz)=\sum_{k=1}^{\widetilde{K}}(\langle\BIy,\widetilde{\BIw}_k\rangle+\widetilde{z}_k)\max\left\{\langle\widetilde{\BIa}_{k,i},\BIz\rangle:1\le i\le \widetilde{I}_k\right\},$$ 
we have for fixed $(i_k)\in\widetilde{\CI}$ and all 
\begin{align*}
\BIz\in\left\{\BIz\in\R^d_+:\langle\widetilde{\BIa}_{k,i_k},\BIz\rangle\ge\langle\widetilde{\BIa}_{k,i},\BIz\rangle \;\forall i\in\{1,\ldots,\widetilde{I}_k\}\setminus\{i_k\},k=1,\ldots,\widetilde{K}\right\}
\end{align*}
that
\begin{align*}
\widetilde{s}_{\BIy}(\BIz)=\sum_{k=1}^{\widetilde{K}}(\langle\BIy,\widetilde{\BIw}_k\rangle+\widetilde{z}_k)\langle\widetilde{\BIa}_{k,i_k},\BIz\rangle.
\end{align*}
But
\begin{align}
\begin{split}
&\!\!\!\!\!\!\!\!\!\!\!\! \left\{\BIz\in\R^d_+:\langle\widetilde{\BIa}_{k,i_k},\BIz\rangle\ge\langle\widetilde{\BIa}_{k,i},\BIz\rangle \;\forall i\in\{1,\ldots,\widetilde{I}_k\}\setminus\{i_k\},k=1,\ldots,\widetilde{K}\right\}\\
&= \left\{\BIz\in\R^d_+:\langle\widetilde{\BIa}_{k,i_k}-\widetilde{\BIa}_{k,i},\BIz\rangle\ge0 \;\forall i\in\{1,\ldots,\widetilde{I}_k\}\setminus\{i_k\},k=1,\ldots,\widetilde{K}\right\}\\
&= \dual\left(\left\{\widetilde{\BIa}_{k,i_k}-\widetilde{\BIa}_{k,i}:i\in\{1,\ldots,\widetilde{I}_k\}\setminus\{i_k\},k=1,\ldots,\widetilde{K}\right\}\setminus\{\veczero\}\right)\cap\R^d_+\\
&= \dual(A_{(i_k)})\cap\R^d_+.
\end{split}
\label{eqn:dualconerep}
\end{align}
By (\ref{eqn:dualconerep}), it is straightforward to verify that $\bigcup_{(i_k)\in\widetilde{\CI}}\dual(A_{(i_k)})\cap\R^d_+=\R^d_+$, and that if $(i_k)\in\widetilde{\CI}$, $(i'_k)\in\widetilde{\CI}$, $(i_k)\ne(i'_k)$, then $\inter(\dual(A_{(i_k)})\cap\R^d_+)\cap\inter(\dual(A_{(i'_k)})\cap\R^d_+)=\emptyset$. 
We have thus expressed $\widetilde{s}_{\BIy}(\BIz)$ by a local representation as in Lemma~\ref{lem:cpwaproperties}(ii) over the partition $\big\{\dual(A_{(i_k)})\cap\R^d_+:(i_k)\in\widetilde{\CI}\big\}$ of $\R^d_+$. 
The last part of statement~(i) is due to Lemma~\ref{lem:cpwaproperties}(iii). 
We have thus proved statement~(i).

To prove statement~(ii), let us fix $(i_k)\in\widetilde{\CI}$. 
By statement~(i), $\widetilde{s}_{\BIy}(\BIz)\ge0$ for all $\BIz\in\dual(A_{(i_k)})\cap\R^d_+$ is equivalent to $\sum_{k=1}^{\widetilde{K}}(\langle\BIy,\widetilde{\BIw}_k\rangle+\widetilde{z}_k)\langle\widetilde{\BIa}_{k,i_k},\BIz\rangle\ge0$ for all $\BIz\in\dual(A_{(i_k)})\cap\R^d_+$. 
This is further equivalent to
\begin{align*}
\sum_{k=1}^{\widetilde{K}}(\langle\BIy,\widetilde{\BIw}_k\rangle+\widetilde{z}_k)\widetilde{\BIa}_{k,i_k}
	\in\dual(\dual(A_{(i_k)})\cap\R^d_+)
	&=\dual(\dual(A_{(i_k)}+\R^d_+))\\
	&=\cone(A_{(i_k)})+\R^d_+.
\end{align*}
Therefore,
\begin{align*}
\widetilde{s}_{\BIy}(\BIz)\ge0\ \forall \BIz\in\dual(A_{(i_k)})\cap\R^d_+\quad\Longleftrightarrow\quad &\exists\BIu\in\cone(A_{(i_k)}) \text{ s.t. }\sum_{k=1}^{\widetilde{K}}(\langle\BIy,\widetilde{\BIw}_k\rangle+\widetilde{z}_k)\widetilde{\BIa}_{k,i_k}\ge\BIu,
\end{align*}
which is equivalent to the existence of $\eta^{(i_k)}_\BIv\ge0$ for each $\BIv\in A_{(i_k)}$ such that (\ref{eqn:lpconsdir}) holds. 
We have now proved statement~(ii). 

Let us now prove statement~(iii).
Suppose first that $\inter(\dual(A_{(i_k)})\cap\R^d_+)=\emptyset$. 
Since $\dual(A_{(i_k)})\cap\R^d_+$ is a polyhedral cone, it is contained in a $(d-1)$-dimensional subspace, say $\{\BIw\in\R^n:\langle\BIu_0,\BIw\rangle=0\}$ for some $\BIu_0\ne\veczero$, which implies that both $\BIu_0$ and $-\BIu_0$ are elements of $\dual(\dual(A_{(i_k)})\cap\R^d_+)=\cone(A_{(i_k)})+\R^d_+$. 
This implies that there exist $\BIu_1,\BIu_2\in\cone(A_{(i_k)})$ such that $\BIu_1\le\BIu_0,\BIu_2\le-\BIu_0$. 
Thus, $\BIu:=\BIu_1+\BIu_2\le\veczero$, and $\BIu\in\cone(A_{(i_k)})$. 
This further implies that for each $\BIv\in A_{(i_k)}$ there exists $\eta^{(i_k)}_{\BIv}\ge0$ such that $\BIu=\sum_{\BIv\in A_{(i_k)}}\eta^{(i_k)}_{\BIv} \BIv$, while $(\eta^{(i_k)}_{\BIv})_{\BIv\in A_{(i_k)}}$ are not all identically zero (otherwise we will have $\BIu_1=\BIu_2=\veczero$, and thus $\BIu_0=\veczero$). 
We are allowed to scale $(\eta^{(i_k)}_{\BIv})_{\BIv\in A_{(i_k)}}$ by the same positive factor to make $\sum_{\BIv\in A_{(i_k)}}\eta^{(i_k)}_{\BIv}=1$, thus condition (\ref{eqn:lpconsdircheck}) holds.

Conversely, suppose that the condition (\ref{eqn:lpconsdircheck}) holds, and let $(\eta^{(i_k)}_{\BIv})_{\BIv\in A_{(i_k)}}$ be given as in (\ref{eqn:lpconsdircheck}). 
Define $\BIu:=\sum_{\BIv\in A_{(i_k)}}\eta^{(i_k)}_{\BIv} \BIv$. By (\ref{eqn:lpconsdircheck}), $\BIu\le\veczero$. If $\BIu\ne\veczero$, then $\BIu\in\cone(A_{(i_k)})+\R^d_+$, $-\BIu\in\R^d_+\subset\cone(A_{(i_k)})+\R^d_+$, and $\dual(A_{(i_k)})\cap\R^d_+=\dual(\cone(A_{(i_k)}))\cap\R^d_+=\dual(\cone(A_{(i_k)})+\R^d_+)$ is a subset of the subspace $\{\BIw\in\R^d:\langle\BIu,\BIw\rangle=0\}$, thus $\inter(\dual(A_{(i_k)})\cap\R^d_+)=\emptyset$. If $\BIu=\veczero$, then let $\BIv'\in A_{(i_k)}$ such that $\eta^{(i_k)}_{\BIv'}>0$. 
Thus, $\eta^{(i_k)}_{\BIv'}\BIv'\ne\veczero$ since $\veczero\notin A_{(i_k)}$, $\eta^{(i_k)}_{\BIv'}\BIv'\in\cone(A_{(i_k)})+\R^d_+$, and
\begin{align*}
-\eta^{(i_k)}_{\BIv'}\BIv'=\BIu-\eta^{(i_k)}_{\BIv'}\BIv'=\sum_{\BIv\in A_{(i_k)},\BIv\ne\BIv'}\eta^{(i_k)}_{\BIv}\BIv\in\cone(A_{(i_k)})+\R^d_+.
\end{align*}
Since both $\eta^{(i_k)}_{\BIv'}\BIv'$ and $-\eta^{(i_k)}_{\BIv'}\BIv'$ are elements of $\cone(A_{(i_k)})+\R^d_+$, one can deduce that $\dual(A_{(i_k)})\cap\R^d_+=\dual(\cone(A_{(i_k)}))\cap\R^d_+=\dual(\cone(A_{(i_k)})+\R^d_+)$ is a subset of the subspace $\{\BIw\in\R^d:\langle\BIv',\BIw\rangle=0\}$, thus $\inter(\dual(A_{(i_k)})\cap\R^d_+)=\emptyset$. 
The proof is now complete.
\endproof

\proof{Proof of Proposition~\ref{prop:radconsalgo}}
Algorithm~\ref{alg:radcons} guarantees that the linear constraint (\ref{eqn:lpconsdir}) with auxiliary variables is included in the system of linear inequalities $\widetilde{\sigma}$ for every $(i_k)\in\widetilde{\CI}$ such that the condition (\ref{eqn:lpconsdircheck}) fails to hold. Proposition~\ref{prop:radcons}(iii) states that condition (\ref{eqn:lpconsdircheck}) fails to hold if and only if $\inter(\dual(A_{(i_k)})\cap\R^d_+)\ne\emptyset$. Hence, (\ref{eqn:lpconsdir}) is included in $\widetilde{\sigma}$ for every $(i_k)\in\widetilde{\CI}$ such that $\inter(\dual(A_{(i_k)})\cap\R^d_+)\ne\emptyset$. 
Moreover, Proposition~\ref{prop:radcons}(i) implies that 
\begin{align*}
\bigcup_{(i_k)\in\widetilde{\CI},\;\inter(\dual(A_{(i_k)})\cap\R^d_+)\ne\emptyset}\left(\dual(A_{(i_k)})\cap\R^d_+\right)=\R^d_+.
\end{align*}
Consequently, Proposition~\ref{prop:radcons}(ii) and Proposition~\ref{prop:slackfunc}(iii) imply that 
\begin{align*}
\BIy\text{ satisfies all constraints in }\widetilde{\sigma}\quad\Longleftrightarrow\quad \widetilde{s}_{\BIy}(\BIz)\ge0 \quad \forall\BIz\in\R^d_+ \quad \Longleftrightarrow\quad \inf_{\BIx\in\R^d_+}s_{\BIy}(\BIx)>-\infty.
\end{align*}
This completes the proof. 
\endproof

\proof{Proof of Lemma~\ref{lem:milp}}
Let $p$ denote the objective function of problem (\ref{eqn:milp}). 
To prove Lemma~\ref{lem:milp}, we show that each minimizer of problem~(\ref{eqn:milp}) gives a feasible point of problem~(\ref{eqn:milpglobalmin}) with identical objective value and vice versa. 
Suppose that $(\BIx^\star,(\lambda^\star_k),(\zeta^\star_k),(\delta^\star_{k,i}),(\iota^\star_{k,i}))$ is a minimizer of problem~(\ref{eqn:milp}) (which exists since (\ref{eqn:milp}) is feasible and bounded, and can be considered as the minimum of a finite collection of LP problems, each corresponding to a feasible combination of binary variables $(\iota_{k,i})$). 
By the constraints of (\ref{eqn:milp}), for each $k$ such that $\xi_k=1$, $\lambda^\star_k\ge\max_{1\le i\le I_k}\left\{\langle\BIa_{k,i},\BIx^\star\rangle+b_{k,i}\right\}$, 
and since $\lambda^\star_k$ is the minimum of such $\lambda_k$, $\lambda^\star_k=\max_{1\le i\le I_k}\left\{\langle\BIa_{k,i},\BIx^\star\rangle+b_{k,i}\right\}$. 
For each $k$ such that $\xi_k=-1$, $\zeta^\star_k\ge\max_{1\le i\le I_k}\left\{\langle\BIa_{k,i},\BIx^\star\rangle+b_{k,i}\right\}$. By $\iota_{k,i}\in\{0,1\}$ for all $1\le i\le I_k$ and $\sum_{i=1}^{I_k}\iota_{k,i}=1$, there exists a unique $i'\in\{1,\ldots,I_k\}$ such that $\iota_{k,i'}=1$. 
Subsequently, $\delta_{k,i'}=0$ and $\zeta^\star_k=\langle\BIa_{k,i},\BIx^\star\rangle+b_{k,i'}$ and thus $\zeta^\star_k=\max_{1\le i\le I_k}\left\{\langle\BIa_{k,i},\BIx^\star\rangle+b_{k,i}\right\}$. 
We have
\begin{align*}
p(\BIx^\star,(\lambda^\star_k),(\zeta^\star_k),(\delta^\star_{k,i}),(\iota^\star_{k,i}))=&\sum_{k=1,\ldots,K,\xi_k=1}\lambda^\star_k+\sum_{k=1,\ldots,K,\xi_k=-1}-\zeta^\star_k\\
=&\sum_{k=1,\ldots,K,\xi_k=1}\max_{1\le i\le I_k}\left\{\langle\BIa_{k,i},\BIx^\star\rangle+b_{k,i}\right\}\\
&+\sum_{k=1,\ldots,K,\xi_k=-1}-\max_{1\le i\le I_k}\left\{\langle\BIa_{k,i},\BIx^\star\rangle+b_{k,i}\right\}\\
=&\sum_{k=1,\ldots,K}\xi_k\max_{1\le i\le I_k}\left\{\langle\BIa_{k,i},\BIx^\star\rangle+b_{k,i}\right\}\\
=&h(\BIx^\star).
\end{align*}
Hence, $\BIx^\star$ is a feasible point of problem~(\ref{eqn:milpglobalmin}) with identical objective value. 
Conversely, let $\BIx^\star$ be a minimizer of problem~(\ref{eqn:milpglobalmin}). 
For each $k$ such that $\xi_k=1$, define $\lambda^\star_k:=\max_{1\le i\le I_k}\left\{\langle\BIa_{k,i},\BIx^\star\rangle+b_{k,i}\right\}$. 
For each $k$ such that $\xi_k=-1$, define $\zeta^\star_k:=\max_{1\le i\le I_k}\left\{\langle\BIa_{k,i},\BIx^\star\rangle+b_{k,i}\right\}$ and for $i=1,\ldots,I_k$, define $\delta^\star_{k,i}:=\zeta^\star_k-\langle\BIa_{k,i},\BIx^\star\rangle-b_{k,i}$. 
It is clear that $\delta^\star_{k,i}\ge0$. $\delta^\star_{k,i}\le M_{k,i}$ is guaranteed by the definition of $M_{k,i}$ in (\ref{eqn:milpM}). Let $i'\in\{1,\ldots,I_k\}$ such that $\delta^\star_{k,i'}=0$, and define $\iota^\star_{k,i'}:=1$ and $\iota^\star_{k,i}:=0$ for all $i\ne i'$. 
It is straightforward to verify that $(\BIx^\star,(\lambda^\star_k),(\zeta^\star_k),(\delta^\star_{k,i}),(\iota^\star_{k,i}))$ is a feasible point of problem~(\ref{eqn:milp}). 
Moreover, by our definitions
\begin{align*}
h(\BIx^\star)=&\sum_{k=1,\ldots,K,\xi_k=1}\lambda^\star_k+\sum_{k=1,\ldots,K,\xi_k=-1}-\zeta^\star_k\\
=&p(\BIx^\star,(\lambda^\star_k),(\zeta^\star_k),(\delta^\star_{k,i}),(\iota^\star_{k,i})).
\end{align*}
Hence, $(\BIx^\star,(\lambda^\star_k),(\zeta^\star_k),(\delta^\star_{k,i}),(\iota^\star_{k,i}))$ is a feasible point of problem~(\ref{eqn:milp}) with identical objective value. 
The proof is now complete. 
\endproof

\section{Proofs related to the arbitrage removal process}
\label{ecsec:proofarbitrageremoval}

\proof{Proof of Proposition~\ref{prop:no-arbitrage}}
Suppose, for the sake of contradiction, that there exist $y^{\CALL}_1$, $\ldots$, $y^{\CALL}_m$, $y^{\PUT}_1$, $\ldots$, $y^{\PUT}_m$ such that (\ref{eqn:no-arbitrage-nonneg}) holds but (\ref{eqn:no-arbitrage-zero}) does not hold. Let $h:\R_+\to\R$ denote the left-hand side of the inequality in (\ref{eqn:no-arbitrage-nonneg}) as a function of $x$, \textit{i.e.}
\begin{align*}
h(x):=\left[\sum_{j=1}^m y^{\CALL}_j(x-\kappa_j)^++y^{\PUT}_j(\kappa_j-x)^+\right]-\pi\big(y^{\CALL}_1,\ldots,y^{\CALL}_m,y^{\PUT}_1,\ldots,y^{\PUT}_m\big)\quad\forall x\in\R_+.
\end{align*}
Observe that $h$ is a continuous function which is piece-wise affine on intervals $[0,\kappa_1]$, $[\kappa_1,\kappa_2]$, $\ldots$, $[\kappa_{m-1},\kappa_m]$, and $[\kappa_m,\infty)$. 
By our assumption, $h$ is non-negative on $\R_+$ but is not everywhere zero. Hence, there exists $\widehat{x}\in\{0,\kappa_1,\ldots,\kappa_m,\overline{x}\}$ such that $h(\widehat{x})=\alpha>0$. 
By the continuity of $h$, this implies that there exists an open interval $E\subset\R$ with $\widehat{x}\in E$, such that $h(x)>\frac{\alpha}{2}$ for all $x\in E\cap\R_+$. Moreover, by the assumption that $\{0,\kappa_1,\ldots,\kappa_m,\overline{x}\}\subseteq\mathrm{supp}(\mu)$, we have $\widehat{x}\in\mathrm{supp}(\mu)$ and thus $\mu\big(E\cap[0,\overline{x}]\big)>0$. Hence,
\begin{align}
\int_{[0,\overline{x}]}h\,\Td\mu\ge\int_{E\cap[0,\overline{x}]}h\,\Td\mu\ge\frac{\alpha}{2}\mu\big(E\cap[0,\overline{x}]\big)>0.
\label{eqn:no-arbitrage-contradiction}
\end{align}
On the other hand, by (\ref{eqn:no-arbitrage-measure}), it holds for $j=1,\ldots,m$ that
\begin{align}
y^{\CALL}_j\int_{[0,\overline{x}]}(x-\kappa_j)^+\,\mu(\Td x)\le\big(y^{\CALL}_j \vee 0\big)\overline{\pi}^{\CALL}_j-\big(-y^{\CALL}_j \vee 0\big)\underline{\pi}^{\CALL}_j, \label{eqn:no-arbitrage-integral1}\\
y^{\PUT}_j\int_{[0,\overline{x}]}(\kappa_j-x)^+\,\mu(\Td x)\le\big(y^{\PUT}_j \vee 0\big)\overline{\pi}^{\PUT}_j-\big(-y^{\PUT}_j \vee 0\big)\underline{\pi}^{\PUT}_j.\label{eqn:no-arbitrage-integral2}
\end{align}
Adding up both sides of (\ref{eqn:no-arbitrage-integral1}) and (\ref{eqn:no-arbitrage-integral2}) and then summing over $j=1,\ldots,m$, we get $\int_{[0,\overline{x}]}h\,\Td\mu\le 0$, which contradicts (\ref{eqn:no-arbitrage-contradiction}). The proof is now complete. 
\endproof

\proof{Proof of Proposition~\ref{prop:arbitrage-discrete}}
Let $\kappa_0:=0$ and $\kappa_{m+1}:=\overline{x}$. 
Let $\widetilde{g}_{0},\ldots,\widetilde{g}_{m+1}:[0,\overline{x}]\to\R$ be defined as follows:
\begin{align*}
\begin{split}
\widetilde{g}_0(x):=&\frac{(\kappa_1-x)^+}{\kappa_1},\\
\widetilde{g}_j(x):=&\frac{(x-\kappa_{j-1})^+}{\kappa_{j}-\kappa_{j-1}}\wedge \frac{(\kappa_{j+1}-x)^+}{\kappa_{j+1}-\kappa_j}\quad \text{for }j=1,\ldots,m,\\
\widetilde{g}_{m+1}(x):=&\frac{(x-\kappa_m)^+}{\overline{x}-\kappa_m}.
\end{split}
\end{align*}
One can check that for $i=0,\ldots,m+1$ and $j=0,\ldots,m+1$, it holds that 
\begin{align}
\widetilde{g}_i(\kappa_j)=\begin{cases}
1 & \text{if }i=j,\\
0 & \text{if }i\ne j.
\end{cases}
\label{eqn:arbitrage-discrete-basis-prop}
\end{align}
Consequently, for each $j=1,\ldots,J$, it holds that 
\begin{align}
g_j(x)=\sum_{i=0}^{m+1}g_j(\kappa_i)\widetilde{g}_i(x) \quad \forall x\in [0,\overline{x}].
\label{eqn:arbitrage-discrete-basis-rep}
\end{align} 

Now, let us fix some $\mu\in\CP([0,\overline{x}])$ with $\{0,\kappa_1,\ldots,\kappa_m,\overline{x}\}\subseteq\mathrm{supp}(\mu)$ and $\int_{[0,\overline{x}]}g_j\,\Td\mu=\pi_j$ for $j=1,\ldots,J$. For $i=0,\ldots,m+1$, let $\widetilde{\pi}_i:=\int_{[0,\overline{x}]}\widetilde{g}_i\,\Td\mu$. 
It follows from the assumption $\{0$, $\kappa_1$, $\ldots$, $\kappa_m$, $\overline{x}\}\subseteq\mathrm{supp}(\mu)$ that $\widetilde{\pi}_i>0$ for $i=0,\ldots,m+1$. 
Moreover, since $\sum_{i=0}^{m+1}\widetilde{g}_i(x)=1$ for all $x\in[0,\overline{x}]$, it holds that $\sum_{i=0}^{m+1}\widetilde{\pi}_i=1$. 
Subsequently, let us define $\widehat{\pi}\in\CP([0,\overline{x}])$ as follows:
\begin{align}
\widehat{\pi}(\Td x):=\sum_{i=0}^{m+1}\widetilde{\pi}_i\delta_{\kappa_i}(\Td x),
\label{eqn:arbitrage-discrete-meas-def}
\end{align}
where $\delta_{\kappa_i}(\Td x)$ denotes the Dirac measure at $\kappa_i$. Thus, $\mathrm{supp}(\widehat{\mu})=\{0,\kappa_1,\ldots,\kappa_m,\overline{x}\}$. 
Finally, combining (\ref{eqn:arbitrage-discrete-basis-rep}), (\ref{eqn:arbitrage-discrete-meas-def}), (\ref{eqn:arbitrage-discrete-basis-prop}), and the assumption on $\mu$, we have for $j=1,\ldots,J$ that
\begin{align*}
\int_{[0,\overline{x}]}g_j\,\Td\widehat{\mu}=&\sum_{i=0}^{m+1}g_j(\kappa_i)\int_{[0,\overline{x}]}\widetilde{g}_i\,\Td\widehat{\mu}\\
=&\sum_{i=0}^{m+1}g_j(\kappa_i)\left(\sum_{l=0}^{m+1}\widetilde{g}_i(\kappa_l)\widetilde{\pi}_l\right)\\
=&\sum_{i=0}^{m+1}g_j(\kappa_i)\widetilde{\pi}_i\\
=&\sum_{i=0}^{m+1}g_j(\kappa_i)\int_{[0,\overline{x}]}\widetilde{g}_i\,\Td\mu\\
=&\int_{[0,\overline{x}]}\textstyle\sum_{i=0}^{m+1}g_j(\kappa_i)\widetilde{g}_i\,\Td\mu\\
=&\int_{[0,\overline{x}]}g_j\,\Td\mu\\
=&\pi_j.
\end{align*}
The proof is now complete. 
\endproof



\bibliographystyle{informs2014} 

\input{ModelFree_arXiv_EC.bbl}
\end{document}

%% file: mathsymbols.tex
\usepackage{amsfonts,mathrsfs,bbm}

\newcommand{\Ta}{\mathrm{a}}

\newcommand{\Tc}{\mathrm{c}}
\newcommand{\Td}{\mathrm{d}}

\newcommand{\Tl}{\mathrm{l}}

\newcommand{\Tp}{\mathrm{p}}

\newcommand{\Tt}{\mathrm{t}}
\newcommand{\Tu}{\mathrm{u}}

\newcommand{\TA}{\mathrm{A}}
\newcommand{\TB}{\mathrm{B}}

\newcommand{\TD}{\mathrm{D}}

\newcommand{\TI}{\mathrm{I}}

\newcommand{\TL}{\mathrm{L}}

\newcommand{\TU}{\mathrm{U}}

\newcommand{\BC}{\mathbf{C}}
\newcommand{\BD}{\mathbf{D}}

\newcommand{\BIa}{{\boldsymbol{a}}}

\newcommand{\BIe}{{\boldsymbol{e}}}

\newcommand{\BIg}{{\boldsymbol{g}}}

\newcommand{\BIq}{{\boldsymbol{q}}}

\newcommand{\BIu}{{\boldsymbol{u}}}
\newcommand{\BIv}{{\boldsymbol{v}}}
\newcommand{\BIw}{{\boldsymbol{w}}}
\newcommand{\BIx}{{\boldsymbol{x}}}
\newcommand{\BIy}{{\boldsymbol{y}}}
\newcommand{\BIz}{{\boldsymbol{z}}}

\newcommand{\CB}{\mathcal{B}}
\newcommand{\CC}{\mathcal{C}}

\newcommand{\CI}{\mathcal{I}}

\newcommand{\CL}{\mathcal{L}}
\newcommand{\CM}{\mathcal{M}}
\newcommand{\CN}{\mathcal{N}}
\newcommand{\CO}{\mathcal{O}}
\newcommand{\CP}{\mathcal{P}}
\newcommand{\CQ}{\mathcal{Q}}

\newcommand{\FC}{\mathfrak{C}}

\newcommand{\FF}{\mathfrak{F}}

\newcommand{\BLambda}{{\boldsymbol{\Lambda}}}

\newcommand{\BPsi}{{\boldsymbol{\Psi}}}

\newcommand{\Bgamma}{{\boldsymbol{\gamma}}}

\newcommand{\Blambda}{{\boldsymbol{\lambda}}}

\newcommand{\Bpi}{{\boldsymbol{\pi}}}

\newcommand{\R}{\mathbb{R}} 
\newcommand{\N}{\mathbb{N}} 
\newcommand{\vecone}{\mathbf{1}} 
\newcommand{\veczero}{\mathbf{0}} 
\newcommand{\INDI}{\mathbbm{1}} 

\DeclareMathOperator*{\inter}{\mathrm{int}} 
\DeclareMathOperator*{\relint}{\mathrm{relint}} 
\DeclareMathOperator*{\cone}{\mathrm{cone}} 
\DeclareMathOperator*{\dual}{\mathrm{dual}} 
\DeclareMathOperator*{\dom}{\mathrm{dom}} 


%% file: ModelFree_arXiv.bbl
\begin{thebibliography}{75}
\providecommand{\natexlab}[1]{#1}
\providecommand{\url}[1]{\texttt{#1}}
\providecommand{\urlprefix}{URL }

\bibitem[{Acciaio et~al.(2016)Acciaio, Beiglb{\"o}ck, Penkner,
  \protect\BIBand{} Schachermayer}]{acciaio2016model}
Acciaio B, Beiglb{\"o}ck M, Penkner F, Schachermayer W (2016) A model-free
  version of the fundamental theorem of asset pricing and the super-replication
  theorem. \emph{Mathematical Finance} 26(2):233--251.

\bibitem[{Auslender et~al.(2015)Auslender, Ferrer, Goberna, \protect\BIBand{}
  L\'{o}pez}]{auslender2015comparative}
Auslender A, Ferrer A, Goberna MA, L\'{o}pez MA (2015) Comparative study of
  {RPSALG} algorithm for convex semi-infinite programming. \emph{Comput. Optim.
  Appl.} 60(1):59--87.

\bibitem[{Auslender et~al.(2009)Auslender, Goberna, \protect\BIBand{}
  L{\'o}pez}]{auslender2009penalty}
Auslender A, Goberna MA, L{\'o}pez MA (2009) Penalty and smoothing methods for
  convex semi-infinite programming. \emph{Mathematics of Operations Research}
  34(2):303--319.

\bibitem[{Bartl et~al.(2017)Bartl, Kupper, Lux, Papapantoleon,
  \protect\BIBand{} Eckstein}]{bartl2017marginal}
Bartl D, Kupper M, Lux T, Papapantoleon A, Eckstein S (2017) Marginal and
  dependence uncertainty: bounds, optimal transport, and sharpness. Preprint,
  arXiv:1709.00641.

\bibitem[{Bartl et~al.(2020)Bartl, Kupper, \protect\BIBand{}
  Neufeld}]{bartl2020pathwise}
Bartl D, Kupper M, Neufeld A (2020) Pathwise superhedging on prediction sets.
  \emph{Finance and Stochastics} 24(1):215--248.

\bibitem[{Bartl et~al.(2019)Bartl, Kupper, Pr{\"o}mel, \protect\BIBand{}
  Tangpi}]{bartl2019duality}
Bartl D, Kupper M, Pr{\"o}mel DJ, Tangpi L (2019) Duality for pathwise
  superhedging in continuous time. \emph{Finance and Stochastics}
  23(3):697--728.

\bibitem[{Bayraktar et~al.(2015)Bayraktar, Huang, \protect\BIBand{}
  Zhou}]{bayraktar2015hedging}
Bayraktar E, Huang YJ, Zhou Z (2015) On hedging {A}merican options under model
  uncertainty. \emph{SIAM Journal on Financial Mathematics} 6(1):425--447.

\bibitem[{Bayraktar \protect\BIBand{} Zhang(2016)}]{Bayrakter_Zhang_2016}
Bayraktar E, Zhang Y (2016) Fundamental theorem of asset pricing under
  transaction costs and model uncertainty. \emph{Mathematics of Operations
  Research} 41:1039--1054.

\bibitem[{Bayraktar et~al.(2014)Bayraktar, Zhang, \protect\BIBand{}
  Zhou}]{Bayrakter_Zhang_Zhou_2014}
Bayraktar E, Zhang Y, Zhou Z (2014) A note on the fundamental theorem of asset
  pricing under model uncertainty. \emph{Risks} 2:425--433.

\bibitem[{Beiglb\"{o}ck et~al.(2013)Beiglb\"{o}ck, Henry-Labord{\`e}re,
  \protect\BIBand{} Penkner}]{Beiglboeck_HenryLabordere_Penkner_2013}
Beiglb\"{o}ck M, Henry-Labord{\`e}re P, Penkner F (2013) Model-independent
  bounds for option prices---a mass transport approach. \emph{Finance and
  Stochastics} 17:477--501.

\bibitem[{Beissner(2017)}]{beissner2017equilibrium}
Beissner P (2017) Equilibrium prices and trade under ambiguous volatility.
  \emph{Economic Theory} 64(2):213--238.

\bibitem[{Beissner \protect\BIBand{} Riedel(2019)}]{beissner2019equilibria}
Beissner P, Riedel F (2019) Equilibria under {K}nightian price uncertainty.
  \emph{Econometrica} 87(1):37--64.

\bibitem[{Bertsimas \protect\BIBand{} Bushueva(2006)}]{bertsimas2006option}
Bertsimas D, Bushueva N (2006) Option pricing without price dynamics: a
  probabilistic approach. Preprint, arXiv:math/0612075.

\bibitem[{Bertsimas \protect\BIBand{} Popescu(2002)}]{bertsimas2002on}
Bertsimas D, Popescu I (2002) On the relation between option and stock prices:
  a convex optimization approach. \emph{Operations Research} 50(2):358--374.

\bibitem[{Betr{\`o}(2004)}]{betro2004accelerated}
Betr{\`o} B (2004) An accelerated central cutting plane algorithm for linear
  semi-infinite programming. \emph{Mathematical Programming} 101(3):479--495.

\bibitem[{Biagini et~al.(2017)Biagini, Bouchard, Kardaras, \protect\BIBand{}
  Nutz}]{Biagini_Bouchard_Kardaras_Nutz_2017}
Biagini S, Bouchard B, Kardaras C, Nutz M (2017) Robust fundamental theorem for
  continuous processes. \emph{Mathematical Finance} 27:963--987.

\bibitem[{Borwein \protect\BIBand{}
  Lewis(1991{\natexlab{a}})}]{borwein1991convergence}
Borwein JM, Lewis AS (1991{\natexlab{a}}) Convergence of best entropy
  estimates. \emph{SIAM Journal on Optimization} 1(2):191--205.

\bibitem[{Borwein \protect\BIBand{}
  Lewis(1991{\natexlab{b}})}]{borwein1991duality}
Borwein JM, Lewis AS (1991{\natexlab{b}}) Duality relationships for
  entropy-like minimization problems. \emph{SIAM Journal on Control and
  Optimization} 29(2):325--338.

\bibitem[{Bouchard et~al.(2019)Bouchard, Deng, \protect\BIBand{}
  Tan}]{bouchard2019superreplication}
Bouchard B, Deng S, Tan X (2019) Superreplication with proportional transaction
  cost under model uncertainty. \emph{Mathematical Finance} 29(3):837--860.

\bibitem[{Bouchard \protect\BIBand{} Nutz(2015)}]{bouchard2015arbitrage}
Bouchard B, Nutz M (2015) Arbitrage and duality in nondominated discrete-time
  models. \emph{The Annals of Applied Probability} 25(2):823--859.

\bibitem[{Burzoni et~al.(2019)Burzoni, Frittelli, Hou, Maggis,
  \protect\BIBand{} Ob{\l}{\'o}j}]{burzoni2019pointwise}
Burzoni M, Frittelli M, Hou Z, Maggis M, Ob{\l}{\'o}j J (2019) Pointwise
  arbitrage pricing theory in discrete time. \emph{Mathematics of Operations
  Research} 44(3):1034--1057.

\bibitem[{Burzoni et~al.(2016)Burzoni, Frittelli, \protect\BIBand{}
  Maggis}]{burzoni2016universal}
Burzoni M, Frittelli M, Maggis M (2016) Universal arbitrage aggregator in
  discrete-time markets under uncertainty. \emph{Finance and Stochastics}
  20(1):1--50.

\bibitem[{Burzoni et~al.(2021)Burzoni, Riedel, \protect\BIBand{}
  Soner}]{burzoni2021viability}
Burzoni M, Riedel F, Soner HM (2021) Viability and arbitrage under {K}nightian
  uncertainty. \emph{Econometrica} 89(3):1207--1234.

\bibitem[{Chen et~al.(2008)Chen, Deelstra, Dhaene, \protect\BIBand{}
  Vanmaele}]{Chen_Deelstra_Dhaene_Vanmaele_2008}
Chen X, Deelstra G, Dhaene J, Vanmaele M (2008) Static super-replicating
  strategies for a class of exotic options. \emph{Insurance Math. Econom.}
  42:1067--1085.

\bibitem[{Cheridito et~al.(2017)Cheridito, Kupper, \protect\BIBand{}
  Tangpi}]{cheridito2017duality}
Cheridito P, Kupper M, Tangpi L (2017) Duality formulas for robust pricing and
  hedging in discrete time. \emph{SIAM Journal on Financial Mathematics}
  8(1):738--765.

\bibitem[{Cho et~al.(2016)Cho, Kim, \protect\BIBand{} Lee}]{Cho_Kim_Lee_2016}
Cho H, Kim KK, Lee K (2016) Computing lower bounds on basket option prices by
  discretizing semi-infinite linear programming. \emph{Optimization Letters}
  10:1629--1644.

\bibitem[{Cohen et~al.(2020)Cohen, Reisinger, \protect\BIBand{}
  Wang}]{cohen2020detecting}
Cohen SN, Reisinger C, Wang S (2020) Detecting and repairing arbitrage in
  traded option prices. Preprint, arXiv:2008.09454.

\bibitem[{Coope \protect\BIBand{} Price(1998)}]{coope1998exact}
Coope ID, Price CJ (1998) Exact penalty function methods for nonlinear
  semi-infinite programming. \emph{Semi-infinite programming}, volume~25 of
  \emph{Nonconvex Optim. Appl.}, 137--157 (Kluwer Acad. Publ., Boston, MA).

\bibitem[{Dana \protect\BIBand{} Riedel(2013)}]{dana2013intertemporal}
Dana RA, Riedel F (2013) Intertemporal equilibria with {K}nightian uncertainty.
  \emph{Journal of Economic Theory} 148(4):1582--1605.

\bibitem[{d'Aspremont \protect\BIBand{}
  El~Ghaoui(2006)}]{dAspremont_ElGhaoui_2006}
d'Aspremont A, El~Ghaoui L (2006) Static arbitrage bounds on basket option
  prices. \emph{Math. Program.} 106(3, Ser. A):467--489.

\bibitem[{Daum \protect\BIBand{} Werner(2011)}]{daum2011novel}
Daum S, Werner R (2011) A novel feasible discretization method for linear
  semi-infinite programming applied to basket option pricing.
  \emph{Optimization} 60(10-11):1379--1398.

\bibitem[{Davis et~al.(2014)Davis, Ob{\l}{\'o}j, \protect\BIBand{}
  Raval}]{davis2014arbitrage}
Davis M, Ob{\l}{\'o}j J, Raval V (2014) Arbitrage bounds for prices of weighted
  variance swaps. \emph{Mathematical Finance} 24(4):821--854.

\bibitem[{De~Gennaro~Aquino \protect\BIBand{} Bernard(2020)}]{aquino2019bounds}
De~Gennaro~Aquino L, Bernard C (2020) Bounds on multi-asset derivatives via
  neural networks. \emph{Int. J. Theor. Appl. Finance} 23(8):2050050, 31.

\bibitem[{Dhaene et~al.(2002{\natexlab{a}})Dhaene, Denuit, Goovaerts, Kaas,
  \protect\BIBand{} Vyncke}]{Dhaene_etal_2002_b}
Dhaene J, Denuit M, Goovaerts MJ, Kaas R, Vyncke D (2002{\natexlab{a}}) The
  concept of comonotonicity in actuarial science and finance: applications.
  \emph{Insurance Math. Econom.} 31:133--161.

\bibitem[{Dhaene et~al.(2002{\natexlab{b}})Dhaene, Denuit, Goovaerts, Kaas,
  \protect\BIBand{} Vyncke}]{Dhaene_etal_2002_a}
Dhaene J, Denuit M, Goovaerts MJ, Kaas R, Vyncke D (2002{\natexlab{b}}) The
  concept of comonotonicity in actuarial science and finance: theory.
  \emph{Insurance Math. Econom.} 31:3--33.

\bibitem[{Dolinsky \protect\BIBand{} Neufeld(2018)}]{dolinsky2018super}
Dolinsky Y, Neufeld A (2018) Super-replication in fully incomplete markets.
  \emph{Mathematical Finance} 28(2):483--515.

\bibitem[{Dolinsky \protect\BIBand{}
  Soner(2014{\natexlab{a}})}]{dolinsky2014martingale}
Dolinsky Y, Soner HM (2014{\natexlab{a}}) Martingale optimal transport and
  robust hedging in continuous time. \emph{Probability Theory and Related
  Fields} 160(1-2):391--427.

\bibitem[{Dolinsky \protect\BIBand{}
  Soner(2014{\natexlab{b}})}]{Dolinsky_Soner_2014_b}
Dolinsky Y, Soner HM (2014{\natexlab{b}}) Robust hedging with proportional
  transaction costs. \emph{Finance and Stochastics} 18:327--347.

\bibitem[{Eckstein et~al.(2021)Eckstein, Guo, Lim, \protect\BIBand{}
  Ob\l\'{o}j}]{eckstein2019robust}
Eckstein S, Guo G, Lim T, Ob\l\'{o}j J (2021) Robust pricing and hedging of
  options on multiple assets and its numerics. \emph{SIAM J. Financial Math.}
  12(1):158--188.

\bibitem[{Eckstein \protect\BIBand{} Kupper(2019)}]{eckstein2019computation}
Eckstein S, Kupper M (2019) Computation of optimal transport and related
  hedging problems via penalization and neural networks. \emph{Applied
  Mathematics \& Optimization} 1--29.

\bibitem[{Eckstein et~al.(2020)Eckstein, Kupper, \protect\BIBand{}
  Pohl}]{eckstein2018robust}
Eckstein S, Kupper M, Pohl M (2020) Robust risk aggregation with neural
  networks. \emph{Mathematical Finance} 30:1229--1272.

\bibitem[{Epstein \protect\BIBand{} Ji(2013)}]{epstein2013ambiguous}
Epstein LG, Ji S (2013) Ambiguous volatility and asset pricing in continuous
  time. \emph{The Review of Financial Studies} 26(7):1740--1786.

\bibitem[{Galichon et~al.(2014)Galichon, Henry-Labord{\`e}re, \protect\BIBand{}
  Touzi}]{galichon2014stochastic}
Galichon A, Henry-Labord{\`e}re P, Touzi N (2014) A stochastic control approach
  to no-arbitrage bounds given marginals, with an application to lookback
  options. \emph{The Annals of Applied Probability} 24(1):312--336.

\bibitem[{Goberna \protect\BIBand{} L{\'o}pez(1998)}]{goberna1998linear}
Goberna MA, L{\'o}pez MA (1998) \emph{Linear Semi-Infinite Optimization} (John
  Wiley \& Sons).

\bibitem[{Goberna \protect\BIBand{} L{\'o}pez(2018)}]{goberna2018recent}
Goberna MA, L{\'o}pez MA (2018) Recent contributions to linear semi-infinite
  optimization: an update. \emph{Annals of Operations Research}
  271(1):237--278.

\bibitem[{{Gurobi Optimization, LLC}(2020)}]{gurobi}
{Gurobi Optimization, LLC} (2020) Gurobi optimizer reference manual.
  \urlprefix\url{http://www.gurobi.com}.

\bibitem[{Henry-Labord{\`e}re(2013)}]{henry2013automated}
Henry-Labord{\`e}re P (2013) Automated option pricing: {N}umerical methods.
  \emph{International Journal of Theoretical and Applied Finance}
  16(08):1350042.

\bibitem[{Hettich \protect\BIBand{} Kortanek(1993)}]{Hettich_Kortanek_1993}
Hettich R, Kortanek KO (1993) Semi-infinite programming: theory, methods, and
  applications. \emph{SIAM Review} 35:380--429.

\bibitem[{Hobson et~al.(2005{\natexlab{a}})Hobson, Laurence, \protect\BIBand{}
  Wang}]{Hobson_Laurence_Wang_2005_2}
Hobson D, Laurence P, Wang TH (2005{\natexlab{a}}) Static-arbitrage optimal
  subreplicating strategies for basket options. \emph{Insurance Math. Econom.}
  37:553--572.

\bibitem[{Hobson et~al.(2005{\natexlab{b}})Hobson, Laurence, \protect\BIBand{}
  Wang}]{Hobson_Laurence_Wang_2005_1}
Hobson D, Laurence P, Wang TH (2005{\natexlab{b}}) Static-arbitrage upper
  bounds for the prices of basket options. \emph{Quant. Finance} 5:329--342.

\bibitem[{Hobson(1998)}]{hobson1998robust}
Hobson DG (1998) Robust hedging of the lookback option. \emph{Finance and
  Stochastics} 2(4):329--347.

\bibitem[{Hu et~al.(2019)Hu, Li, \protect\BIBand{} Mehrotra}]{hu2019data}
Hu J, Li J, Mehrotra S (2019) A data-driven functionally robust approach for
  simultaneous pricing and order quantity decisions with unknown demand
  function. \emph{Operations Research} 67(6):1564--1585.

\bibitem[{Kahal{\'e}(2017)}]{kahale2017superreplication}
Kahal{\'e} N (2017) Superreplication of financial derivatives via convex
  programming. \emph{Management Science} 63(7):2323--2339.

\bibitem[{Karlin \protect\BIBand{} Studden(1966)}]{Karlin_Studden_1966}
Karlin S, Studden WJ (1966) \emph{Tchebycheff systems: {W}ith applications in
  analysis and statistics} (John Wiley \& Sons).

\bibitem[{Lin et~al.(1998)Lin, Fang, \protect\BIBand{}
  Wu}]{lin1998unconstrained}
Lin CJ, Fang SC, Wu SY (1998) An unconstrained convex programming approach to
  linear semi-infinite programming. \emph{SIAM Journal on Optimization}
  8(2):443--456.

\bibitem[{L\'{o}pez \protect\BIBand{} Still(2007)}]{lopez2007semi}
L\'{o}pez M, Still G (2007) Semi-infinite programming. \emph{European J. Oper.
  Res.} 180(2):491--518.

\bibitem[{L{\"u}tkebohmert \protect\BIBand{}
  Sester(2019)}]{lutkebohmert2019tightening}
L{\"u}tkebohmert E, Sester J (2019) Tightening robust price bounds for exotic
  derivatives. \emph{Quantitative Finance} 19(11):1797--1815.

\bibitem[{Lux \protect\BIBand{} Papapantoleon(2017)}]{lux2016}
Lux T, Papapantoleon A (2017) Improved {F}r\'echet--{H}oeffding bounds on
  $d$-copulas and applications in model-free finance. \emph{Ann. Appl. Probab.}
  27:3633--3671.

\bibitem[{Maruhn(2009)}]{maruhn2009robust}
Maruhn JH (2009) \emph{Robust Static Super-Replication of Barrier Options}
  (Walter de Gruyter GmbH \& Co. KG, Berlin).

\bibitem[{Neufeld \protect\BIBand{} Nutz(2013)}]{neufeld2013superreplication}
Neufeld A, Nutz M (2013) Superreplication under volatility uncertainty for
  measurable claims. \emph{Electronic Journal of Probability} 18.

\bibitem[{Nishihara et~al.(2007)Nishihara, Yagiura, \protect\BIBand{}
  Ibaraki}]{nishihara2007duality}
Nishihara M, Yagiura M, Ibaraki T (2007) Duality in option pricing based on
  prices of other derivatives. \emph{Oper. Res. Lett.} 35(2):165--171.

\bibitem[{Papapantoleon \protect\BIBand{}
  Yanez~Sarmiento(2020)}]{papapantoleon2020detection}
Papapantoleon A, Yanez~Sarmiento P (2020) Detection of arbitrage opportunities
  in multi-asset derivatives markets. Preprint, arXiv:2002.06227.

\bibitem[{Pe{\~n}a et~al.(2010{\natexlab{a}})Pe{\~n}a, Saynac, Vera,
  \protect\BIBand{} Zuluaga}]{pena2010computing}
Pe{\~n}a J, Saynac X, Vera JC, Zuluaga LF (2010{\natexlab{a}}) Computing
  general static-arbitrage bounds for {European} basket options via
  {Dantzig--Wolfe} decomposition. \emph{Algorithmic Operations Research}
  5(2):65--74.

\bibitem[{Pe{\~n}a et~al.(2010{\natexlab{b}})Pe{\~n}a, Vera, \protect\BIBand{}
  Zuluaga}]{Pena_Vera_Zuluaga_2010}
Pe{\~n}a J, Vera JC, Zuluaga LF (2010{\natexlab{b}}) Static-arbitrage lower
  bounds on the prices of basket options via linear programming. \emph{Quant.
  Finance} 10:819--827.

\bibitem[{Pe{\~n}a et~al.(2012)Pe{\~n}a, Vera, \protect\BIBand{}
  Zuluaga}]{pena2012computing}
Pe{\~n}a J, Vera JC, Zuluaga LF (2012) Computing arbitrage upper bounds on
  basket options in the presence of bid--ask spreads. \emph{European Journal of
  Operational Research} 222(2):369--376.

\bibitem[{Possama{\"\i} et~al.(2013)Possama{\"\i}, Royer, \protect\BIBand{}
  Touzi}]{possamai2013robust}
Possama{\"\i} D, Royer G, Touzi N (2013) On the robust superhedging of
  measurable claims. \emph{Electronic Communications in Probability} 18.

\bibitem[{Puccetti et~al.(2016)Puccetti, R{\"u}schendorf, \protect\BIBand{}
  Manko}]{puccetti2016}
Puccetti G, R{\"u}schendorf L, Manko D (2016) {VaR} bounds for joint portfolios
  with dependence constraints. \emph{Depend. Model.} 4:368--381.

\bibitem[{Reemtsen \protect\BIBand{}
  R\"uckmann(1998)}]{Reemtsen_Rueckmann_1998}
Reemtsen R, R\"uckmann JJ (1998) \emph{Semi-infinite programming}, volume~25
  (Springer Science \& Business Media).

\bibitem[{Riedel(2015)}]{riedel2015financial}
Riedel F (2015) Financial economics without probabilistic prior assumptions.
  \emph{Decisions in Economics and Finance} 38(1):75--91.

\bibitem[{Rigotti \protect\BIBand{} Shannon(2005)}]{rigotti2005uncertainty}
Rigotti L, Shannon C (2005) Uncertainty and risk in financial markets.
  \emph{Econometrica} 73(1):203--243.

\bibitem[{Rockafellar(1970)}]{rockafellar1970convex}
Rockafellar RT (1970) \emph{Convex Analysis} (Princeton University Press).

\bibitem[{Stein(2012)}]{stein2012how}
Stein O (2012) How to solve a semi-infinite optimization problem.
  \emph{European J. Oper. Res.} 223(2):312--320.

\bibitem[{Tankov(2011)}]{tankov}
Tankov P (2011) Improved {Fr\'echet} bounds and model-free pricing of
  multi-asset options. \emph{J. Appl. Probab.} 48:389--403.

\bibitem[{Tavin(2015)}]{Tavin_2015}
Tavin B (2015) Detection of arbitrage in a market with multi-asset derivatives
  and known risk-neutral marginals. \emph{J. Banking Finance} 53:158--178.

\bibitem[{Yan et~al.(2018)Yan, Cheng, Natarajan, \protect\BIBand{}
  Teo}]{yan2018marginal}
Yan Z, Cheng C, Natarajan K, Teo CP (2018) {``M}arginal estimation+ price
  optimization'' for multi-product pricing problems. Technical report, Working
  paper.

\end{thebibliography}

\begin{thebibliography}{13}
\providecommand{\natexlab}[1]{#1}
\providecommand{\url}[1]{\texttt{#1}}
\providecommand{\urlprefix}{URL }

\bibitem[{mainpaperref(2000)}]{mainpaperrefs}
{\fs.{10}.{12}.{\it See references list in the main paper.}}

\bibitem[{Bartl et~al.(2017)Bartl, Kupper, Lux, Papapantoleon,
  \protect\BIBand{} Eckstein}]{bartl2017marginal}
Bartl D, Kupper M, Lux T, Papapantoleon A, Eckstein S (2017) Marginal and
  dependence uncertainty: bounds, optimal transport, and sharpness. Preprint,
  arXiv:1709.00641.

\bibitem[{Betr{\`o}(2004)}]{betro2004accelerated}
Betr{\`o} B (2004) An accelerated central cutting plane algorithm for linear
  semi-infinite programming. \emph{Mathematical Programming} 101(3):479--495.

\bibitem[{Borwein \protect\BIBand{} Lewis(1992)}]{borwein1992partially}
Borwein JM, Lewis AS (1992) Partially finite convex programming, part {I}:
  Quasi relative interiors and duality theory. \emph{Mathematical Programming}
  57(1-3):15--48.

\bibitem[{Bouchard \protect\BIBand{} Nutz(2015)}]{bouchard2015arbitrage}
Bouchard B, Nutz M (2015) Arbitrage and duality in nondominated discrete-time
  models. \emph{The Annals of Applied Probability} 25(2):823--859.

\bibitem[{Ciripoi et~al.(2019)Ciripoi, L{\"o}hne, \protect\BIBand{}
  Wei{\ss}ing}]{ciripoibensolve}
Ciripoi D, L{\"o}hne A, Wei{\ss}ing B (2019) Bensolve tools, version 1.3. {Gnu
  Octave/Matlab} toolbox for calculus of convex polyhedra, calculus of
  polyhedral convex functions, global optimization, vector linear programming.

\bibitem[{Dellacherie \protect\BIBand{}
  Meyer(1982)}]{dellacherie1982probabilities}
Dellacherie C, Meyer PA (1982) Probabilities and potential. b.
  \emph{North-Holland Mathematics Studies}.

\bibitem[{Goberna \protect\BIBand{} L{\'o}pez(1998)}]{goberna1998linear}
Goberna MA, L{\'o}pez MA (1998) \emph{Linear Semi-Infinite Optimization} (John
  Wiley \& Sons).

\bibitem[{Kellerer(1984)}]{kellerer1984duality}
Kellerer HG (1984) Duality theorems for marginal problems. \emph{Zeitschrift
  f{\"u}r Wahrscheinlichkeitstheorie und verwandte Gebiete} 67(4):399--432.

\bibitem[{Le~Verge(1992)}]{le1992note}
Le~Verge H (1992) A note on {Chernikova's} algorithm.

\bibitem[{Rachev \protect\BIBand{}
  R\"{u}schendorf(1994)}]{Rachev_Rueschendorf_1994}
Rachev ST, R\"{u}schendorf L (1994) Solution of some transportation problems
  with relaxed or additional constraints. \emph{SIAM Journal on Control and
  Optimization} 32:673--689.

\bibitem[{Rockafellar(1970)}]{rockafellar1970convex}
Rockafellar RT (1970) \emph{Convex Analysis} (Princeton University Press).

\bibitem[{Zaev(2015)}]{Zaev_2015}
Zaev DA (2015) On the {M}onge-{K}antorovich problem with additional linear
  constraints. \emph{Matematicheskie Zametki} 98:664--683.

\bibitem[{Zhen \protect\BIBand{} Den~Hertog(2018)}]{zhen2018computing}
Zhen J, Den~Hertog D (2018) Computing the maximum volume inscribed ellipsoid of
  a polytopic projection. \emph{INFORMS Journal on Computing} 30(1):31--42.

\end{thebibliography}
